\magnification=1200
\baselineskip=12pt
\input amssym
\input xypic
\font\bigly=cmbx14
\font\bagly=cmbx12
\overfullrule=0pt

\def\ARI{{\rm ARI}}
\def\B\ARI{{\rm B\ARI}}
\def\G\ARI{{\rm G\ARI}}
\def\GAXI{{\rm GAXI}}
\def\GB\ARI{{\rm GB\ARI}}
\def\GV\ARI{{\overline{\G\ARI}}}
\def\sec{{\rm sec}}
\def\ret{{\rm ret}}
\def\ds{{\frak{ds}}}
\def\mt{{\frak{mt}}}
\def\ls{{\frak{ls}}}
\def\nfz{{\frak{nfz}}}

\def\x{{\bf x}}

\def\u{{\bf u}}
\def\v{{\bf v}}
\def\w{{\bf w}}
\def\x{{\bf x}}
\def\a{{\bf a}}
\def\b{{\bf b}}
\def\c{{\bf c}}
\def\d{{\bf d}}

\def\Q{{\Bbb Q}}
\def\F{{\Bbb Q\langle C\rangle}}
\def\Fo{{\Bbb Q\langle\langle C\rangle\rangle}}
\def\Fzero{{\Bbb Q_0\langle C\rangle}}
\def\Fzeroo{{\Bbb Q_0\langle\langle C\rangle\rangle}}

\ \vskip 2cm
\centerline{\bigly ARI, GARI, ZIG and ZAG}
\vskip .5cm
\centerline{\bagly An introduction to Ecalle's theory of multiple zeta values}
\vskip  1cm
\centerline{\bf Leila Schneps}
\vskip .5cm
\centerline{\it with contributions by}
\vskip .5cm
\centerline{\bf Samuel Baumard, Nao Komiyama, Adriana Salerno}
\vskip  1cm
{\narrower{\narrower{The text has two goals. The first is to give an 
introduction to Ecalle's work on mould theory, multiple zeta values and double 
shuffle theory and relate this work explicitly to the classical theory of 
multiple zeta values and double shuffle expressed in the usual terms of 
two non-commutative variables. The second is to provide complete proofs of 
those of his main statements and identities which are useful in the
context of (non-colored) multiple zeta values.  Many of these proofs were
never written down by \'Ecalle.  Some of them are difficult, laborious and 
not enlightening, yet it is clearly necessary to have them in order to be 
able to apply with confidence a theory that, once in place, forms an 
astonishingly powerful toolbox with many applications.  Of these laborious 
proofs, some have been relegated to appendices and others, which
appear in full in separate publications, have simply been cited.

The emphasis in this text is to provide an easily approachable 
introduction to Ecalle's language while placing it almost from the start 
in the context of multiple zeta value theory.
\vskip .5cm
\noindent {\bf Disclaimer:} This text is not final and is not submitted
for publication.  The intention is to continue to add to and complete it
over time.}\par}\par}
\vfill\eject
\ \vskip 1cm
\centerline{\bf Contents}
\ \vskip .5cm
\noindent {\bf Chapter 1: Real and formal multiple zeta values}\hfill 3
\ \vskip .5cm
\S 1.1. Multiple zeta values and their regularizations

\S 1.2. Formal multiple zeta values

\S 1.3. The double shuffle Lie algebra $\ds$

\S 1.4. The linearized double shuffle space
\ \vskip .5cm
\noindent {\bf Chapter 2: Introduction to \ARI and its panoply of tools}\hfill 14
\ \vskip .5cm
\S 2.1. Moulds and bimoulds

\S 2.2. The Lie algebra $\ARI$

\S 2.3. Symmetrality, alternality, symmetrility, alternility

\S 2.4. $Swap$ commutation in $\ARI$

\S 2.5. Special subspaces of $\ARI$

\S 2.6. Circ-neutrality

\S 2.7. The group $\G\ARI$

\S 2.8. The group law on $\G\ARI$

\S 2.9. \'Ecalle's first fundamental identity: $swap$ commutation in 
$\G\ARI$
\ \vskip .5cm
\noindent {\bf Chapter 3: From double shuffle to \ARI}\hfill 38
\ \vskip .5cm
\S 3.1. The ring $\F$

\S 3.2. Associating moulds to elements $f\in \F$

\S 3.3. The Poisson bracket and the \ARI bracket

\S 3.4. The $ma$ map from $\ds$ to $\ARI$

\S 3.5. The group $\G\ARI$ and the twisted Magnus group
\ \vskip .5cm
\noindent {\bf Chapter 4: The mould pair $pal/pil$ and its symmetries}\hfill 49
\ \vskip .5cm
\S 4.1. Diffeomorphisms and the mould $pil$

\S 4.2. Two definitions of the mould $pal$

\S 4.3. Symmetrality of $pal$

\S 4.4. The identity $crash(pal)=pac$

\S 4.5. Ecalle's second fundamental identity

\S 4.6. Double shuffle is a Lie algebra

\S 4.7. The $\Delta$-denominator
\ \vskip .5cm
\noindent {\bf Chapter 5: Elliptic mould theory}\hfill 69
\ \vskip .5cm
\S 5.1. The operator $Ad_{ari}(invpal)$ and the denominator $\Delta$

\S 5.2. $\Delta$ as a Lie algebra isomorphism, and the
$Dari$-bracket

\S 5.3. Adding the mould $a$ to $\ARI$

\S 5.4. Closed subspaces of $\ARI_{Dari}$

\S 5.5. The real function of the moulds $pal$ and $invpal$

\S 5.6. The real meaning of the operator $\Delta\circ Ad_{ari}(invpal)$
\ \vskip .5cm
\noindent {\bf Appendix}\hfill 86
\ \vskip .5cm
\S A.1. Proof of Proposition 2.2.1

\S A.2. Proofs of (2.4.7) and (2.4.8)

\S A.3. Proof of Lemma 3.2.1

\S A.4. Proof of Proposition 3.3.2

\S A.5. Proof of Lemma 3.4.1

\S A.6. Proof of Proposition 4.2.6
\ \vfill\eject
\ \vskip .5cm
\centerline{\bf Chapter 1} 
\vskip .5cm
\centerline{\bf Real and formal multiple zeta values}
\vskip .8cm
In this first chapter, we introduce some of the basic objects of study in the 
classical theory; the algebras of real and formal multiple zeta values, 
the real and formal Drinfel'd associators, the double shuffle Lie algebra, and 
the weight grading and depth filtrations. Everything in this chapter is well-known
and has been written in detail elsewhere, so we content ourselves with recalling
the main definitions and facts without proof.
\vskip .8cm
\noindent {\bf \S 1.1. Multiple zeta values and their regularizations}
\vskip .5cm
For every sequence ${\bf k}=(k_1,\ldots,k_r)$ of strictly positive integers
with $k_1\ge 2$, let $\zeta(k_1,\ldots,k_r)$ be the {\it multiple zeta
value} defined by
$$\zeta(k_1,\ldots,k_r)=\sum_{n_1>\cdots>n_r>0} {{1}\over{n_1^{k_1}
\cdots n_r^{k_r}}}.\eqno(1.1.1)$$
For every word in $\Q\langle x,y\rangle$, we define a multiple zeta value
$\zeta(w)$ as follows.  If $w$ starts with $x$ and ends with $y$,
we write $w=x^{k_1-1}y\cdots x^{k_r-1}y$ with $k_1\ge 2$, and set
$\zeta(w)=\zeta(k_1,\ldots,k_r)$.

For general $w$, we write $w=y^rvx^s$ and set
$$\zeta(w)=\sum_{a=0}^r \sum_{b=0}^s (-1)^{a+b}\zeta\bigl(\pi(sh(y^a,y^{r-a}v
x^{s-b},x^b))\bigr),\eqno(1.1.2)$$
where $\pi$ is the projection of a polynomial onto the {\it convergent} words,
i.e. those starting with $x$ and ending with $y$, and $\zeta$ is considered
to be additive.  This way of extending the real multizeta values of convergent
words (called {\it convergent multizeta values}) to all words is called the 
{\it shuffle regularization}, because of the following property that 
characterizes it.
\vskip .2cm
\noindent {\bf Definition.} The {\it shuffle product} of two words $u$ and $v$ in an alphabet ${\cal X}$ is defined recursively by 
$sh(u,1)=sh(1,u)=u$ and $sh(Xu,Yv)=X\,sh(u,Yv)+Y\,sh(Xu,v)$ for any letters $X,Y\in {\cal X}$.
\vskip .3cm
The path leading to the formula given in (1.1.2) is not a short one, starting
as it does by using standard regularization techniques to give regularized
values to the non-convergent multizeta values in the form of integrals over
simplices ([LM]). The explicit formula (1.1.2) was established by H.~Furusho
in [F] (Prop. 3.2.3).
\vskip .3cm
\noindent {\bf Examples.} We use the notation in which
the shuffle of two words is written as a formal sum of words.
Taking ${\cal X}=\{a,b,c,d\}$, we have 
$$sh((ab),(cd))=abcd+acbd+acdb+cabd+cadb+cdab.$$
Taking ${\cal X}=\{x,y\}$, we thus have
$$sh((x,y),(x,y))=4xxyy+2xyxy.$$
\vskip .3cm
\noindent {\bf Theorem 1.1.1.} {\it For all words $u$, $v\in \Q\langle x,y
\rangle$, the regularized $\zeta$ values defined in (1.1.2) satisfy the 
shuffle relations
$$\zeta\bigl(sh(u,v)\bigr)=\zeta(u)\zeta(v)\eqno(1.1.3)$$
in the alphabet ${\cal X}=\{x,y\}$.}
\vskip .3cm
Multiple zeta values possess a second interesting multiplicative property.  
\vskip .2cm
\noindent {\bf Definition.} Let ${\cal Y}$ be an additive alphabet, i.e. a set equipped with an addition rule 
such that for every pair of letters $X,Y\in {\cal Y}$, $X+Y$ is also an element of
${\cal Y}$. The stuffle product in the additive alphabet
${\cal Y}$ is defined recursively by $st(u,1)=st(1,u)=u$ and
$$st(Xu,Yv)=X\,st(u,Yv)+Y\,st(Xu,v)+ (X+Y)\,st(u,v)\eqno(1.1.4)$$ 
for all letters $X,Y\in {\cal Y}$.

An equivalent formulation of the stuffle product is given by
$$st(u,v)=\sum_{\sigma\in Sh^{\le}(r,s)} c^\sigma(u,v)\eqno(1.1.5)$$
where $u$ is a word in $r$ letters and $v$ in $s$ letters, 
$Sh^{\le}(r,s)$ is the set of surjective maps
$$\sigma:\{1,\ldots,r+s\}\rightarrow\!\!\!\rightarrow \{1,\ldots,N\}$$
for all $1\le N\le r+s$ such that
$$\sigma(1)<\cdots<\sigma(r)\ \ {\rm and}\ \ \sigma(r+1)<\ldots<\sigma(r+s),$$
and for each $\sigma\in Sh^{\le}(r,s)$, we set
$c^\sigma(u,v)=(c_1,\ldots,c_N)$ with
$$c_i=\sum_{k\in \sigma^{-1}(i)} a_k.\eqno(1.1.6)$$
By the definition of $Sh^{\le}(r,s)$, $c_i$ is either a single letter $a_k$ or 
a sum of two letters $a_k+a_l$ with $k\le r<l$.
\vskip .3cm
\noindent {\bf Examples.} Let ${\cal A}$ be an additive alphabet; then we have
$$\eqalign{st(a,b)&=(a,b)+(b,a)+(a+b)\cr
st((a,b),(c))&=abc+acb+cab+(a+b,c)+(a,b+c)\cr
st((a,b),(b))&=2(a,b,b)+(b,a,b)+(a+b,b)+(a,2b).}$$
Considering the additive alphabet ${\Bbb N}^+$, we have for example
$$st((2,1),(2))=2(2,2,1)+(2,1,2)+(4,1)+(2,3).$$
In a different notation that will be used often below, let ${\cal Y}=\{y_1,y_2,y_3,
\ldots\}$ with the addition rule $y_i+y_j=y_{i+j}$. This is identical to considering
the alphabet ${\Bbb N}^+$ except that the numbers now appear as indices.  We have
for example
$$st((y_1),(y_2,y_3))=(y_1,y_2,y_3)+(y_2,y_1,y_3)+(y_2,y_3,y_1)+(y_3,y_3)+
(y_2,y_4).\eqno(1.1.7)$$
For all convergent words $u$, $v$, considered to be written in the
variables $y_i=x^{i-1}y$, the convergent multizeta values satisfy the {\it stuffle
relations} $\zeta\bigl(st(u,v)\bigr)=\zeta(u)\zeta(v)$ in the alphabet ${\cal Y}=\{y_i|i\ge 0\}$,
considered to be additive via the rule $y_i+y_j=y_{i+j}$.  This result
follows easily from the expression of $\zeta(k_1,\ldots,k_r)$ as a power
series.  But there is a second regularization of the zeta values, called
the {\it stuffle regularization}, extending the stuffle relation to all
words in the $y_i$.  It is defined as follows.
\vskip .2cm
\noindent {\bf Definition.} The {\it Drinfel'd associator} $\Phi_{KZ}$ is 
given by 
$$\Phi_{KZ}=1+\sum_{w\in \Q\langle x,y\rangle} (-1)^{d(w)}\zeta(w)w,\eqno(1.1.8)$$
where for each monomial $w$ in $x,y$, $d(w)$ denotes the {\it depth} of $w$, which is
the number of $y$'s occurring in the word $w$.
Let $\Phi$ denote the {\it double shuffle power series} defined by
$\Phi(x,y)=\Phi_{KZ}(x,-y)$, so
$$\Phi(x,y)=1+\sum_w \zeta(w)w.$$
Let $\pi_y$ denote the projection of power series onto their words ending
in $y$, rewritten in the $y_i$.  Set
$$\Phi_*=exp\Bigl(\sum_{n\ge 1} {{(-1)^{n-1}}\over{n}}\zeta(y_n)y_1^n\Bigr)
\pi_y(\Phi),\eqno(1.1.9)$$
and for every word $v$ in the $y_i$, define $\zeta^*(v)$ to be the coefficient 
of the word $v$ in $\Phi_*$, denoted $\bigl(\Phi_*|v)$.
Since the exponential ``correction'' factor is a power series in $y_1$, it 
follows that for any convergent word $v$ (i.e. any word in the $y_i$
not starting with $y_1$), we have $\zeta^*(v)=\zeta(v)$.  Inversely, the
stuffle-regularized values $\zeta^*(1,\ldots,1)$ come entirely from the
correction factor and are all polynomials in the single zeta values 
$\zeta(n)$; we see for instance that 
$$\zeta^*(1)=\zeta(1)=0,\ \ 
\zeta^*(1,1)=-{{1}\over{2}}\zeta(2),\ \ 
\zeta^*(1,1,1)={{1}\over{3}}\zeta(3),$$
$$\zeta^*(1,1,1,1)=-{{1}\over{4}}\zeta(4)+{{1}\over{8}}\zeta(2)^2=
-{{1}\over{4}}\zeta(4)+{{5}\over{16}}\zeta(4)={{1}\over{16}}\zeta(4);$$
thus, we can write the correction factor as
$$exp\Bigl(\sum_{n\ge 1} {{(-1)^{n-1}}\over{n}}\zeta(y_n)y_1^n\Bigr)
=\sum_{n\ge 1} \zeta^*(\underbrace{1,\ldots,1}_n)y_1^n.\eqno(1.1.10)$$
For words of the form $w=y_1^iv$ with $v$ a word in the $y_i$ not starting with $y_1$,
the {\it stuffle regularized multizeta values} are given by the formula
$$\zeta^*(w)=\Bigl(\Phi_*|v\Bigr)=
\sum_{j=0}^i \zeta^*(\underbrace{1,\ldots,1}_j) \bigl(\Phi|y_1^{i-j}v\bigr).\eqno(1.1.11)  $$

The values $\zeta^*(v)$ are called the {\it stuffle regularization} 
of the convergent multizeta values, because of the following theorem.
\vskip .2cm
\noindent {\bf Theorem 1.1.2.} {\it For all words $u$, $v$ in the variables $y_i$,
the values $\zeta^*(v)$ satisfy the {\rm stuffle relations} 
$$\zeta^*\bigl(st(u,v)\bigr)=\zeta^*(u)\zeta^*(v).\eqno(1.1.12)$$}
\vskip .1cm
\noindent {\bf Remark.} Theorems 1.1.1 and 1.1.2 are part of the classical
theory of multizeta values, proved originally by Drinfel'd in the form of the
two following statements on $\Phi_{KZ}$:
\vskip .1cm
\noindent (i) $\Phi_{KZ}\in \Q\langle\langle x,y\rangle\rangle$ is group-like for
the coproduct $\Delta$ defined by $\Delta(x)=x\otimes 1+1\otimes x$,
$\Delta(y)=y\otimes 1+1\otimes y$.
\vskip .1cm
\noindent (ii) $\Phi^*\in \Q\langle\langle y_1,y_2,\ldots\rangle\rangle$ is group-like
for the coproduct $\Delta^*$ defined by
$$\Delta^*(y_i)=\sum_{k+l=i} y_k\otimes y_l.$$
Theorems 1.1.1 and 1.1.2 are direct translations of these two properties on
power series into multiplicative properties of the coefficients of those
power series (cf. [R] for a detailed exploration of these facts).
\vskip .3cm
\noindent {\bf Definition.} Let ${\cal Z}$ denote the $\Q$-algebra generated 
by the convergent multizeta values under the multiplication law (1.1.3).  By
(1.1.2) and (1.1.11), ${\cal Z}$ contains all the shuffle and 
stuffle regularized multizeta values.  For every word $w\in \Q\langle 
x,y\rangle$ of length §(i.e.~degree) $n$ containing $r$ $y$'s, 
the corresponding multiple 
zeta value $\zeta(w)$ is said to be of {\it weight} $n$ and {\it depth} $r$.
For each $n\ge 0$, let ${\cal Z}_n$ denote the $\Q$-vector space generated by the convergent
multiple zeta values of weight $n$.  We have ${\cal Z}_0=\Q$, 
${\cal Z}_1=\langle 0 \rangle$, ${\cal Z}_2=\langle\zeta(2)\rangle$.
\vskip .2cm
The algebra ${\cal Z}$  has a rich structure of which the shuffle and stuffle
families of algebraic relations (known as the double shuffle relations) are only
one aspect.  There are many other known algebraic relations between elements
of ${\cal Z}$, and also, of course, difficult problems of transcendence
and irrationality.  Few results are known on the transcendence;
the fundamental conjecture that all multiple zeta values are transcendent 
still seems far out of reach.

The transcendence conjecture can be subsumed into the following 
seemingly simple structural conjecture on ${\cal Z}$.
\vskip .2cm
\noindent {\bf Main transcendence conjecture.} {\it The weight provides a
grading of the $\Q$-algebra ${\cal Z}$; in other words, there are no
linear relations between multizeta values of different weights.}
\vskip .2cm
This assumption indeed implies that every multizeta value is transcendent,
since otherwise, if some $\zeta$ of weight $n$ were algebraic, there would 
be a minimal polynomial $P(x)$ such that $P(\zeta)=0$; each term of the 
polynomial would be a $\zeta^i$, which when expanded out as a sum
by the shuffle multiplication rule would yield a non-zero linear combination
of multizetas of weight $in$, and the sum of all these terms of different
weights would be zero, contradicting the main conjecture.

The conjectures concerning transcendence seem unprovable for the
time being, but the combinatorial/algebraic structure of the multizeta
algebra is still a rich subject of study, with another conjecture specifically
concerning algebraic relations.

\vskip .2cm
\noindent {\bf Main algebraic conjecture.}  {\it The ``regularized'' double shuffle 
relations (1.1.3) and (1.1.12) generate all algebraic relations between multizeta 
values.}

\vskip .2cm
This conjecture makes it natural to focus attention on the double shuffle relations.
For this purpose, it is useful to define a {\it formal multiple zeta algebra} of 
transcendent symbols satisfying only the regularized double shuffle relations,
and investigate its structure. This algebra, defined in the next section, 
is one of the main objects of study in the theory of multiple zeta values.
\vskip .8cm
\noindent {\bf \S 1.2. Formal multiple zeta values}
\vskip .5cm
For every word $w$ in $x$ and $y$, let $\overline Z(w)$ denote a formal symbol associated to $w$, and let
$\Q[\overline Z(w)]$ be the commutative $\Q$-algebra generated as a vector space by these symbols, equipped with the
multiplication law
$$\overline Z(u)\overline Z(v)=\overline Z\bigl(sh(u,v)\bigr).\eqno(1.2.1)$$
Let ${\cal SH}$ be the quotient of $\Q[\overline Z(w)]$ by the 
linear relations analogous to (1.1.2)
$$\overline Z(w)=\sum_{a=0}^r \sum_{b=0}^s (-1)^{a+b}\overline Z\bigl(\pi(sh(y^a,y^{r-a}v
x^{s-b},x^b))\bigr)\eqno(1.2.2)$$
for every non-convergent word $w$.  As in theorem 1.1, this definition ensures that the multiplication law
(1.2.1) passes to the quotient ${\cal SH}$.  We write $\widetilde Z(w)$ for the image of $\overline Z(w)$ in ${\cal SH}$.

In analogy with (1.1.9), we define $\widetilde Z^*(\underbrace{1,\ldots,1}_n)$ to be the coefficient of 
$y_1^n$ in the formal power series with coefficients in ${\cal SH}$
$$exp\Bigl(\sum_{n\ge 1} {{(-1)^{n-1}}\over{n}}\widetilde Z(y_n)y_1^n\Bigr),$$
so they are polynomials in the $\widetilde Z(y_i)$; note that all polynomials in the $\widetilde Z(w)$ can be 
expressed as linear combinations of convergent multizetas by using the multiplication rule (1.2.1) and then (1.2.2).  
In analogy with (1.1.11), we set
$$\widetilde Z^*(w)=
\sum_{j=0}^i \widetilde Z^*(\underbrace{1,\ldots,1}_j) \bigl(\Phi|y_1^{i-j}v\bigr)=
\sum_{j=0}^i \widetilde Z^*(\underbrace{1,\ldots,1}_j) \widetilde Z(y^{i-j}v),\eqno(1.2.3)$$
for every word $w=y_1^iv$ where $v$ is a word in the $y_i$ not starting
with $y_1$; thus these values can also be expressed as linear combinations of convergent $\widetilde Z(w)$.
Therefore, ${\cal SH}$ is generated as a vector space by the $\widetilde Z(w)$ for convergent $w$.
\vskip .2cm
Let ${\cal FZ}$, the {\it formal multizeta algebra}, be the vector space quotient of ${\cal SH}$ by the relations 
$$\widetilde Z^*\bigl(st(u,v)\bigr)=\widetilde Z^*(u)\widetilde Z^*(v),$$
which although they appear algebraic, can be written as above as linear relations between the convergent
$\widetilde Z(w)$.  The multiplication (1.2.1) passes to ${\cal FZ}$, making it into a $\Q$-algebra.
We write $Z(w)$ for the image of $\widetilde Z(w)$ in ${\cal FZ}$.
\vskip .2cm
By definition, we have a surjection ${\cal FZ}\rightarrow {\cal Z}$.  But
the space ${\cal FZ}$ is easier to study than ${\cal Z}$ because the real
multizeta values satisfy unknown numbers of other relations, including, as explained in 1.1,
the fact that it is not even known whether they are transcendent, or whether
there are any linear relations between real multizeta values of different
weights.  It is tempting to conjecture that ${\cal FZ}\simeq {\cal Z}$, but
pending any kind of knowledge about the transcendence properties of real multizeta
values, we adopt the strategy of replacing the real value algebra 
by the formal multizeta algebra ${\cal FZ}$ as the main object of study in the 
combinatorial/algebraic theory of multizetas. 

By definition, ${\cal FZ}$ is a graded algebra, with ${\cal FZ}_0=\Q$, ${\cal FZ}_1=0$
and ${\cal FZ}_2$ a one-dimensional space generated by $Z(2)=Z(xy)$ (as for real
multizetas, we use the notation $Z(k_1,\ldots,k_r)=Z(x^{k_1-1}y\cdots x^{k_r-1}y)$).
Let $\overline {\cal FZ}$ denote the quotient of ${\cal FZ}$ by the ideal generated
by $Z(2)$.

Let $\nfz$ denote the quotient of $\overline{\cal FZ}$ 
modulo the ideal generated by ${\cal FZ}_0$ and products ${\cal FZ}_{>0}^2$. 
Known as the {\it new formal zeta space}, lifts of its generators to 
$\overline{\cal FZ}$ form a set of ring generators.  In fact, $\nfz$ is more than
just a vector space.  An important and difficult theorem due to Racinet states that 
the dual of $\nfz$ is a Lie algebra, known as the double shuffle Lie algebra $\ds$
(see next section). Thus $\nfz$ is a Lie coalgebra, and $\overline{\cal FZ}$ is a 
Hopf algebra. In Chapter 4, we give the neat and simple theoretical 
proof of Racinet's theorem that emerges easily from Ecalle's theory.

The following section is devoted to the Lie algebra $\ds$, which is one of the main points 
of focus of the entire theory, thanks to the simplicity of its definition and the concrete 
nature of its elements, which make it into a valuable and attractive ``way in'' to
the theory, accessible to explicit computation.
\vskip .8cm
\noindent {\bf \S 1.3. The double shuffle Lie algebra $\ds$}
\vskip .4cm 
\noindent {\bf Definition 1.3.1.} The Lie algebra $\ds$ is the dual of the Lie 
coalgebra $\nfz$ of new formal multizeta values.  It can be defined directly
as the set
of polynomials $f\in \Q\langle x,y\rangle$ having the two following properties.
\vskip .2cm
\noindent (1) The coefficients of $f$ satisfy the {\it shuffle relations}
$$\sum_{w\in sh(u,v)} (f|w)=0,\eqno(1.3.1)$$
where $u,v$ are words in $x,y$ and $sh(u,v)$ is the set of words obtained
by shuffling them.  This condition is equivalent to the assertion that
$f\in {\rm Lie}[x,y]$.
\vskip .2cm
\noindent (2) Let $f_*=\pi_y(f)+f_{\rm corr}$, where $\pi_y(f)$ is the 
projection of $f$ onto just the words ending in $y$, and 
$$f_{\rm corr}=\sum_{n\ge 1} {{(-1)^{n-1}}\over{n}}(f|x^{n-1}y)y^n.\eqno(1.3.2)$$
(When $f$ is homogeneous of degree $n$, which we usually assume, then
$f_{\rm corr}$ is just the monomial ${{(-1)^n}\over{n}}(f|x^{n-1}y)y^n$.)
The coefficients of $f_*$ satisfy the {\it stuffle relations}:
$$\sum_{w\in st(u,v)} (f_*|w)=0,\eqno(1.3.3)$$
where now $u$, $v$ and $w$ are words ending in $y$, considered as rewritten
in the variables $y_i=x^{i-1}y$, and $st(u,v)$ is the stuffle of two such 
words.
\vskip .3cm
For every $f\in {\rm Lie}[x,y]$, define a derivation $D_f$ of 
${\rm Lie}[x,y]$ by setting it to be
$$D_f(x)=0,\ \ D_f(y)=[y,f]$$
on the generators.  Define the {\it Poisson bracket} on (the underlying
vector space of) ${\rm Lie}[x,y]$ by 
$$\{f,g\}=[f,g]+D_f(g)-D_g(f).\eqno(1.3.4)$$
This definition corresponds naturally to the Lie bracket on the space of
derivations of ${\rm Lie}[x,y]$; indeed, it is easy to check that
$$[D_f,D_g]=D_f\circ D_g-D_g\circ D_f=D_{\{f,g\}}.\eqno(1.3.5)$$
\vskip .1cm
\noindent {\bf Definition 1.3.2.} Let ${\bf L}$ denoted the Lie algebra generated by 
the polynomials $C_i=ad(x)^{i-1}(y)$, $i\ge 1$ inside $\Q\langle x,y\rangle$.
We have ${\rm Lie}[x,y]=\Q x\oplus {\bf L}$, and it is a standard result of 
Lazard elimination that the $ad(x)^{i-1}(y)$ generate ${\bf L}$ freely.
The {\it twisted Magnus Lie algebra $\mt$} is defined to be
the Lie algebra whose underlying vector space is ${\bf L}$,
but equipped with the Poisson bracket (1.3.4). 
\vskip .3cm
In his 2000 Ph.D. thesis, G. Racinet proved the following theorem, using 
a complicated series of arguments later condensed and reworked in the appendix
to [Furusho].  In Chapter 4 of this text, we show how this result drops naturally
and easily out of Ecalle's theory once the basic machinery has been established.
\vskip .3cm 
\noindent {\bf Theorem 1.3.3.} {\it The double shuffle space $\ds$ is a Lie
algebra under the Poisson bracket, i.e. $\ds$ is a Lie subalgebra of $\mt$.}
\vskip .2cm
This theorem raises the question of the Lie algebra structure of $\ds$, which has
given rise to a great deal of conjectures and computations.  
\vskip .3cm
\noindent {\bf Structure conjecture for $\ds$.} {\it The Lie algebra $\ds$ is
freely generated by one generator of weight $n$ for each odd $n\ge 3$.}
\vskip .3cm
In 2010, an important breakthrough by F. Brown concerning motivic multiple zeta values
had, as a consequence, the result that the free Lie algebra on one generator
in each odd weight $\ge 3$ does have a canonical injection into $\ds$. For the
rest, this is still a wide open question.

The double shuffle Lie algebra inherits a grading from ${\rm Lie}[x,y]$, corresponding
to the degree (weight) of polynomials.  We write $\ds_n$ for the graded part of 
weight $n$. It is also equipped with an increasing depth filtration 
$$\ds^1\subset \ds^2\subset \cdots$$
where $f\in \ds$ lies in $\ds^d$ if the smallest number of $y$'s appearing
in any monomial of $f$ is greater than or equal to $d$.  The depth filtration
is not a grading because there are known (so-called ``period polynomial'') linear 
combinations of elements of depth $d$ which are themselves in depth $>d$.  
This filtration is dual to the decreasing filtration on ${\cal Z}$ given by letting 
the depth of $\zeta(k_1,\ldots,k_r)$ be equal to $r$.  Again, this is a filtration 
rather than a grading since there can be linear relations mixing depths.  The first
example was already known to Euler: $\zeta(2,1)=\zeta(3)$.  

The following theorem is more or less ``folklore'', but the only published proof
so far appears to be the one in [IKZ] (which actually proves the slightly
stronger Theorem 1.4.1 in the next section), which uses some rather astute
combinatorics.
\vskip .2cm
\noindent {\bf Theorem 1.3.4} {\it Let $n\ge 3$, $d\ge 1$.  Then the quotient space 
$\ds_n^d/\ds_n^{d+1}$ is equal to 0 if $d\not\equiv n$ mod 2.}
\vskip .2cm
In Chapter 3, \S 3.4, we show how the proof of this result 
(or rather, of Theorem 1.4.1 below) falls out as an easy consequence of 
Ecalle's methods.  

Theorem 1.3.4 is just one special case of
another structure conjecture for $\ds$, that is much finer than the previous one.
Let $BK(X,Y)$ denote the Broadhurst-Kreimer function of two commutative variables
defined by
$$BK(X,Y)={{1}\over{1-{\cal O}(X)Y+{\cal S}(X)Y^2-{\cal S}(X)Y^4}},\eqno(1.3.6)$$
where ${\cal O}(X)=X^3/(1-X^2)$ and ${\cal S}(X)=X^{12}/(1-X^4)(1-X^6)$. Let
${\cal U}\ds$ denote the universal enveloping algebra of $\ds$.  Then ${\cal U}\ds$
is automatically equipped with a weight grading and depth filtration corresponding
to those of $\ds$.  The following conjecture was formulated by Broadhurst and
Kreimer for real multiple zetas, but it applies just as well to formal ones.
\vskip .3cm
\noindent {\bf Broadhurst-Kreimer structure conjecture for $\ds$.} {\it  For
all $n\ge 3$ and $d\ge 1$, the coefficient of $X^nY^d$ in the Taylor expansion
of $BK(X,Y)$ is the dimension of the graded quotient space 
${\cal U}\ds_n^d/{\cal U}\ds_n^{d+1}$.}
\vskip .3cm
Note in particular that all terms of the Taylor expansion of ${\cal O}(X)$ are of
odd degree, so in the Taylor expansion of ${\cal O}(X)Y$ the coefficients of
terms where $n\not\equiv d$ mod 2 are all $0$, and the same is even more obvious
for the terms ${\cal S}(X)Y^2$ and ${\cal S}(X)Y^4$ which contain only monomials in 
which $n$ and $d$ are even.  Thus Theorem 1.3.4 would be a corollary of the 
Broadhurst-Kreimer structure conjecture.  Furthermore, ignoring the depth filtration
comes down to setting $Y=1$, so the Broadhurst-Kreimer conjecture can be simplified
to a conjecture purely on the weight-grading of ${\cal U}\ds$, namely the dimension
of the graded piece  ${\cal U}\ds_n$ is given by the coefficient of $X^n$ in the
generating series
$${{1}\over{1-{\cal O}(X)}}={{1-X^2}\over{1-X^2-X^3}}.$$
This is well-known to be the generating series for the graded dimensions of the
free algebra on one generator in each odd weight $n\ge 3$, which is the universal
enveloping algebra of the free Lie algebra on the same generators.  Thus the
Broadhurst-Kreimer conjecture also implies the free-generation structure conjecture 
on $\ds$ given above.

\vskip .8cm
\noindent {\bf \S 1.4. The linearized double shuffle space}
\vskip .3cm
\noindent {\bf Definition 1.4.1.}
The {\it linearized double shuffle space} $\ls$ is defined to be the set of
polynomials in $x,y$ of degree $\ge 3$ satisfying the shuffle relations (1.3.1)
(i.e. belonging to the free Lie algebra ${\rm Lie}[x,y]$) and a second
set of relations given by
$$\sum_{w\in sh(u,v)} (\pi_y(f)|w)=0,\eqno(1.4.1)$$
where $\pi_y(f)$ is the projection of $f$ onto the words ending in $y$, rewritten
in the variables $y_i=x^{i-1}y$, $u,v$ are words in the $y_i$ and $w$ belongs
to their shuffle in the alphabet $y_i$.  However, we {\it exclude} from
$\ls$ all (linear combinations of) the depth 1 even degree
polynomials, namely $ad(x)^{2n+1}(y)$, $n\ge 1$.  Note that the condition 
(1.4.1) is empty on the depth 1 polynomials, so including or excluding them 
is essentially a convention.

The space $\ls$ is not only graded by weight, but also by depth, since unlike the
stuffle relations (1.3.1), the shuffle relations (1.4.1) respect the depth.
We write as usual $\ls_n$ for the graded part of weight $n$ and $\ls^d$ for the
graded part of depth $d$.  
\vskip .3cm
\noindent {\bf Proposition 1.4.1.} {\it The associated graded for the
depth filtration of $\ds$ is contained in $\ls$; i.e. in weight $n\ge 3$ and
depth $d\ge 1$, we have
$$\ds_n^d/\ds_n^{d+1}\subset \ls_n^d.\eqno(1.4.2)$$}
\vskip .3cm
\noindent {\bf Proof.} It is immediate that for any $f\in \ds$, if $\overline f$
is obtained from $f$ by taking only the terms of minimal depth (i.e. minimal 
number of $y$'s), then $\overline f\in \ls$.  Indeed, if $d$ is the 
(minimal) depth of $g$, then the stuffle relations of depth $d$ are actually 
shuffle relations since the additional terms in the stuffle where indices are 
``stuffed'' together are words of smaller depth, and therefore have 
coefficient $0$ in $f$.  Thus the truncations in minimal weight of elements
$f\in\ds$ all satisfy the linearized double shuffle relations, showing
(1.4.2).

The only point that needs some care is the case $d=1$, where the odd degree
polynomials $ad(x)^{2n+1}(y)$ have been excluded from $\ls$.  Therefore
we need a separate argument in order to check (1.4.2) in the case $d=1$;
it is necessary to show that there is no element in $\ds$ of depth 1 and
even weight.  The proof we give here appears in complete detail in
[C, Theorem 2.30 (i)].  By explicitly solving the depth 2 stuffle relations
for $f\in \ds$, given by
$$(f|x^iyx^{n-2-i}y)+(f|x^{n-2-i}yx^iy)+(f|x^{n-1}y)=0,\eqno(1.4.3)$$
one finds that
$$(f|x^{n-2}y^2)={{n-1}\over{2}}(f|x^{n-1}y).\eqno(1.4.4)$$
Now suppose that $f\in\ds$ is of even weight $n$ and
of depth $1$, i.e. the coefficient $(f|x^{n-1}y)\ne 0$. Since every
Lie polynomial satisfies $f=(-1)^{n-1}\overleftarrow f$ where 
$\overleftarrow f$ denotes the polynomial $f$ written backwards (i.e. with
each monomial in $x$ and $y$ written backwards), if $n$ is even then 
$f$ can contain no palindromic words.  Therefore in particular
$(f|yx^{n-2}y)=0$, and so the relation (1.4.3) for $i=n-2$, given by
$$(f|x^{n-2}yy)+(f|yx^{n-2}y)+(f|x^{n-1}y)=0,$$
simplifies to
$$(f|x^{n-2}yy)=-(f|x^{n-1}y),$$
contradicting (1.4.4). This concludes the proof that no depth 1 element of
even weight can exist in $\ds$, and therefore $\ds_n^1/\ds_n^2\subset\ls_n$.
\hfill{$\square$}
\vskip .3cm
The above result is actually the motivation for dropping the even depth 1
Lie polynomials from $\ls$. It is an open question
whether the inclusion (1.4.2) is also a surjection, i.e. whether every element
of the linearized double shuffle space is the lowest-depth part of some
double shuffle element.

The stronger version of Theorem 1.3.4 also holds for $\ls$.
\vskip .3cm
\noindent {\bf Theorem 1.4.1.} {\it The subspace $\ls_n^d$ of $\ls$ is zero if
$n\not\equiv d$ mod 2.}
\vskip .2cm
By (1.4.2), Theorem 1.3.4 is an immediate consequence of 
this one. As explained in the previous section, we give a simple proof of 
Theorem 1.4.1 using Ecalle's methods in Chapter 3, \S 3.4.
\vfill\eject
\ \vskip 1cm
\centerline{\bf Chapter 2} 
\vskip .5cm
\centerline{\bf The Lie algebra \ARI}
\vskip .8cm
\noindent {\bf \S 2.1. Moulds and bimoulds} 
\vskip .4cm
We work over the field ${\Bbb C}$ of complex numbers. Let $u_1,u_2,\ldots $ and 
$v_1,v_2,\ldots$ denote two infinite 
sequences of indeterminates.  A {\it bimould} $M$ is a collection of functions
$$M_r\pmatrix{u_1&u_2&\cdots&u_r\cr v_1&v_2&\cdots&v_r}$$
for each $r\ge 0$, where each $M_r$ is a function of the $2r$ variables 
$u_1,\ldots,u_r,v_1,\ldots,v_r$ (in particular $M_0$ is a constant).  These functions are a priori arbitrary, but later, in the context
of the study of multizeta values, we will restrict our attention to rational functions,
polynomials, and constants.  A {\it mould} is a bimould that is actually only a
function of the $u_i$, and a {\it $v$-mould} is a function only of the $v_i$.
Most of the time, when there is no risk of confusion, we drop the index $r$ and
write $M\pmatrix{u_1&u_2&\cdots&u_r\cr v_1&v_2&\cdots&v_r}$ for
$M_r\pmatrix{u_1&u_2&\cdots&u_r\cr v_1&v_2&\cdots&v_r}$, the {\it depth} $r$ being 
indicated automatically by the number of variables. However, on occasion when
working with a specific mould it may be necessary to use the index for precision; 
for example the mould $M_2(u_1,u_2)=u_2$ is different from the mould 
$M_3(u_1,u_2,u_3)=u_2$. We write $M(\emptyset)$ for $M_0$.  The space of all bimoulds
is denoted BIMU.

\vskip .2cm
Two moulds or bimoulds $M,N\in BIMU$ can be added, multiplied and, if $N(\emptyset)=0$,
composed.  Writing $w_i=\bigl({{u_i}\atop{v_i}}\bigr)$ (or considering the
variables $w_i$ as belonging to an arbitrary alphabet), we have
\vskip .2cm
$$\eqalign{(M+N)(w_1,\ldots,w_r)=M(w_1,\ldots,w_r)+N(w_1,\ldots,w_r)\cr
mu(M,N)(w_1,\ldots,w_r)=\sum_{0\le i\le r} M(w_1,\ldots,w_i)N(w_{i+1},
\ldots,w_r)\cr
(M\circ N)(w_1,\ldots,w_r) =\sum_{{{\w={\bf w_1\cdots w_s}}\atop 
{\w_i\ne \emptyset}}} M(|\w_1|,\cdots,|\w_s|)
N(\w_1)\cdots N(\w_s).}\eqno(2.1.1)$$

Here, $|(w_1,\ldots,w_r)|$ denotes the single-letter word $w_1+\cdots+w_r$,
which is $\bigl({{u_1+\cdots+u_r}\atop{v_1+\cdots+v_r}}\bigr)$ in the bimould case.
\vskip .2cm
\noindent {\bf Remark.}  Moulds are generalizations of power series.  If a 
mould $M$ takes constant values on each word, then it can be identified
with the power series
$$M=\sum_{(w_1,\ldots,w_r)} M(w_1,\ldots,w_r)w_1\cdots w_r.$$
\noindent {\bf Exercise.} Check that in the power series case, the rules for 
addition, multiplication and composition are just the usual ones.
\vskip .3cm
\noindent {\bf Examples.} (1) The first examples are the Log and Exp moulds 
given by $Exp(\emptyset)=Log(\emptyset)=0$,
$$\cases{Log(w_1,\ldots,w_r)={{(-1)^{r+1}}\over{r}}\cr 
Exp(w_1,\ldots,w_r)={{1}\over{r!}}.}$$
\vskip .2cm
\noindent (2) The identity mould for multiplication {\bf 1} is given by 
{\bf 1}$(\emptyset)=1$ and all other values are $0$.  
\vskip .2cm
\noindent (3) The identity mould
{\bf Id} for composition is given by 
$$\hbox{\bf Id}(w_1,\ldots,w_r)=\cases{0&for $r=0$ and all $r>1$\cr
1&for $r=1$.}$$
\vskip .1cm
\noindent {\bf Exercise.} Show that on the one-letter alphabet $T=\{t\}$,
$Exp$ is the mould corresponding to the power series $e^t-1$,
$Log$ to log$(1+t)$ and $Id$ to $t$.  Show that as expected, 
$Exp\circ Log=\ ${\bf Id}.  
\vskip .8cm
\noindent {\bf \S 2.2. The Lie algebra \ARI} 
\vskip .4cm
\noindent {\bf Definition.} Let B\ARI (resp. $\ARI$, $\overline{\ARI}$) denote the set of bimoulds 
(resp. the subspace of moulds, resp. of $v$-moulds) satisfying $A(\emptyset)=0$.  
These spaces are obviously vector spaces, and even Lie algebras under the Lie bracket $lu$ defined
by $lu(A,B)=mu(A,B)-mu(B,A)$.  But Ecalle introduces an alternative bracket, 
the $ari$-bracket, making the same underlying vector space into a different Lie algebra.  
In chapter 3, we will explore the analogy between the two brackets on \ARI and the two 
different Lie brackets on the free Lie algebra ${\rm Lie}[x,y]$ seen in Chapter 1. Let
us define some necessary notation for the $ari$-bracket and other operators in 
Ecalle's theory.
\vskip .2cm
\noindent {\bf Flexions.}  Let $\w=\pmatrix{u_1&\cdots&u_r\cr v_1&\cdots&v_r}.$
For every possible way of cutting the word $\w$ into three (possibly empty) subwords 
$\w=\a\b\c$ with
$$\a=\bigl({{u_1,\ldots,u_k}\atop{v_1,\ldots,v_k}}\bigr),\ \ , 
\b=\bigl({{u_{k+1},\ldots,u_{k+l}}\atop{v_{k+1},\ldots,v_{k+l}}}\bigr),\ \ , 
\c=\bigl({{u_{k+l+1},\ldots,u_r}\atop{v_{k+l+1},\ldots,v_r}}\bigr),$$
set 
$$\cases{\lceil\c=\c&if $\b=\emptyset$\cr
\a\rceil=\a&if $\b=\emptyset$\cr
\b\rfloor=\b&if $\c=\emptyset$\cr
\lfloor\b=\b&if $\a=\emptyset$,}$$
otherwise
$$\cases{\lceil\c=\pmatrix{u_{k+1}+\cdots+u_{k+l+1}&u_{k+l+2}&\cdots&u_r\cr
v_{k+l+1}&v_{k+l_2}&\cdots&v_r}&if $\b\ne \emptyset$\cr
\a\rceil=\pmatrix{u_1&u_2&\cdots&u_k+u_{k+1}+\cdots+u_{k+l}\cr
v_1&v_2&\cdots&v_k}&if $\b\ne\emptyset$\cr
\b\rfloor=\pmatrix{u_{k+1}&u_{k+2}&\cdots&u_{k+l}\cr
v_{k+1}-v_{k+l+1}&v_{k+2}-v_{k+l+1}&\cdots&v_{k+l}-v_{k+l+1}}&if $\c\ne\emptyset$\cr
\lfloor\b=\pmatrix{u_{k+1}&u_{k+2}&\cdots&u_{k+l}\cr
v_{k+1}-v_k&v_{k+2}-v_k&\cdots&v_{k+l}-v_k}&if $\a\ne\emptyset$.}$$

\noindent {\bf Definition.} 
For every bimould $B\in \B\ARI$, we define operators $amit(B)$ and $anit(B)$ on B\ARI as follows: 
$$amit(B)\cdot A({\bf w})=\sum_{{{\w=\a\b\c}\atop{\b,\c\ne\emptyset}}}
A(\a\lceil\c)B(\b\rfloor),\eqno(2.2.1)$$
$$anit(B)\cdot A({\bf w})=\sum_{{{\w=\a\b\c}\atop{\a,\b\ne\emptyset}}} 
A(\a\rceil\c)B(\lfloor\b).\eqno(2.2.2)$$
For every pair of moulds $B,C\in B\ARI$, we set 
$$axit(B,C)\cdot A=amit(B)\cdot A+anit(C)\cdot A\eqno(2.2.3)$$
and 
$$arit(B)\cdot A=axit(B,-B)\cdot A=amit(B)\cdot A-anit(B)\cdot A.\eqno(2.2.4)$$
We have the following explicit expression for $arit(B)$:
$$\bigl({arit}(B)\cdot A\bigr)({\bf w})=
\sum_{{{\w=\a\b\c}\atop{\b,\c\ne \emptyset}}}
A(\a\lceil\c)B(\b\rfloor)-
\sum_{{{\w=\a\b\c}\atop{\a,\b\ne \emptyset}}}
A(\a\rceil\c)B(\lfloor\b).\eqno(2.2.5)$$ 
For $A\in \ARI$ (resp. $\overline{\ARI}$) we define the analogous operators on
$\ARI$ (resp. $\overline{\ARI}$) by dropping the lower (resp. upper) flexion signs in 
(2.2.1), (2.2.2) and (2.2.5).
\vskip .3cm
\noindent {\bf Proposition 2.2.1.} {\it For all bimoulds $B\in \B\ARI$ (resp. moulds 
$B\in \ARI$, resp. $v$-moulds $B\in \overline{\ARI}$), the operators $amit(B)$, $anit(B)$ and
$arit(B)$ are derivations for the $lu$-bracket.}
\vskip .2cm
The proof of this proposition is given in \S A.1 of the Appendix.
\vskip .3cm
Define a ``pre-Lie'' operation on $\B\ARI$ by
$$\eqalign{preari(A,B)({\bf w})&=\bigl(arit(B)\cdot A+mu(A,B)\bigr)({\bf w})\cr
&= \sum_{{{\w=\a\b\c}\atop{\b\ne\emptyset}}} 
A(\a\lceil\c)B(\b\rfloor)
-\sum_{{{\w=\a\b\c}\atop{\a,\b\ne \emptyset}}} A(\a\rceil\c)B(\lfloor\b),}\eqno(2.2.6)$$
Then the $ari$-bracket is defined on $\B\ARI$ by the formula
$$ari(A,B)=preari(A,B)-preari(B,A),\eqno(2.2.7)$$
so it is given explicitly by the formula
$$ari(A,B)({\bf w})=
\sum_{{{\w=\a\b\c}\atop{\b\ne\emptyset}}} \Bigl(A(\a\lceil\c)B(\b\rfloor) 
- B(\a\lceil\c)A(\b\rfloor)\Bigr)\qquad\qquad$$ 
$$\qquad\qquad\qquad\qquad -\sum_{{{\w=\a\b\c}\atop{\a,\b\ne \emptyset}}} \Bigl(A(\a\rceil\c)B(\lfloor\b)
-B(\a\rceil\c)A(\lfloor\b)\Bigr) .\eqno(2.2.8)$$

Notice that we then have the ``Poisson bracket'' type identity\footnote{$^*$}{cf. 
{\it \ARI/G\ARI et la d\'ecomposition des multiz\^etas en irr\'eductibles}, p. 28 (75) 
and p. 29 (84)).}  
$$ari(A,B)=arit(B)\cdot A-arit(A)\cdot B+lu(A,B).\eqno(2.2.9)$$
This analogy with the situation of two non-commutative free variables $x,y$ as
in Chapter 1, and further analogies with the group laws in the next section,
will be explained in Chapter 3.  As above, the operators $preari$ and $ari$ can 
be defined on $\ARI$ resp. $\overline{\ARI}$ by dropping the lower 
resp.~upper flexion 
signs from their defining formulas.
\vskip .3cm
\noindent {\bf Proposition 2.2.2.} {\it The $ari$-bracket is a Lie 
bracket, therefore $\B\ARI$ (and a fortiori $\ARI$ and $\overline{\ARI}$) 
are Lie algebras under $ari$.}
\vskip .2cm
\noindent {\bf Proof.} Let $\B\ARI_{lu}$ denote the vector space $\B\ARI$ 
made into a Lie algebra by equipping it with the Lie bracket $lu$.
Let $D_{arit}$ denote the image in the space ${\rm Der}\,\B\ARI_{lu}$
of the map 
$$\eqalign{\B\ARI&\rightarrow {\rm Der}\,\B\ARI_{lu}\cr
P&\mapsto arit(P).}$$
Then we have a linear isomorphism $\B\ARI\rightarrow D_{arit}$, 
and the identity
$$arit\bigl(ari(A,B)\bigr)=arit(A)\circ arit(B)-arit(B)\circ arit(A)
=[arit(A),arit(B)]$$
shows that the $ari$-bracket on $\B\ARI$ is nothing other than the
restriction to $D_{arit}$ of the usual bracket of derivations
on ${\rm Der}\,\B\ARI_{lu}$.
\hfill{$\square$}
\vskip .8cm
\noindent {\bf \S 2.3. Symmetrality, alternality, symmetrility, alternility} 
\vskip .4cm
For the study of multizeta values, Ecalle introduces four fundamental 
symmetries.
\vskip .3cm
\noindent {\bf Symmetrality and alternality.}
The first two symmetries are based on the shuffle product defined in \S 1.1.
\vskip .4cm
\noindent {\bf Definition.}
A bimould (resp. mould resp. $v$-mould) is said to be {\it symmetral} if it has
constant term 1 and
$$M\bigl(sh(u,v)\bigr)=M(u)M(v)\ {\rm for\ all\ words\ } u,v,\eqno(2.3.1)$$
and {\it alternal} if it has constant term 0 and
$$M\bigl(sh(u,v)\bigr)=0\ {\rm for\ all\ words\ } u,v.\eqno(2.3.2)$$
Note that it suffices to check both relations for the pairs
$(u,v)=\bigl(w_1,\ldots,w_s),(w_{s+1},\ldots,w_r)\bigr)$
for $1\le s\le [r/2]$ since all shuffle relations can be deduced from these
by variable change.
\vskip .3cm
%
\noindent {\bf Examples.} The alternality condition in depth 2 is
$$M\bigl(sh((u_1),(u_2))\bigr)=M(u_1,u_2)+M(u_2,u_1)=0.$$  
In depth 3, there is again only one condition to check,
namely 
$$M\bigl(sh((u_1),(u_2,u_3))\bigr)=M(u_1,u_2,u_3)+M(u_2,u_1,u_3)+M(u_2,u_3,u_1).$$
The other shuffle condition $M\bigl(sh((u_1,u_2),(u_3))\bigr)=0$ is automatically
satisfied if this one is, by the variable change $u_3\mapsto u_1$, $u_1\mapsto u_2$,
$u_2\mapsto u_3$.  
In depth 4, there are two necessary conditions for alternality, namely
$$M\bigl(sh((u_1),(u_2,u_3,u_4))\bigr) =M(u_1,u_2,u_3,u_4)+M(u_2,u_1,u_3,u_4)$$
$$+M(u_2,u_3,u_1,u_4)+M(u_2,u_3,u_4,u_1)=0$$
and
$$M\bigl(sh((u_1,u_2),(u_3,u_4))\bigr)
=M(u_1,u_2,u_3,u_4)+M(u_1,u_2,u_3,u_4)+ M(u_1,u_3,u_4,u_2)$$
$$+M(u_3,u_1,u_2,u_4)+M(u_3,u_1,u_4,u_2)+M(u_3,u_4,u_1,u_2)=0.$$
\vskip .2cm
\noindent {\bf Symmetrility and alternility.}
In this text we only define the second set of symmetries for moulds in the 
$v_i$, although Ecalle's flexion unit definition 
works for all bimoulds (cf.  {\it Flexion structure...}, p. 64-68.). 
These relations are deduced from the stuffle product introduced
in \S 1.1. Recall that on an additive alphabet ${\cal X}$
the stuffle product is given by (1.1.5).
To establish the symmetrility/alternility relations, we do not need to work 
with actual sequences; only the lengths of the sequences count. Let us write
$u=(v_1,\ldots,v_r)$, $v=(v_{r+1},\ldots,v_{r+s})$ for indetermines $v_i$,
and set
$$st(r,s)=st(u,v).$$

Let $M$ be a mould. For each stuffle sum $st(r,s)$, we define a 
symmetrality/alternility sum of terms in $M$, by associating a specific 
term to each word in (1.1.5) as follows.  For each $\sigma\in Sh^{\le}(r,s)$,
let $I_\sigma\subset \{1,\ldots,N\}$ be the set of indices $i$ such that
$|\sigma^{-1}(i)|=2$.  To each word $c^\sigma(u,v)$ as in (1.1.6), we
associate a set of $2^{|I_\sigma|}$ words indexed by the subsets
$J\subset I_\sigma$ (including the empty set), defined as follows:
$$C^\sigma_J=(d_1,\ldots,d_N)$$
where we write $\sigma^{-1}(i)=\{k_\sigma, l_\sigma\}$ with 
$k_\sigma<l_\sigma$ for all $i\in I_\sigma$, and
$$d_i=\cases{v_{\sigma^{-1}(i)}&if $|\sigma^{-1}(i)|=1$\cr
v_{k_\sigma}&if $|\sigma^{-1}(i)|=2$ and $i\not\in I_\sigma$\cr
v_{l_\sigma}&if $|\sigma^{-1}(i)|=2$ and $i\in I_\sigma$.}$$
Note that if $I_\sigma=\emptyset$ then $C^\sigma_\emptyset=c^\sigma(u,v)$.
We set
$$M_{r,s}=\sum_{\sigma\in Sh^{\le}(r,s)} M^\sigma_{r,s}\eqno(2.3.3)$$
where 
$$M^\sigma_{r,s}={{1}\over{\prod_{i\in I_\sigma}(v_{k_\sigma}-v_{l_\sigma})}} \sum_{J\subset I_\sigma} (-1)^{|J|}M(C^\sigma_J).\eqno(2.3.4)$$
\vskip .3cm
\noindent {\bf Low depth.} In depth 2, 
The set $Sh^{\le}(r,s)$ contains only three maps: the identity map 
$\sigma_1$, the map $\sigma_2$ exchanging $1$ and $2$,
the map $\sigma_3:\{1,2\}\rightarrow \{1\}$ sending $1$ and $2$ to $1$. 
The corresponding words are
$$c^{\sigma_1}((v_1),(v_2))=(v_1,v_2),\
c^{\sigma_2}((v_1),(v_2))=(v_2,v_1),\ c^{\sigma_3}((v_1),(v_2))=(v_1+v_2),$$
so the stuffle sum is 
$st(1,1)=st((v_1),(v_2))=(v_1,v_2)+(v_2,v_1)+(v_1+v_2)$.
We have $I_{\sigma_1}=I_{\sigma_2}=\emptyset$, $I_{\sigma_3}=\{1\}$, and
$\sigma_3^{-1}(1)=\{k_{\sigma_3},l_{\sigma_3}\}$ with $k_{\sigma_3}=1$,
$l_{\sigma_3}=2$.  The words $C^\sigma_J$ corresponding to the two subsets
$J=\emptyset$ and $J=I_{\sigma_3}$ of $I_{\sigma_3}=\{1\}$ are
$C^{\sigma_3}_\emptyset=(v_1)$ and $C^{\sigma_3}_{I_{\sigma_3}}=(v_2)$.
The corresponding alternility terms are

$$\cases{M^{\sigma_1}_{r,s}=M(c^{\sigma_1}((v_1),(v_2))=M(v_1,v_2)\cr
M^{\sigma_2}_{r,s}=M(c^{\sigma_2}((v_1),(v_2))=M(v_2,v_1)\cr
M^{\sigma_3}_{r,s}={{1}\over{(v_1-v_2)}}\Bigl(M(v_1)-M(v_2)\Bigr),}$$
so the alternility sum in depth 2 is given by
$$M_{1,1}(v_1,v_2)=M(v_1,v_2)+M(v_2,v_1)+{{1}\over{v_1-v_2}}\bigl(M(v_1)-M(v_2)\bigr).\eqno(2.3.5)$$
In depth 3 the condition corresponding to 
$st(1,2)=st((v_1),(v_2,v_3))=(v_1,v_2,v_3)+(v_2,v_1,v_3)+(v_2,v_3,v_1)+
(v_1+v_2,v_3)+(v_2,v_1+v_3)$ is given by
$$M_{1,2}(v_1,v_2,v_3)=M(v_1,v_2,v_3)+
M(v_2,v_1,v_3)+
M(v_2,v_3,v_1)\qquad\qquad\qquad$$
$$+{{1}\over{v_1-v_2}}\Bigl(M(v_1,v_3)- M(v_2,v_3)\Bigr)
+{{1}\over{v_1-v_3}}\Bigl(M(v_2,v_1)-
M(v_2,v_3)\Bigr).$$
In depth 4, the term in $M_{2,2}$ corresponding to the word $(v_1+v_3,v_2+v_4)$ 
in the stuffle sum $st(2,2)=st((v_1,v_2),(v_3,v_4))$ is given by
$${{1}\over{(v_1-v_3)(v_2-v_4)}}\Bigl(M(v_1,v_2)-M(v_3,v_2)-M(v_1,v_4)+M(v_3,v_4)\Bigr).\eqno(2.3.6)$$
\noindent {\bf Definition.} The mould $M\in \ARI$ is said to be {\it symmetril} 
if it has constant term 1 and for all pairs $1\le r\le s$ we have
$$M_{r,s}(v_1,\ldots,v_{r+s})=M_r(v_1,\ldots,v_r)M_s(v_{r+1},\ldots,v_{r+s}),
\eqno(2.3.7)$$
and {\it alternil} if it has constant term 0 and for all pairs we have
$$M_{r,s}(v_1,\ldots,v_{r+s})=0.\eqno(2.3.8)$$
\vskip .3cm
\noindent {\bf Remark.} If $M$ is a polynomial-valued mould, then the
alternility sums are polynomials. To see this, it suffices to note that setting
$v_{k_\sigma}=v_{l_\sigma}$ for any $\sigma\in I_\sigma$,
in the numerator of $M_{r,s}^\sigma$ yields zero, canceling out the pole
in (2.3.4).
\vskip .8cm
\noindent{\bf \S 2.4. $Swap$ commutation in \ARI}
\vskip .5cm
We begin this section by defining some of the main {\it mould operators}.
Let $push$, $neg$, $anti$, $mantar$, $circ$, and
$swap$ be the operators on bimoulds defined as follows: 
$$\eqalign{push(M)\pmatrix{u_1&u_2&\cdots&u_r\cr v_1&v_2&\cdots&v_r}&=
M\pmatrix{-u_1-\cdots-u_r&u_1&\cdots&u_{r-1}\cr -v_r&v_1-v_r&\cdots&v_{r-1}-v_r}\cr
neg(M)\pmatrix{u_1&u_2&\cdots&u_r\cr v_1&v_2&\cdots&v_r}&=
M\pmatrix{-u_1&-u_2&\cdots&-u_r\cr -v_1&-v_2&\cdots&-v_r}\cr
anti(M)\pmatrix{u_1&u_2&\cdots&u_r\cr v_1&v_2&\cdots&v_r}&=
M\pmatrix{u_r&u_{r-1}&\cdots&u_1\cr v_r&v_{r-1}&\cdots&v_1}\cr
mantar(M)\pmatrix{u_1&u_2&\cdots&u_r\cr v_1&v_2&\cdots&v_r}&=
(-1)^{r-1}M\pmatrix{u_r&\cdots&u_1\cr v_r&\cdots&v_1}\cr
circ(M)\pmatrix{u_1&u_2&\cdots&u_r\cr v_1&v_2&\cdots&v_r}&=
M\pmatrix{u_r&u_1&\cdots&u_{r-1}\cr v_r&v_1&\cdots&v_{r-1}}\cr
swap(M)\pmatrix{u_1&u_2&\cdots&u_r\cr v_1&v_2&\cdots&v_r}&=
M\pmatrix{v_r&v_{r-1}-v_r&\cdots&v_2-v_3&v_1-v_2\cr u_1+\cdots +u_r&u_1+\cdots+u_{r-1}&\cdots&
u_1+u_2&u_1}.}$$
The first four operators can be considered as operators only on $\ARI$ 
(resp. $\overline{\ARI}$) by ignoring the $v_i$ (resp. the $u_i$).  The $swap$, however,
exchanges the two spaces $\ARI$ and $\overline{\ARI}$.  We will make use below
of the following elementary identity, proved by simple application of the
variables changes above:
$$neg\circ push=anti\circ swap\circ anti\circ swap.\eqno(2.4.1)$$

The purpose of this section and the next one is to prove a set of fundamental 
identities expressing how $swap$ commutes with the \ARI operators
$amit$, $anit$, $arit$, $preari$, $ari$ and $preawi$ (in this section) and
with the G\ARI operators $garit$ and $gari$ (in the next one).  These 
commutations yield a set of {\it fundamental identities} that 
lie at the heart of Ecalle's theory.

Recall the definitions of the operators $amit$ and $anit$ given in \S 2.2, as well as
the definitions of the operators $axit$ and $arit$: 
$$axit(B,C)\cdot A=amit(B)\cdot A+anit(C)\cdot A,\eqno(2.4.2)$$
$$arit(B)\cdot A=axit(B,-B)\cdot A=amit(B)\cdot A-anit(B)\cdot A\eqno(2.4.3)$$
to which we now add the definition of $awit$, as follows:
$$awit(B)\cdot A=axit\bigl(B,anti\circ neg(B)\bigr)=amit(B)\cdot A
-anit\bigl(anti\circ neg(B)\bigr)\cdot A.\eqno(2.4.4)$$
In analogy to the $preari$ law 
$$preari(A,B)= arit(B)\cdot A+mu(A,B),\eqno(2.4.5)$$
we also now define the $preawi$ law
$$preawi(A,B)= awit(B)\cdot A+mu(A,B).\eqno(2.4.6)$$
\vskip .3cm
The key identities are the following ones, which are proven in \S A.2 of the Appendix:
$$swap\Bigl(amit\bigl(swap(B)\bigr)\cdot swap(A)\Bigr)=amit(B)\cdot A+mu(A,B)-
swap\Bigl(mu\bigl(swap(A),swap(B)\bigr)\Bigr),\eqno(2.4.7)$$
$$swap\Bigl(anit\bigl(swap(B)\bigr)\cdot swap(A)\Bigr)=anit\bigl(push(B)\bigr)\cdot A.\eqno(2.4.8)$$
Using these two, it is quite easy to compute the $swap$ commutations with 
$arit$, $preari$, $ari$ and $preawi$.  
Applying the identities (2.4.7) and (2.4.8) to (2.4.3) immediately yields
$$\eqalign{swap&\Bigl(arit\bigl(swap(B)\bigr)\cdot swap(A)\Bigr)\cr
&=swap\Bigl(amit\bigl(swap(B)\bigr)\cdot swap(A)\Bigr)-swap\Bigl(anit\bigl(swap(B)\bigr)\cdot
swap(A)\Bigr)\cr
&=amit(B)\cdot A+mu(A,B)-swap\Bigl(mu\bigl(swap(A),swap(B)\bigr)\Bigr)-anit\bigl(push(B)\bigr)\cdot A\cr
&=axit\bigl(B,-push(B)\bigr)\cdot A+mu(A,B)-swamu(A,B)\,}$$
\vskip -.8cm\hfill(2.4.9)
\vskip .3cm
\noindent where $swamu(A,B)=swap\bigl(mu(swap(A),swap(B))\bigr)$. Applying (2.4.7) and (2.4.8) to
(2.4.5) yields the following computation ($preira$ is defined by the first equality):
$$\eqalign{preira(A,B):&=swap\Bigl(preari\bigl(swap(A),swap(B)\bigr)\Bigr)\cr
&=swap\Bigl(arit\bigl(swap(B)\bigr)\cdot A\Bigr)+swamu(A,B)\cr
&=axit\bigl(B,-push(B)\bigr)\cdot A+mu(A,B)\cr
&=amit\bigl(B\bigr)\cdot A+anit\bigl(-push(B)\bigr)\cdot A+mu(A,B)\cr
&=arit\bigl(B\bigr)\cdot A+anit\bigl(B-push(B)\bigr)\cdot A+mu(A,B)\cr
&=preari(A,B)+anit\bigl(B-push(B)\bigr)\cdot A\cr
&=irat(B)\cdot A+mu(A,B),}$$
\vskip -.8cm\hfill(2.4.10)
\vskip .3cm
where the last line introduces the operator $irat(B)\cdot A$ given by
$$irat(B)\cdot A=axit\bigl(B,-push(B)\bigr)\cdot A.\eqno(2.4.11)$$
\vskip .3cm
\noindent Applying the same method to $ari$ yields the operator $ira$ computed as:
$$\eqalign{ira(A,B):&=swap\Bigl(ari\bigl(swap(A),swap(B)\bigr)\Bigr)\cr
&=axit\bigl(B,-push(B)\bigr)\cdot A+mu(A,B)-axit\bigl(A,-push(A)\bigr)\cdot B-mu(B,A).}$$
\hfill(2.4.12)
\vskip .3cm
\noindent Finally, we define and compute $preiwa$ as follows:
$$\eqalign{preiwa(A,B):&=swap\Bigl(preawi\bigl(swap(A),swap(B)\bigr)\Bigr)\cr
&= swap\Bigl(amit(swap(B))\cdot swap(A)\Bigr)\cr
&\ \ \ \ \ \ \ \ \ \ \ \ +swap\Bigl(anit\bigl(anti\cdot neg(swap(B))\bigr)\cdot swap(A)\Bigr)\cr
&\ \ \ \ \ \ \ \ \ \ \ \ \ \ \ \ \ \ \ \ + swap\Bigl(mu(swap(A),swap(B))\Bigr)\cr
&=amit(B)\cdot A
+anit\bigl(push\cdot swap\cdot anti\cdot neg\cdot swap(B)\bigr)\cdot A
+mu(A,B)\cr 
&=amit(B)\cdot A +anit\bigl(anti(B)\bigr)\cdot A+mu(A,B)\cr
&=iwat(B)\cdot A+mu(A,B)}$$
\vskip -.8cm\hfill(2.4.13)
\vskip .3cm
\noindent where the last line introduces the definition
$$iwat(B)\cdot A=axit\bigl(B,anti(B)\bigr)\cdot A.\eqno(2.4.14)$$
Note that an easy corollary of (2.4.12) is the following result.
\vskip .3cm
\noindent {\bf Lemma 2.4.1.} {\it If $A,B$ are $push$-invariant moulds in $\ARI$,
then
$$swap\bigl(ari(swap(A),swap(B))\bigr)=ari(A,B).\eqno(2.4.15)$$}
\noindent {\bf Proof.} By (2.2.4), we have $arit(B)=axit(B,-B)$.  If
$A$ and $B$ are $push$-invariant, then by (2.4.12) we have
$$swap\bigl(ari(swap(A),swap(B))\bigr)=arit(B)\cdot A+lu(A,B)-arit(A)\cdot B,$$
which is nothing but $ari(A,B)$ by (2.2.9).\hfill{$\square$}
\vfill\eject
\noindent {\bf \S 2.5. Special subspaces of \ARI} 
\vskip .4cm
There are many interesting subspaces of \ARI, containing only moulds having
special symmetry properties or {\it dimorphic symmetries} to use Ecalle's
term, which is to say moulds in $\ARI$ having a special symmetry property and 
whose swap, in $\overline{\ARI}$, has another.  
\vskip .3cm
\noindent {\bf Definition.} We write
\vskip .2cm
\noindent $\bullet$ $\ARI^{pol}$ (resp. $\overline{\ARI}^{pol}$, $\B\ARI^{pol}$) for the 
subspace of polynomial-valued (bi)moulds;
\vskip .2cm
\noindent $\bullet$ $\ARI_{al}$ (resp. $\overline{\ARI}_{al}$, $\B\ARI_{al}$) for the 
subspace of alternal (bi)moulds.
\vskip .3cm
Following Ecalle, we also use the notation $\ARI_{a/b}$ for moulds in $\ARI$
having the property $a$ and/or whose swap has the property $b$; for instance
we may write $\ARI_{\bullet,al}$ for moulds in $\ARI$ with alternal swap.
The most important {\it dimorphy spaces} we will consider are the following:
\vskip .3cm
\noindent $\bullet$ $\ARI_{al/al}$, the subspace of alternal moulds in
$\ARI$ whose swap is alternal in $\overline{\ARI}$, and $\ARI_{\underline{al}/
\underline{al}}$, the subspace of $\ARI_{al/al}$ of moulds that are even
functions of $u_1$ in depth 1;
\vskip .3cm
\noindent $\bullet$ $\ARI_{al*al}$, the subspace of alternal moulds in
$\ARI$ whose swap is alternal in $\overline{\ARI}$ up to addition of a constant-valued
mould, and the corresponding subspace $\ARI_{\underline{al}*\underline{al}}$
of moulds that are even functions in depth 1;
\vskip .3cm
\noindent $\bullet$ $\ARI_{al/il}$, the subspace of alternal moulds in
$\ARI$ whose swap is alternil;
\vskip .3cm
\noindent $\bullet$ $\ARI_{al*il}$, the subspace of alternal moulds in
$\ARI$ whose swap is alternil up to addition of a constant-valued mould.
\vskip .3cm
In this section we are concerned with studying the Lie algebra properties of
some of these subspaces. In particular the following result
follows immediately from the definition of the $ari$-bracket, which is made
up of operations and flexions that preserve polynomials.
\vskip .3cm
\noindent {\bf Proposition 2.5.1.} {\it The subspace $\ARI^{pol}$ is a
Lie algebra under the $ari$-bracket.}
\vskip .3cm
We also have the next, significantly more difficult result, whose detailed 
proof is given in [SS, Appendix A].
\vskip .3cm
\noindent {\bf Proposition 2.5.2.} {\it $\ARI_{al}$ and $\overline{\ARI}_{al}$ are 
Lie algebras under the ari bracket. More generally, if $A$ and $B$ are
alternal moulds, then $arit(B)\cdot A$ is alternal.}
\vskip .3cm
The main result of this section is that $\ARI_{\underline{al}/\underline{al}}$
and $\ARI_{\underline{al}*\underline{al}}$ are Lie algebras under the 
$ari$-bracket.  This result is given in Theorem 2.5.6 below. We first need
three lemmas. 
\vskip .3cm
\noindent {\bf Lemma 2.5.3.} {\it If $A\in \ARI_{al}$, then
$$anti(A)(w_1,\ldots,w_r)=(-1)^{r-1}A(w_1,\ldots,w_r),\eqno(2.5.1)$$
in other words, $A$ is $mantar$-invariant.}
\vskip .1cm
\noindent {\bf Proof.} We first show the following equality on sums of 
shuffle relations:
$$sh\bigl((1),(2,\ldots,r)\bigr)-sh\bigl((2,1),(3,\ldots,r)\bigr)
+sh\bigl((3,2,1),(4,\ldots,r)\bigr)+\cdots$$
$$\ \ \ \ \ +(-1)^{r-1}sh\bigl((r-1,\ldots,2,1),(r)\bigr)=(1,\ldots,r)+(-1)^{r-1}
(r,\ldots,1).$$
Indeed, using the recursive formula for shuffle, we can write the above sum with two terms for each shuffle, as
$$\eqalign{(1,\ldots,r)&+2\cdot sh\bigl((1),(3,\ldots,r)\bigr)\cr
&-2\cdot sh\bigl((1),(3,\ldots,r)\bigr)-3\cdot sh\bigl((2,1),(4,\ldots,r)\bigr)\cr
&+3\cdot sh\bigl((2,1),(4,\ldots,r)\bigr)+4\cdot sh\bigl((3,2,1),(5,\ldots,r)\bigr)\cr
&+\cdots+(-1)^{r-2}(r-1)\cdot sh\bigl((r-2,\ldots,1),(r)\bigr)\cr
&+(-1)^{r-1}(r-1)\cdot sh\bigl((r-2,\ldots,1),(r)\bigr)+(-1)^{r-1}(r,r-1,\ldots,1)\cr
&=(1,\ldots,r)+(-1)^{r-1}(r,\ldots,1).}$$
Using this, we conclude that if $A$ satisfies the shuffle relations, then 
$$A(w_1,\ldots,w_r)+(-1)^{r-1}A(w_r,\ldots,w_1),$$
which is the desired result.\hfill{$\square$}
\vskip .3cm
\noindent {\bf Lemma 2.5.4.} {\it $\ARI_{\underline{al}*\underline{al}}$ 
is $(neg\circ push)$-invariant.}
\vskip .3cm
\noindent {\bf Proof.}  We first deal with the case
$A\in \ARI_{\underline{al}/\underline{al}}$. Using (2.4.1) and (2.5.1), we have
$$\eqalign{neg\circ push(A)(w_1,\ldots,w_r)&=
anti\circ swap\circ anti\circ swap(A)(w_1,\ldots,w_r)\cr
&=(-1)^{r-1}anti\circ swap\circ swap(A)(w_1,\ldots,w_r)\cr
&=(-1)^{r-1}anti(A)(w_1,\ldots,w_r)\cr
&=A(w_1,\ldots,w_r),}\eqno(2.5.2)$$
which proves the result. 

To extend the argument from $\ARI_{\underline{al}/\underline{al}}$ to
$\ARI_{\underline{al}*\underline{al}}$ takes some extra arguments, that
we take here directly from [SS].
Suppose that $A\in\ARI_{\underline{al}*\underline{al}}$, so $A$ is
alternal and $swap(A)+A_0$ is alternal for some constant mould $A_0$.
By additivity, we may assume that $A$ is concentrated in depth $r$.  First
suppose that $r$ is odd.  Then $mantar(A_0)(v_1,\ldots,v_r)=(-1)^{r-1}
A_0(v_r,\ldots,v_1)$, so since $A_0$ is a constant mould, it is
mantar-invariant.  But $swap(A)+A_0$ is alternal, so it is also
mantar-invariant by Lemma B.1; thus $swap(A)$ is mantar-invariant, and
the identity $neg\circ push=mantar\circ swap\circ mantar\circ swap$
shows that $A$ is $neg\circ push$-invariant as in (B.2).

Finally, we assume that $A$ is concentrated in even depth $r$. Here
we have $mantar(A_0)=-A_0$, so we cannot use the argument above;
indeed $swap(A)+A_0$ is mantar-invariant, but
$$mantar(swap(A))=swap(A)+2A_0.\eqno(2.5.3)$$
Instead, we note that if $A$ is alternal then so is $neg(A)=A$.  Thus we can
write $A$ as a sum of an even and an odd function of the $u_i$ via the formula
$$A={{1}\over{2}}(A+neg(A))+{{1}\over{2}}(A-neg(A)).\eqno(2.5.4)$$
So it is enough to prove the desired result for all moulds concentrated in
even depth $r$ such that either $neg(A)=A$ (even functions) or
$neg(A)=-A$ (odd functions).  First suppose that $A$ is
even. Then since $neg$ commutes with $push$ and $push$ is of odd order
$r+1$ and $neg$ is of order 2, we have
$$(neg\circ push)^{r+1}(A)=neg(A)=A.\eqno(2.5.5)$$
However, we also have
$$\eqalign{neg\circ push(A)&=mantar\circ swap\circ mantar\circ swap(A)\cr
&=mantar\circ swap\bigl(swap(A)+2A_0\bigr)\ \ {\rm by\ (2.5.3)}\cr
&=mantar\bigl(A+2A_0\bigr)\cr
&=A-2A_0.}$$
Thus $(neg\circ push)^{r+1}(A)=A-2(r+1)A_0$, and this is equal to $A$
by (2.5.5), so $A_0=0$; thus in fact $A\in\ARI_{\underline{al}/\underline{al}}$
and that case is already proven.

Finally, if $A$ is odd, i.e. $neg(A)=-A$, the same argument as above
gives $A-2(r+1)A_0=-A$, so $A=(r+1)A_0$, so $A$ is a constant-valued
mould concentrated in depth $r$, but this contradicts the assumption
that $A$ is alternal since constant moulds are not alternal, unless
$A=A_0=0$.  Note that this argument shows that all moulds in
$\ARI_{\underline{al}*\underline{al}}$ that are not in
$\ARI_{\underline{al}/\underline{al}}$ must be concentrated in odd depths.
\hfill{$\square$}
\vskip .3cm
\noindent {\bf Lemma 2.5.5.} {\it $\ARI_{\underline{al}*\underline{al}}$ is $neg$-invariant and $push$-invariant.}
\vskip .3cm
\noindent {\bf Proof.}  Let $A\in \ARI_{\underline{al}*\underline{al}}$.
Because $neg(A)=push(A)$ by Lemma 2.5.4, it is enough to prove that
$neg(A)=A$. As before, we may assume that $A$ is concentrated in a fixed
depth $d$, meaning that  $A(w_1,\ldots,w_d)=0$ for all $r\ne d$.  
If $d=1$, then $A=neg(A)$ is just the assumption on 
$A$.  If $d=2s$ is even, then since $neg$ is of order 2 and commutes with 
$push$ and $push$ is of order $d+1=2s+1$, we have
$$A=(neg\circ push)^{2s+1}(A)=neg^{2s+1}(A)=neg(A).$$
If $d=2s+1$ is odd, we can write $A$ as a sum of an even and an odd part
$$A={{1}\over{2}}\bigl(A(w_1,\ldots,w_d)+ A(-w_1,\ldots,-w_d)\bigr)+
{{1}\over{2}}\bigl(A(w_1,\ldots,w_d)-A(-w_1,\ldots,-w_d)\bigr),$$ 
so we may assume that $A(w_1,\ldots,w_d)$ is odd, i.e. $neg(A)=-A$.
Then, since $A$ is alternal, using the shuffle
$sh\bigl((w_1,\ldots,w_{2s})(w_{2s+1})\bigr)$, we have
$$\sum_{i=0}^{2s}A(w_1,\ldots,w_i,w_{2s+1},w_{i+1},\ldots,w_{2s})=0.$$
Making the variable change $w_0\leftrightarrow w_{2s+1}$ gives
$$\sum_{i=0}^{2s}A(w_1,\ldots,w_i,w_0,w_{i+1},\ldots,w_{2s})=0,$$ 
which we write out as
$$\sum_{i=0}^{2s}A\pmatrix{u_1&\ldots&u_i&u_0&u_{i+1}&\ldots&u_{2s}\cr
v_1&\ldots&v_i&v_0&v_{i+1}&\ldots&v_{2s}}=0.\eqno(2.5.6)$$
Now consider the shuffle relation $sh((w_1)(w_2,\ldots,w_{2s+1}))$, which
gives
$$\sum_{i=1}^{2s+1} A(w_2,\ldots,w_i,w_1,w_{i+1},\ldots,w_{2s+1})=0.\eqno(2.5.7)$$
Set $u_0=-u_1-\cdots -u_{2s+1}$.  Since $neg\circ push$ acts
like the identity on $A$, we can apply it to each term of (2.5.7) to obtain
$$\sum_{i=1}^{2s} -A\pmatrix{u_0&u_2&\ldots&u_i&u_1&u_{i+1}&\ldots&u_{2s}\cr
v_{2s+1}&v_2-v_{2s+1}&\ldots&v_i-v_{2s+1}&v_1-v_{2s+1}&v_{i+1}-v_{2s+1}&\ldots&v_{2s}-v_{2s+1}}$$
$$-A\pmatrix{u_0&u_2&\ldots&u_{2s}&u_{2s+1}\cr -v_1&v_2-v_1&\ldots&v_{2s}-v_1&v_{2s+1}-v_1}=0.$$
We apply $neg\circ push$ again to the final term of this sum in order to
get the $u_{2s+1}$ and $v_{2s+1}$ to disappear, obtaining
$$\sum_{i=1}^{2s} -A\pmatrix{u_0&u_2&\ldots&u_i&u_1&u_{i+1}&\ldots&u_{2s}\cr
-v_{2s+1}&v_2-v_{2s+1}&\ldots&v_i-v_{2s+1}&v_1-v_{2s+1}&v_{i+1}-v_{2s+1}&\ldots&v_{2s}-v_{2s+1}}$$
$$+A\pmatrix{u_1&u_0&u_2&\ldots&u_{2s-1}&u_{2s}\cr v_1-v_{2s+1}&-v_{2s+1}&v_2-v_{2s+1}&\ldots&v_{2s-2}-v_{2s-1}&v_{2s-1}-v_{2s}}=0.$$
\vskip .1cm
Making the variable changes $u_0\leftrightarrow u_1$ and $v_1\mapsto v_0-v_1$, $v_i\mapsto v_i-v_1$ for
$2\le i\le 2s$, $v_{2s+1}\mapsto -v_1$ in this identity yields
$$\sum_{i=1}^{2s} -A\pmatrix{u_1&u_2&\ldots&u_i&u_0&u_{i+1}&\ldots&u_{2s}\cr
v_1&v_2&\ldots&v_i&v_0&v_{i+1}&\ldots&v_{2s}}
+A\pmatrix{u_0&u_1&u_2&\ldots&u_{2s-1}&u_{2s}\cr v_0&v_1&v_2&\ldots&v_{2s-1}&v_{2s}}=0.\eqno(2.5.8)$$
Finally, adding (2.5.6) and (2.5.8) yields 
$$2A\pmatrix{u_0&u_1&\ldots&u_{2s}\cr v_1&v_2&\ldots&v_{2s}}=0,$$
so $A=0$.  This concludes the proof that if $A\in \ARI_{al/al}$, then $A(w_1,\ldots,w_d)$ is an even function for
all $d>1$; thus if we assume in addition that $A$ is even for $d=1$, then $neg(A)=A$, and by Lemma 2.5.4, we have $push(A)=A$.  \hfill{$\square$}
\vskip .5cm
Finally, to prove Theorem 2.5.6, we will also need the following important 
identity that appears in Chapter 4 as Lemma 2.4.1. 
For all $push$-invariant moulds $A,B\in\ARI$, we have
$$swap\bigl(ari(A,B)\bigr)=ari\bigl(swap(A),swap(B)\bigr),\eqno(2.5.9)$$
\vskip .4cm
\noindent {\bf Theorem 2.5.6.} {\it $\ARI_{\underline{al}/\underline{al}}$ is a 
Lie algebra under the $ari$-bracket.}
\vskip .3cm
\noindent {\bf Proof.}  Let $A$, $B\in \ARI_{\underline{al}/\underline{al}}$ and 
set $C=ari(A,B)$.  The mould $C$ is alternal by Proposition 2.5.2. By Lemma
2.5.5, $A$ and $B$ are $push$-invariant, so by (2.5.9) we have
$swap(C)=swap\bigl(ari(A,B)\bigr)=ari\bigl(swap(A),swap(B)\bigr)$, which is
also alternal by Proposition 2.5.2. 
It remains only to check that $C$ is even 
in depth 1.  But in fact, $C\bigl({{u_1}\atop{v_1}}\bigr)=0$, as the depth 1 
part of an ari-bracket is always zero, which follows directly from its 
definition in (2.2.8).\hfill{$\square$}
\vskip .8cm
\noindent {\bf \S 2.6. Circ-neutrality}
\vskip .3cm
In this section we work with the circ-operator defined on bimoulds in 
\S 2.4. We say that a mould $A\in \ARI$ is {\it circ-invariant} if
$circ(A)=A$ and {\it circ-neutral} if for each depth $r>1$ we have
$$A+circ(A)+circ^2(A)+\cdots+circ^{r-1}(A)=0.$$
Most applications concern moulds in the variables $v_i$. We use the
notation $\overline{\ARI}$ to denote the space of moulds in 
$\ARI$ that are functions only of the variables $v_i$, 
$\overline{\ARI}_{circneut}$ for the space of these moulds that are 
circ-neutral, and $\overline{\ARI}_{*circneut}$ for the space of these moulds 
that are circ-neutral up to addition of a constant-valued mould. 
The following Proposition is also proved in [FK, Prop. 1.30].
\vskip .3cm
\noindent {\bf Proposition 2.6.1} {\it The space $\overline{\ARI}_{circneut}$ 
forms a Lie algebra under the $ari$-bracket.}
\vskip .2cm
\noindent {\bf Proof.} Let $A,B\in \overline{\ARI}_{circneut}$.
We need to show that 
$$\sum_{i=1}^r ari(A,B)(v_i,\ldots,v_r,v_1,\ldots,v_{i-1})=0,$$
where the formula for the $ari$-bracket is given in (2.2.9) as
$$\eqalign{ari(A,B)&=
lu(A,B)+arit(B)\cdot A-arit(A)\cdot B\cr
&=lu(A,B)+amit(B)\cdot A-anit(B)\cdot A
-amit(A)\cdot B+anit(A)\cdot B,}$$
where $lu(A,B)=mu(A,B)-mu(B,A)$ for $mu$ as in (2.1.1), and 
$arit$ and $amit$ are defined explicitly in (2.2.1) and
(2.2.2).
We will show that this expression is circ-neutral because in fact, each of the five terms in
the sum is individually circ-neutral.  

Let us start by showing this for
the first term, $lu(A,B)$.
Let $\sigma$ denote the cyclic permutation of $\{1,\ldots,r\}$ defined by
$$\sigma(i)=i+1\ \ {\rm for}\ \ 1\le i\le r-1,\ \
\sigma(r)=1.$$
By additivity, since the circ-neutrality property is depth-by-depth,
we may assume that $A$ is concentrated in depth $s$ and $B$ in
depth $t$, with $s\le t$, $s+t=r$.
In this simplified situation, we have
$$lu(A,B)(v_1,\ldots,v_r)=A(v_1,\ldots,v_s)B(v_{s+1},\ldots,v_r)
-B(v_1,\ldots,v_t)A(v_{t+1},\ldots,v_r).$$
If $s,t>1$, we have
$$\eqalign{&\sum_{i=0}^{r-1}lu(A,B)(v_{\sigma^{i}(1)},
\ldots,v_{\sigma^{i}(r)})\cr
&= \sum_{i=0}^{r-1} \Bigl(A(v_{\sigma^i(1)},\ldots,v_{\sigma^i(s)})
B(v_{\sigma^i(s+1)},\ldots,v_{\sigma^i(r)})-
B(v_{\sigma^i(1)},\ldots,v_{\sigma^i(t)})
A(v_{\sigma^i(t+1)},\ldots,v_{\sigma^i(r)})\Bigr)\cr
&= \sum_{i=0}^{r-1} \Bigl(A(v_{\sigma^i(1)},\ldots,v_{\sigma^i(s)})
B(v_{\sigma^i(s+1)},\ldots,v_{\sigma^i(r)})
-A(v_{\sigma^{i+t}(1)},\ldots,v_{\sigma^{i+t}(s)})
B(v_{\sigma^{i+t}(s+1)},\ldots,v_{\sigma^{i+t}(r)})\Bigr)\cr
&=0}$$
as the terms cancel out pairwise.

We now prove that the second term
$$\bigl(amit(B)\cdot A\bigr)(v_1,\ldots,v_r)=\sum_{i=1}^s
A(v_1,\ldots,v_{i-1}, v_{i+t},\ldots,v_r)B(v_i-v_{i+t},\ldots,v_{i+t-1}-v_{i+t})
$$
is circ-neutral.
Fix $j\in \{1,\ldots,s\}$ and consider the term
$$A(v_1,\ldots,v_{j-1},v_{j+t},\ldots,v_r)
B(v_j-v_{j+t},\ldots,v_{j+t-1}-v_{j+t}).$$
Thus for each of the other terms
$$A(v_1,\ldots,v_{i-1},v_{i+t},\ldots,v_r)
B(v_i-v_{i+t},\ldots,v_{i+t-1}-v_{i+t})$$
in the sum, with $i\in \{1,\ldots,s\}$, there is exactly one cyclic
permutation, namely $\sigma^{j-i}$, that maps this term to
$$A(v_{\sigma^{j-i}(1)},\ldots,v_{\sigma^{j-i}(i-1)},v_{\sigma^{j-i}(i+t)},\ldots,v_{\sigma^{j-i}(r)})
B(v_j-v_{j+t},\ldots,v_{j+t-1}-v_{j+t}).$$
For fixed $j\in \{1,\ldots,s\}$, the values of $k=j-i$ mod $s$ as
$i$ runs through $\{1,\ldots,s\}$ are exactly $\{0,\ldots,s-1\}$.
Therefore, the coefficient of the
term $B(v_j-v_{j+t},\ldots,v_{j+t-1}-v_{j+t})$ in the sum of
the cyclic permutations of $amit(B)\cdot A$ is equal to
$$\sum_{k=0}^{s-1} A(v_{\sigma^{k}(1)},\ldots,v_{\sigma^{k}(i-1)},v_{\sigma^{k}(i+t)},\ldots,v_{\sigma^{k}(r)}),$$
which is zero due to the circ-neutrality of $A$.
Thus the coefficient of the term
$B(v_j-v_{j+t},\ldots,v_{j+t-1}-v_{j+t})$ in the sum of the cyclic
permutations of $amit(B)\cdot A$ is zero, and this
holds for $1\le j\le s$, so the entire sum is $0$,
i.e.~$amit(B)\cdot A$ is circ-neutral.
\vskip .2cm
\noindent {\bf Example.} $s=3,t=2,r=5$. We have
$$\eqalign{
\bigl(amit(B)\cdot A\bigr)(v_1,v_2,v_3,v_4,v_5)&=
A(v_4,v_5,v_6)B(v_1-v_4,v_2-v_4,v_3-v_4)\cr
&\ \ +A(v_1,v_5,v_6)B(v_2-v_5,v_3-v_5,v_4-v_5)\cr
&\ \ \ +A(v_1,v_2,v_6)B(v_3-v_6,v_4-v_6,v_5-v_6).}\eqno(2.6.1)$$
For $(amit(B)\cdot A)$ to be circ-neutral, the sum of
the images of this expression under the five non-trivial
powers of the six-cycle permutation $\sigma=(123456)$ must be zero.
In particular, the coefficient of every factor of $B$ that occurs in that sum
must sum to zero.  Let us show this for the $B$-factor
$B(v_2-v_5,v_3-v_5,v_4-v_5)$ that arises in the second term of (2.6.1).
The terms in the complete sum containing this factor can only
come from $\sigma$ acting on the first term of (2.6.1), giving
$$A(v_5,v_6,v_1)B(v_2-v_5,v_3-v_5,v_4-v_5)$$
and from $\sigma^5$ acting on the third term of (2.6.1), giving
$$A(v_6,v_1,v_5)B(v_2-v_5,v_3-v_5,v_4-v_5).$$
Therefore the coefficient of $B(v_2-v_5,v_3-v_5,v_4-v_5)$ in the complete
sum is equal to
$$A(v_1,v_5,v_6)+A(v_5,v_6,v_1)+A(v_6,v_1,v_5)$$
which is equal to zero by the circ-neutrality of $A$. The same holds for
every $B$-factor that occurs in the sum; there will always be exactly
three possible ways to obtain it by a unique permutation acting on each
of the three terms of (2.6.1), and the coefficients will be a
circ-sum of $A$'s that add up to zero.
\vskip.3cm
To conclude the proof of the proposition, we need to prove that
the term $anit(B)\cdot A$ is also circ-neutral, but the proof
is analogous to the case of $amit$. Finally, by exchanging $A$
and $B$, this also shows that
$amit(A)\cdot B$ and $anit(A)\cdot B$ are circ-neutral.
This concludes the proof of Proposition 2.6.1.  \hfill{$\square$}
\vskip .8cm
\noindent {\bf \S 2.7. The group G\ARI}
\vskip .4cm
The last two sections of this chapter are devoted to the group G\ARI. We
begin by defining G\ARI to be the set of moulds in the variables $u_i$
with constant term $1$; similarly, we define $\GV\ARI$ to be the set of
moulds in the $v_i$ with constant term 1, and $GB\ARI$ the set of bimoulds
with constant term 1.  We will only consider G\ARI in this section, but every 
statement and definition is equally valid for $\GV\ARI$ and GB\ARI. 

We can realize G\ARI as the exponential of \ARI via the exponential map
$exp_{ari}$ defined by
$$exp_{ari}(A)=\sum_{i\ge 0} {{1}\over{i!}}preari(\underbrace{A,\ldots,A}_i)=
{\bf 1}+A+{{1}\over{2!}}preari(A,A)+{{1}\over{3!}}preari(A,A,A)+\cdots,\eqno(2.7.1)$$
where $preari(\underbrace{A,\ldots,A}_i)$ is understood to be taken from left
to right, for example $$preari(A,A,A)=preari(preari(A,A),A).$$
(Note that while in principle $preari$ is an operator on pairs of moulds from 
\ARI, the definition (2.2.6) makes perfect sense even if only the second mould 
is in \ARI and the first is an arbitrary mould.)  Indeed, since the only
condition on elements of G\ARI is to have constant term 1, moulds in the group
$exp_{ari}(\ARI)$ certainly satisfy this condition, and since like all
exponentials $exp_{ari}$ is an isomorphism, its inverse $logari$ takes 
moulds with constant term 1 to moulds with constant term 0, i.e. $logari:
\G\ARI\rightarrow \ARI$.  This shows that $\G\ARI$ is a group.

Naturally, G\ARI has subgroups corresponding to the interesting subalgebras
of $\ARI$.  The most crucial definition is the following.
\vskip .3cm
\noindent {\bf Definition.} A mould $A\in\G\ARI$ is {\it symmetral} if for all
pairs of words $\u$,$\v$ in the $u_i$, we have
$$\sum_{\w\in sh(\u,\v)} A(\w)=A(\u)A(\v).$$
We write $\G\ARI_{as}$ for the set of symmetral moulds in $\G\ARI$.
\vskip .2cm
The following basic result will be useful later on.
\vskip .3cm
\noindent {\bf Proposition 2.7.1.} {\it We have
$$exp_{ari}(\ARI_{al})=\G\ARI_{as}.$$}
\vskip .2cm
\noindent Proof. The proof was worked out completely by N.~Komiyama in
the appendix of [K] (Theorem A.7).\hfill{$\square$}
\vskip .3cm
The following lemma is a good exercise, so we merely sketch the proof.
\vskip .3cm 
\noindent {\bf Lemma 2.7.2.} {\it If $B\in\G\ARI_{as}$ and $A\in \ARI_{al}$,
then the composition $B\circ A$ is symmetral.}
\vskip .2cm
\noindent {\bf Sketch of proof.} Consider the expression for $B\circ A$ in
(2.1.1). Summing up the terms $(B\circ A)(w_1,\ldots,w_r)$ where
$\w=(w_1,\ldots,w_r)$ runs through the shuffles $sh(\u,\v)$ of two
words $\u$ and $\v$, we obtain
$$\sum_{\w\in sh(\u,\v)} \sum_{\w=\w_1\cdots \w_s} B(|\w_1|,\ldots,
|\w_s|)A(\w_1)\cdots A(\w_s).$$
Let $\u=(u_1,\ldots,u_l)$, $\v=(u_{l+1},\ldots,u_s)$. 
There are two types of decomposition $\w=\w_1\cdots \w_r$; for which
$\w_1\cdots \w_m$ is of length $l$ and $\w_{m_1}\cdots \w_r$ is of length
$r-l$, which are called ``compatible with $\u\v$'', and the ``incompatible''
ones for which there is no such division of the decomposition into two 
compatible chunks.

The proof essentially works as follows.  We fix one decomposition
of $\u\v$ into chunks $\u_1\cdots \u_s$, and then consider the corresponding 
decompositions $\w=\w_1\cdots \w_s$ of all words $\w\in sh(\u,\v)$.
If the fixed decomposition is incompatible, then we can show that
$$\sum_{\w\in sh(\u,\v)} B(|\w_1|,\ldots, |\w_s|)A(\w_1)\cdots A(\w_s)=0,$$
simply because the different shuffles that give the same term
$B(|\w_1|,\ldots,|\w_s|)$ factor out in front of a sum of terms of the
form $A(\w_1)\cdots A(\w_s)$ that is in fact a product of sums of shuffles
and is therefore zero, since $A$ is alternal.

If the fixed decomposition is compatible, then one can show what happens
in two steps.  To start with, all the terms in which 
$|\w_1|,\ldots,|\w_s|$ is not compatible with $\u$ and $\v$ in the
sense that each $|\w_i|$ is either a sum of consecutive letters of $\u$
or consecutive letters of $\v$ sum to zero as above, due to the alternality of 
$A$.  Finally, the remaining terms in the sum are sums of shuffles of the 
$|u_i|$ in the decompositions $\u\v=\u_1\cdots \u_s$, and thus they
they simplify to products due to the symmetrality of $M$. 
\hfill{$\square$} 
\vskip 1cm
\noindent {\bf \S 2.8. The group law on G\ARI}
\vskip .3cm
For each mould $B$ in GB\ARI we can associate an automorphism
of GB\ARI denoted $garit_B$ by the formula:
$$garit_B\cdot A=\sum_{{{\w=\a_1\b_1\c_1\cdots\a_s\b_s\c_s}\atop
{\b_i\ne\emptyset, \a_i\c_{i+1}\ne \emptyset}}}
A(\lceil\b_1\rceil\cdots \lceil\b_s\rceil)
B(\a_1\rfloor)\cdots B(\a_s\rfloor)invmu(B)(\lfloor \c_1)\cdots 
invmu(B)(\lfloor \c_s) \eqno(2.8.1)$$
for $s\ge 1$, where the flexions are as defined in \S 2.2 and $invmu(B)$ is of course
the inverse of $B$ for the $mu$-multiplication. Later, another automorphism
will also be very useful:
$$ganit_B\cdot A=\sum_{{{\w=\b_1\c_1\cdots \b_s\c_s}\atop {{\rm only}\ \c_s\ {\rm can\ be\ }0}}} A(\b_1\rceil\cdots \b_s\rceil)B(\lfloor(\c_1)\cdots B(\lfloor(\c_s).\eqno(2.8.2)$$
The expressions for $garit_B$ and $ganit_B$ on G\ARI and $\GV\ARI$ are obtained as 
usual from (2.8.1) by ignoring the lower resp. upper flexions. In Chapter 3,
\S 3.5, we will see the familiar expressions for these automorphisms when we 
consider the very restricted case of moulds that are power series in two 
non-commutative variables with constant term 1, forming the so-called twisted 
Magnus group.

The group law in G\ARI, denoted $gari$, is given by
$$gari(A,B)=mu(garit_B\cdot A,B).\eqno(2.8.3)$$
This law is linear in $A$, so that the product $gari(A,B)$ can be extended
from pairs of moulds in G\ARI to pairs of moulds where $A$ is arbitrary and
$B$ is in G\ARI.  By linearizing $B$, we recover the $preari$ operator. The
linearizing procedure works as follows: we set $B=1+\epsilon C$ for a mould 
$C\in \ARI$, and consider coefficients in the field $k[[\epsilon]]/
(\epsilon^2)$ if $k$ is the base field.  Then we find that
$$garit_{(1+\epsilon C)}\cdot A=A+\epsilon\ arit(C)\cdot A,$$
so 
$$\eqalign{gari(A,1+\epsilon C)&=mu\bigl(garit_{(1+\epsilon C)}\cdot A,
1+\epsilon\ C\bigr)\cr
&= mu(A+\epsilon\ arit(C)\cdot A,1+\epsilon\ C)\cr
&=A+\epsilon\ arit(C)\cdot A +\epsilon\ mu(A,C)=A+\epsilon\ preari(A,C).}
\eqno(2.8.4)$$

The inverse of a mould $B$ for the $gari$-multiplication is written $invgari(B)$.
Since \ARI is a Lie algebra for the Lie bracket $ari$, G\ARI is a pro-unipotent group.
Then $preari$ is the pre-Lie law which expresses multiplication inside the
universal enveloping algebra of \ARI of two elements in \ARI (or more generally
one element in the enveloping algebra and one in \ARI), and $exp_{ari}$ is the 
usual Lie exponential map.  Like exp of any Lie algebra, the group G\ARI acts on
the Lie algebra via an adjoint action known as $Ad_{ari}$ and defined by
$$\eqalign{Ad_{ari}(A)\cdot B
&={{d}\over{dt}}\bigl|_{t=0}\,gari(A,exp_{ari}(tB),invgari(A)\bigr)\cr
&=B+ari\bigl(logari(A),B\bigr)+{{1}\over{2}}ari\bigl(logari(A),ari\bigl(logari(A),B\bigr)+\cdots},\eqno(2.8.5)$$
or equivalently, by
$$Ad_{ari}(A)\cdot B=gari\bigl(preari(A,B),invgari(A)\bigr).\eqno(2.8.6)$$
Writing $adgari$ for the conjugation operator 
$$adgari(A)\cdot B=gari\bigl(A,B,invgari(B)\bigr),\eqno(2.8.7)$$
the following diagram then commutes (as for any Lie algebra):
$$G\ARI \buildrel{adgari(A)}\over\rightarrow G\ARI$$
$$exp_{ari}\uparrow\ \ \qquad\qquad\downarrow logari$$
$$\ARI\  \buildrel{Ad_{ari}(A)}\over\rightarrow \ \ARI,$$
where $logari$ is the inverse of the isomorphism $exp_{ari}$ (cf. [Pisa, p. 47]).

We conclude this section with the definition of the $gaxit$ operator on
GB\ARI and the $gaxi$-multiplication law on the group 
$\GAXI=\GB\ARI\times \GB\ARI$. The very general law $gaxit$, which can be restricted to 
G\ARI and $\GV\ARI$ in the usual way, gives the action of a pair of moulds on a mould, 
whereas $gaxi$ is a multiplication law on pairs of moulds.  Following [Pisa, p. 42],
set
$$gaxit_{B,C}\cdot A=
\sum_{{{\w=\a_1\b_1\c_1\cdots\a_s\b_s\c_s}\atop
{\b_i\ne\emptyset, \a_i\c_{i+1}\ne \emptyset}}}
A(\lceil\b_1\rceil\cdots \lceil\b_s\rceil)
B(\a_1\rfloor)\cdots B(\a_s\rfloor)C(\lfloor \c_1)\cdots 
c(\lfloor \c_s), \eqno(2.8.8)$$
and
$$gaxi\bigl((A,B),(C,D)\bigr)=\bigl(mu(gaxit_{C,D}\cdot A,C),mu(D,gaxit_{C,D}\cdot B\bigr).\eqno(2.8.9)$$
Thus we have $garit_A=gaxit_{A,invmu(A)}$ and $ganit_A=gaxit_{1,A}$.  Then
$$\eqalign{gaxi&\bigl((A,invmu(A)),(C,invmu(C))\bigr)=\cr
&=\Bigl(mu\bigl(garit_C\cdot A,C\bigr),mu\bigl(invmu(C),garit_C\cdot invmu(A)\bigr)
\Bigr)\cr
&=\Bigl(mu(garit_C\cdot A,C),invmu\bigl(mu(garit_C\cdot A,C)\bigr)\Bigr)}\eqno(2.8.10)$$
since $garit_C$ is a group automorphism for $mu$-multiplication. This shows that
$gaxi$ of two pairs of the form $(A,invmu(A))$ is again of that form, and
$gari(A,B)$ is just the left-hand component of (2.8.10).  In other words, G\ARI 
is identified with the subgroup of $GB\ARI\times GB\ARI$ of pairs of the form
$(A,invmu(A))$ and $gari$ is just $gaxi$ restricted to this subgroup. In later 
chapters, other specializations of $gaxit$ and $gaxi$ to specific subgroups
will be useful for certain proofs.  In Chapter 3, \S 3.5, we will also explain 
the connection between G\ARI and $gari$ and the familiar 
twisted Magnus group with its twisted Magnus multiplication.
\vskip .8cm
\noindent {\bf \S 2.9.  \'Ecalle's
first fundamental identity: $swap$ commutation in $\G\ARI$}
\vskip .3cm
In this section we introduce {\it Ecalle's first fundamental identity} (2.9.4), 
which expresses the commutation of $swap$ with $gari$. 

Let $gira(A,B)$ be the swapped $gari$-product, i.e.
$$gira(A,B):=swap\bigl(gari\bigl(swap\cdot A,swap\cdot B)\bigr).$$
By methods similar to those of \S 4.1, we can show that
$$gira(A,B)=gaxi\Bigl(\bigl(A,h(A)\bigr),\bigl(B,h(B)\bigr)\Bigr)\eqno(2.9.1)$$
with $h=push\cdot swap\cdot invmu\cdot swap$. 
\vskip .3cm
We define two operators on moulds following Ecalle ([Pisa, p. 49]):
$$ras\cdot B=invgari\cdot swap\cdot invgari\cdot swap(B)\eqno(2.9.2)$$
$$rash\cdot B=mu\bigl(push\cdot swap\cdot invmu\cdot swap(B),B\bigr).\eqno(2.9.3)$$
\vskip .3cm
\noindent {\bf Theorem 2.8.1.} {\it We have {\rm Ecalle's first fundamental
identity}:
$$gira(A,B)=ganit_{rash(B)}\cdot gari(A,ras\cdot B).\eqno(2.9.4)$$}
The remainder of this chapter is devoted to proving this theorem. Recall the
definitions of $gaxit$, $ganit$, $garit$, $gaxi$ and $gari$ from \S 2.7.
We use the (perhaps slightly doubtful) notation $invgaxi_{A,B}(A)$ to denote the
left-hand component of the pair $invgaxi(A,B)$.
\vskip .3cm
\noindent {\bf Lemma 2.8.2.} {\it We have
$$gaxit_{A,B}\cdot garit_{invgaxi_{A,B}(A)}=ganit_{mu(B,A)}.\eqno(2.9.5)$$\par}
\noindent {\bf Proof.} We have
$$garit_{invgaxi_{A,B}(A)}=gaxit_{invgaxi_{A,B}(A),invmu\cdot
invgaxi_{A,B}(A)},$$
and the composition of two $gaxits$ is given by
$$gaxit_{A,B}\cdot gaxit_{C,D}=gaxit_{gaxit_{A,B}(C)\,A,B\,gaxit_{A,B}(D)},
\eqno(2.9.6)$$
so we can multiply the terms on the LHS of (2.9.5) to obtain
$$gaxit_{gaxit_{A,B}(invgaxi_{A,B}(A))A,B\,gaxit_{A,B}(invmu\cdot
invgaxi_{A,B}(A))}.\eqno(2.9.7)$$
But we have
$$gaxit_{A,B}(invgaxi_{A,B}(A))=invmu\,A,\eqno(2.9.8)$$ 
since by definition of the $gaxi$-multiplication, we have
$$mu(gaxit_{A,B}(invgaxi_{A,B}(A)),A)=gaxi(invgaxi_{A,B}(A),A)=1.$$
Thus we can substitute (2.9.8) into (2.9.7) to obtain
$$gaxit_{1,B\,gaxit_{A,B}(invmu\cdot invgaxi_{A,B}(A))}.\eqno(2.9.9)$$
Similarly, by (2.9.8) and because $gaxit$ is an automorphism for $mu$,
we find that
$$\eqalign{gaxit_{A,B}\bigl(invmu\cdot invgaxi_{A,B}(A)\bigr)&=
invmu\Bigl(gaxit_{A,B}\bigl(invgaxi_{A,B}(A)\bigr)\Bigr)\cr
&= invmu\cdot invmu\cdot A\cr
&=A,}$$
and replacing this into (2.9.9) yields the desired result
$gaxit_{1,mu(B,A)}$, which is equal to $ganit_{mu(B,A)}$.
This concludes the proof of Lemma 2.8.2.\hfill{$\square$}
\vskip .3cm
Let $h=push\cdot swap\cdot invmu\cdot swap$ as in (2.9.1), and
let us introduce the notation $gaxit^h_B=gaxit_{(B,h(B))}$. We also write 
$gaxi^h(A,B)$ for the left-hand component of the pair 
$gaxi\bigl((A,h(A)),(B,h(B))\bigr)$, i.e. $gaxi^h(A,B)=mu\bigl(gaxit^h_B\cdot A,B
\bigr)$ by (2.8.9). Finally, we write $invgaxi^h(A)=invgaxi_{A,h(A)}(A)$, i.e.
the left-hand component of the $gaxi$-inverse of the pair $(A,h(A))$.
\vskip .3cm
\noindent {\bf Lemma 2.8.3.} {\it We have
$$invgaxi^h(B)=swap\cdot invgari\cdot swap\cdot B.\eqno(2.9.10)$$}
\noindent {\bf Proof.} We will show using (2.9.1) that the pair
$\bigl(swap\cdot invgari\cdot swap\cdot B,h(swap\cdot invgari\cdot swap\cdot B)\bigr)$
is the $gaxi$-inverse of $\bigl(B,h(B)\bigr)$. We have
$$\eqalign{gaxi&\Bigl(\bigl(swap\cdot invgari\cdot swap\cdot B,h(swap\cdot invgari\cdot swap\cdot B)\bigr),\bigl(B,h(B)\bigr)\Bigr)\cr
&=gira(swap\cdot invgari\cdot swap\cdot B,B)\cr
&=swap\Bigl(gari\bigl(invgari\cdot swap\cdot B,swap\cdot B\bigr)\Bigr)\cr
&=swap({\bf 1})\cr
&={\bf 1},}$$
where ${\bf 1}$ is the identity mould (that takes the value 1 on the empty set
and $0$ elsewhere). Thus $swap\cdot invgari\cdot swap\cdot B$ is indeed the
left-hand component of the $gaxi$-inverse of $(B,h(B))$, i.e. 
$invgaxi^h(B)$.\hfill{$\square$}
\vskip .3cm
\noindent {\bf Lemma 2.8.4.} {\it We have
$$\cases{gaxit_{A,B}\bigl(invgaxi_{A,B}(A)\bigr)=invmu\cdot A\cr
garit_C\bigl(invgari(C)\bigr)=invmu\cdot C\cr
gaxit^h_C\bigl(invgaxi^h(C)\bigr)=invmu\cdot C.}\eqno(2.9.11)$$}
\vskip .2cm
\noindent {\bf Proof.} Writing $garit_C=gaxit_{C,invmu\cdot C}$ and
$gaxit^h_C=gaxit_{C,h(C)}$ shows that the first equality implies the
second and third, so we only need to prove the first one.  To prove it,
we simply note that the left-hand component of 
$gaxi(invgaxi(A,B),(A,B))$ is the identity mould ${\bf 1}$, and it is given by 
$mu\bigl(gaxit_{A,B}(invgaxi_{A,B}(A)),A\bigr)$. This
proves the result.\hfill{$\square$}
\vskip .3cm
\noindent {\bf Lemma 2.8.5}. {\it We have
$$ganit_{rash\cdot C}(ras\cdot C)=C.\eqno(2.9.12)$$}
\vskip .2cm
\noindent {\bf Proof.} Recall that $rash\cdot B=mu(h(B),B)$. By (2.9.10) we have
$$ras\cdot B=invgari\cdot swap\cdot invgari\cdot swap\cdot B
=invgari\cdot invgaxi^h(B).\eqno(2.9.13)$$
Let us apply (2.9.5) with $A=C$ and $B=h(C)$, so that
$$gaxit_{C,h(C)}\cdot garit_{invgaxi_{C,h(C)}(C)}=ganit_{rash\cdot C}.
\eqno(2.9.14)$$
The LHS of (2.9.12) is the RHS of (2.9.14) applied to $ras\cdot C$, so to 
compute it, we will study the LHS of (2.9.14) applied to $ras\cdot C$.
Using the fact that $gaxit$ is a $mu$-automorphism, we obtain
$$\eqalign{gaxit^h_C\cdot &garit_{invgaxi^h(C)}\bigl(
invgari\cdot invgaxi^h(C)\bigr)\cr
&=gaxit^h_C\cdot invmu\cdot invgaxi^h(C)\ \ \ {\rm by}\ (2.9.11)\cr
&=invmu\cdot gaxit^h_C\cdot invgaxi^h(C)\cr
&=invmu\cdot invmu\cdot C\ \ \ {\rm by}\ (2.9.11)\cr
&=C.}$$
This completes the proof.\hfill{$\square$}
\vskip .3cm
We can now prove Theorem 2.8.1.  By (2.9.1) we have
$$gira(A,B)=gaxi^h(A,B).$$
With this, the desired (2.9.4) becomes
$$gaxi^h(A,B)=ganit_{rash\cdot B}\cdot gari(A,ras\cdot B).\eqno(2.9.15)$$
By (2.9.5), we have
$$gaxit_{A,B}\cdot garit_{invgaxi_{A,B}(A)}=ganit_{mu(B,A)}.$$
Replacing the couple $(A,B)$ by $(B,h(B))$ and recalling that
$mu(h(B),B)=rash\cdot B$, this gives
$$gaxit^h_B\cdot garit_{invgaxi^h(B)}=ganit_{rash\cdot B},$$
which, given that the inverse automorphism of $garit_B$ is $garit_{invgari(B)}$,
we can rewrite as 
$$gaxit^h_B =ganit_{rash\cdot B}\cdot garit_{invgari\cdot invgaxi^h(B)}
=ganit_{rash\cdot B}\cdot garit_{ras\cdot B}\eqno(2.9.16)$$
since by definition of $ras$ and (2.9.10) we have
$$ras\cdot B=invgari\cdot swap\cdot invgari\cdot swap\cdot B
=invgari\cdot invgaxi^h(B).$$ 
We will prove (2.9.4) by applying each side of (2.9.16) to a mould $A$,
then $mu$-multiplying the result with $B$.

The LHS of (2.9.16) yields 
$$mu\bigl(gaxit^h_B(A),B\bigr)=gaxi^h(A,B).$$
The RHS yields
$$\eqalign{
mu&\bigl(ganit_{rash\cdot B}\cdot garit_{ras\cdot B}(A),B\bigr)\cr
&=mu\bigl(ganit_{rash\cdot B}\cdot garit_{ras\cdot B}(A),
ganit_{rash\cdot B}(ras\cdot B)\bigr)\ \ {\rm by\ (2.9.12)}\cr 
&=ganit_{rash\cdot B}\cdot mu\bigl(garit_{ras\cdot B}(A),ras\cdot B)\bigr)\cr
&=ganit_{rash\cdot B}\cdot gari\bigl(A,ras\cdot B)\bigr).}$$
This completes the proof of Theorem 2.8.1.\hfill{$\square$}
\vskip .3cm
The following corollary of Theorem 2.8.1 containing the equality (2.9.17)
will be useful in Chapter 4, when we come
to prove Ecalle's second fundamental identity.  Let 
$fragari(A,B)=gari(A,invgari(B))$. Then (2.9.17) is proved simply by 
substituting $C=invgari\cdot ras\cdot B=
swap\cdot invgari\cdot swap\cdot B$ into (2.9.4).
\vskip .3cm
\noindent {\bf Corollary 2.8.6.} {\it We have
$$swap\cdot fragari(swap\cdot A,swap\cdot C)=ganit_{crash\cdot C}\cdot
fragari(A,C),
\eqno(2.9.17)$$
where $crash\cdot C=rash\cdot swap\cdot invgari\cdot swap\cdot C$.}
\vfill\eject
\centerline{\bf Chapter 3} 
\vskip .5cm
\centerline{\bf From double shuffle to \ARI}
\vskip .6cm
In this chapter, we define a map from the twisted Magnus Lie algebra
$\mt$ (introduced in \S 1.3) to $\ARI^{pol}_{al}$, and prove that it is a Lie
algebra isomorphism. We further show that the images of the two Lie subalgebras
$\ls$ and $\ds$ of $\mt$ defined in \S 1.3 and \S 1.4 map isomorphically
onto $\ARI^{pol}_{\underline{al}*\underline{al}}$ and
$\ARI^{pol}_{al*il}$.
In \S 3.4 we use the results of Chapter 2 together with these isomorphisms
to show how Ecalle's methods give a simple proof of some basic results on
double shuffle (Theorems 1.3.2 and 1.4.1), namely that
$\ls_n^d$ is zero if $n\not\equiv d$ mod 2, and hence also
$\ds_n^d/\ds_n^{d+1}$ is zero if $n\not\equiv d$ mod 2.
\vskip .5cm
\noindent {\bf \S 3.1. The ring $\F$}
\vskip .2cm
Consider the ring of polynomials $\Q\langle x,y\rangle$ in non-commutative 
variables $x,y$.  Let $\partial_x$ denote the differential operator with 
respect to $x$.  Set $C_i={\rm ad}(x)^{i-1}(y)$, $i\ge 1$, so
$C_1=y$, $C_2=[x,y]$, $C_3=[x,[x,y]],\ldots$.
\vskip .3cm
\noindent {\bf Definition.} Let $\F$ denote the subspace of $\Q\langle x,y
\rangle$ of polynomials $f$ such that $\partial_x(f)=0$.  
\vskip .2cm
The following well-known result is 
just a standard application of Lazard elimination.
\vskip .2cm
\noindent {\bf Lemma 3.1.1.} {\it The subspace $\F\subset
\Q\langle x,y \rangle$ is equal to the subring generated by the $C_i$, $i\ge 1$.
Moreover the $C_i$ are free generators of this ring. }
\vskip .3cm
Let $\pi_y$ be the projector onto polynomials ending in $y$ (i.e. $\pi_y$ 
forgets all the monomials ending in $x$).  The usefulness of the ring $\F$
is that $\pi_y$ has a section on $\F$.
Indeed, for any polynomial $g$ ending in $y$, define $\sec(g)$ by
$$\sec(g)=\sum_{i\ge 0} {{(-1)^i}\over{i!}}\partial_x^i(g)x^i.\eqno(3.1.1)$$
\vskip .1cm
\noindent {\bf Lemma 3.1.2.} [R, Prop IV.2.8] {\it (1) $\sec\circ\pi_y={\rm id}$ on $\F$.
\vskip .1cm
(2) $\pi_y\circ \sec={\rm id}$ on $\Q\langle x,y\rangle y$.}
\vskip .5cm
\noindent {\bf \S 3.2. Associating moulds to elements $f\in \F$}
\vskip .2cm
\noindent {\bf Definitions.}
Let $\F_n$ denote the vector subspace of polynomials in $\F$ of homogeneous
degree $n$ in $x$ and $y$, $\F^r$ the subspace of polynomials of
homogeneous degree $r$ (i.e. linear combinations of monomials of the
form $C_{a_1}\cdots C_{a_r}$), and $\F_n^r$ the intersection.  The
space $\F$ is bigraded, i.e. $\F=\oplus_{n,r\ge 0} \Q\langle C\rangle_n^r$.
If $f\in \F$, we write $f_n$ for its weight $n$ part and $f^r$ for its
depth $r$ part.

Let $\pi_y(f)$ denote the projection of $f$ onto the monomials ending in $y$
as above, and let $f_y$ denote $\pi_y(f)$ rewritten in 
the variables $y_i=x^{i-1}y$, $i\ge 1$, and $f_y^r$ the depth 
$r$ part, i.e. $\pi_y(f^r)$ written in the $y_i$.  Similarly, let
$\pi_Y(f)$ denote the projection of $f$ onto the monomials starting
with $y$.
Let $\ret_X:\Q\langle x,y\rangle\rightarrow \Q\langle x,y\rangle$ 
denote the ``backwards writing'' map 
$$\ret_X(x^{a_0}y\cdots yx^{a_{r-1}}yx^{a_r})=
x^{a_r}yx^{a_{r-1}}y\cdots yx^{a_0}.\eqno(3.2.1)$$
Note that ${\rm Lie}[x,y]\subset \F$.  If $f\in \F_n$ is actually a Lie 
element, we have 
$$\ret_X(f)=(-1)^{n-1}f.\eqno(3.2.2)$$
Finally, let $f_Y^r$ denote the polynomial $\ret_X\bigl(\pi_Y(f^r))$ written
in the variables $y_i$ and $f_Y=\sum_r f_Y^r$.
\vskip .4cm
We note here that by a result in [CS], the introduction of $f_Y$ 
gives an equivalent formulation of the definition of $\ds$ that will be useful below.
\vskip .3cm
We saw in Lemma 3.1.1 that $\F$ is the set of polynomials in $\Q\langle x,y\rangle$
that can be written as polynomials in the $C_i$, and that such a
writing is unique.  Let $f_C$ denote $f$ written in this way.

Define three maps from monomials in the variables $x,y$ (resp.
$y_1,y_2,\ldots$ resp. $C_1,C_2,\ldots$) to monomials in commutative variables 
$z_0,z_1,\ldots$ (resp.  $u_1,u_2,\ldots$ resp. $v_1,v_2,\ldots$) as follows:
$$\eqalign{\iota_X&:x^{a_0-1}y\cdots x^{a_{r-1}-1}yx^{a_r-1} \mapsto 
z_0^{a_0-1}\cdots z_r^{a_r-1} \cr
\iota_C&:C_{a_1}\cdots C_{a_r}\mapsto u_1^{a_1-1}\cdots u_r^{a_r-1} \cr
\iota_Y&:y_{a_1}\ldots y_{a_r}\mapsto v_1^{a_r-1}\cdots v_r^{a_1-1}.}
\eqno(3.2.3)$$ 

Then we define a mould in commutative variables $z_0,z_1,\ldots$ associated to $f\in \F_n$ 
as follows:
$$vimo_f(z_0,z_1,\ldots,z_r)=\iota_X(f^r),\eqno(3.2.4)$$
and also a mould and a $v$-mould associated to $f$ by
$$ma_f(u_1,\ldots,u_r)=(-1)^{r+n}\iota_C(f_C^r),\ \ 
mi_f(v_1,\ldots,v_r)=\iota_Y\bigl(f_Y^r).\eqno(3.2.5)$$
All other values of these moulds are $0$.
\vskip .3cm
\noindent {\bf Remark.} Note that by (3.2.2), if $f\in {\rm Lie}[x,y]$, we have 
$$\pi_y(f)=(-1)^{n-1}\ret_X\bigl(\pi_Y(f)\bigr),$$
so $f_y^r=(-1)^{n-1}f_Y^r$.
Thus, if $f\in {\rm Lie}[x,y]$, the $v$-mould $mi$ can also be defined by
$$mi_f(v_1,\ldots,v_r)=(-1)^{n-1}\iota_Y(f_y^r).\eqno(3.2.6)$$
When we turn our attention to the twisted Magnus Lie algebra $\mt$ and its
double shuffle and linearized double shuffle subspaces, in \S\S 3.3-3.4,
we will be in this situation.
\vskip .3cm
Since the maps $\iota_X$, $\iota_C$ and $\iota_Y$ are obviously invertible,
we recover $f$ from $vimo_f$, $f_C$ from $ma$ and $f_Y$ from $mi$.  But of
course, we easily recover $f$ from $f_C$ by expanding out the $C_i$, and
we also recover $f$ from $f_Y$ by setting $f={\rm sec}(f_Y)$, as we have 
assumed that $f\in \F_n$.  Thus, for any element $f\in \F_n$, $f$ itself, 
$f_C$, $f_Y$ and $vimo_f$ are all different encodings of the same information.
The moulds $ma$ and $mi$ are also equivalent encodings, related to $vimo_f$ as follows.
\vskip .3cm
\noindent {\bf Lemma 3.2.1.} {\it The mould $ma$ and the $v$-mould $mi$ are obtained
from $vimo_f$ by the formulas
$$ma_f(u_1,\ldots,u_r)=vimo_f(0,u_1,u_1+u_2,\ldots,
u_1+\cdots+u_r)\eqno(3.2.7)$$
$$mi_f(v_1,\ldots,v_r)=vimo_f(0,v_r,v_{r-1},\ldots,
v_1).\eqno(3.2.8)$$}
\vskip .1cm
The proof of this lemma is given in \S A.3 of the Appendix.
\vskip .3cm
\noindent {\bf Remark.} If $vimo(z_0,\ldots,z_r)$ for $r\ge 0$ is an arbitrary family
of polynomials, then there is a unique $f\in \Q\langle x,y\rangle$ associated to it 
by (3.2.3).  It is natural to ask what condition on the family $vimo$ ensures that
$f\in \F$. We leave the following answer as an exercise.
\vskip .3cm
\noindent {\bf Lemma 3.2.2.} {\it If $f\in \Q\langle x,y\rangle$ 
and $vimo_f$ is defined as in (3.2.4), then $f\in \F$ if and only if
$$vimo_f(z_0,\ldots,z_r)=vimo_f(0,z_1-z_0,z_2-z_0,\ldots,z_r-z_0)\eqno(3.2.9)$$
for $r\ge 1$.}
\vskip .3cm
\noindent {\bf Remarks.} (1) Observe that if we apply the variable change $u_1=z_1-z_0$,
$u_2=z_2-z_1$, $u_3=z_3-z_2,\ldots,u_r=z_r-z_{r-1}$ to $ma_f(u_1,\ldots,u_r)$,
obtaining $vimo_f(0,z_1-z_0,\ldots,z_r-z_0)$. Thanks to (3.2.9), if $f\in\F$
then this is equal to $vimo_f(z_0,\ldots,z_r)$, so that $ma_f$ is yet another
equivalent coding for $f\in\F$, and the same holds for $mi_f$ using the
variable change $v_j=z_{r-j+1}-z_0$.
\vskip .2cm
(2) From the expressions (3.2.7) and (3.2.8), it is immediate that for 
$f\in \F$, we have
$$swap(ma_f)=mi_f.\eqno(3.2.10)$$
\vskip .6cm
\noindent {\bf Example.} Let $f$ be the degree 3 Lie polynomial
$$f=[x,[x,y]]+[[x,y],y]=x^2y-2xyx+yx^2+xy^2-2yxy+y^2x.\eqno(3.2.11)$$
Then $\pi_y(f)=x^2y-2yxy+xy^2$, $f_Y=y_3-2y_1y_2+y_2y_1$ and
$f_C=C_3-C_1C_2+C_2C_1$, and we have
$$\cases{vimo_f(z_0)=0\cr
vimo_f(z_0,z_1)=z_0^2-2z_0z_1+z_1^2\cr
vimo_f(z_0,z_1,z_2)=z_0-2z_1+z_2\cr
vimo_f(z_0,z_1,z_2,z_3)=0,}\ 
\cases{ma_f(\emptyset) = 0\cr
ma_f(u_1) = u_1^2\cr
ma_f(u_1,u_2) = -u_1+u_2\cr
ma_f(u_1,u_2,u_3) = 0,}\ 
\cases{mi_f(\emptyset) = 0\cr
mi_f(v_1) = v_1^2\cr
mi_f(v_1,v_2) = -2v_2+v_1\cr
mi_f(v_1,v_2,v_3) = 0.}$$
The results of this section can be summarized by the following theorem. We write $\Fzero$ for the subspace of polynomials in $\F$ with constant term $0$.
\vskip .3cm
\noindent {\bf Theorem 3.2.3.} {\it Let $\Fzeroo$ denote the
degree completion of the polynomial space $\Fzero$, consisting
of power series in the $c_i$ with constant term $0$. Then the map 
$$ma:\Fzeroo\rightarrow \ARI^{pol}\eqno(3.2.12)$$ 
is
a ring isomorphism, where $\Fzeroo$ is equipped with the ordinary (concatenation)
multiplication of polynomials, and $\ARI^{pol}$ with the multiplication $mu$.}
\vskip .2cm
\noindent {\bf Proof.} By Lemma 3.1.1 together with the definition of $\iota_C$ 
in (3.2.3) and the definition of the map $ma$ in (3.2.5),
we see that $ma$ is a vector space isomorphism from $\Fo$ to the set of
polynomial-valued moulds, so it restricts from $\Fzeroo$ to $\ARI^{pol}$.
Thus it remains only to show that 
$$ma_{fg}=mu(ma_f,ma_g).\eqno(3.2.13)$$  
By additivity,
it is enough to assume that $f$ and $g$ are monomials in the $C_i$, say
$f=C_{a_1}\cdots C_{a_r}$ and $g=C_{b_1}\cdots C_{b_s}$; then it is 
immediate that
$$ma_{fg}=u_1^{a_1-1}\cdots u_r^{a_r-1}u_{r+1}^{b_1-1}\cdots u_{r+s}^{b_s-1}
=mu(ma_f,ma_g).$$
This concludes the proof.\hfill{$\square$}
\vskip .8cm
\noindent {\bf \S 3.3.  The Poisson bracket and the \ARI bracket}
\vskip .4cm
In this section we prove that the Poisson bracket is carried over to the
$ari$-bracket under the isomorphism $\F\buildrel{ma}\over
\rightarrow \ARI^{pol}$ of (3.2.12). This result was originally proved in [R, Appendice A,
\S 5]. After introducing the key result in Lemma 3.3.1 (due to Racinet),
we then compare the derivations $D_f$ and $arit(ma_f)$ in Proposition
3.3.3 and deduce the equality $ma_{\{f,g\}}=ari(ma_f,ma_g)$ in Corollary
3.3.4.
\vskip .2cm
Observe that if $f\in \F_n$, then 
$\partial_x([x,f])=0$, so by Lemma 3.1.1,
$[x,f]\in \F_{n+1}$.  By Lemma 3.1.1, we can consider both $f$ and $[x,f]$ as
being polynomials in the $C_i$.  
\vskip .4cm
\noindent {\bf Lemma 3.3.1.} [R] {\it Let $f\in \F_n$.  Then for 
$0\le r\le n$, we have 
$$ma_{[x,f^r]}=-(u_1+\cdots+u_r)ma_{f^r}.\eqno(3.3.1)$$}
\noindent {\bf Proof.}  
Note first that $a\mapsto [x,a]$ is a 
derivation, i.e. $[x,ab]=[x,a]b+a[x,b]$.  Thus,
writing $f^r=\sum_{\a} c_\a C_{a_1}\cdots C_{a_r}$, where $\a=(a_1,\ldots,a_r)$, 
we have
$$\eqalign{[x,f^r]=\sum_{\a}c_\a [x,C_{a_1}\cdots C_{a_r}]&=
\sum_{\a}c_\a\sum_{i=1}^r C_{a_1}\cdots C_{a_{i-1}}[x,C_{a_i}]
C_{a_{i+1}}\cdots C_{a_r}\cr
&=\sum_{\a}\sum_{i=1}^r c_\a\ C_{a_1}\cdots C_{a_{i-1}}C_{a_i+1}
C_{a_{i+1}}\cdots C_{a_r}.}$$
Thus, the left-hand side of (3.3.1) is equal to
$$(-1)^{r+n+1}\sum_{\a} \sum_{i=1}^r c_\a\  u_1^{a_1-1}\cdots 
u_i^{a_i}\cdots u_r^{a_r-1}.\eqno(3.3.2)$$
But since $ma_{f^r}=(-1)^{r+n}\sum_{\a} c_\a\  u_1^{a_1-1}\cdots u_r^{a_r-1}$,
(3.3.2) is equal to $ma_{f^r}$ multiplied by $-(u_1+\cdots+u_r)$, proving
(3.3.1).\hfill{$\square$}
\vskip .4cm
\noindent {\bf Proposition 3.3.2.} {\it For any mould $A$, the operator 
$arit(A)$ is a derivation for the $mu$-multiplication.}
\vskip .3cm
The proof is given in the Appendix, \S A.4.
\vskip .3cm
\noindent {\bf Proposition 3.3.3.} {\it Let $f\in\F_n$ be of homogeneous
depth $r$ and $g\in \F_m$ of homogeneous depth $s$.  Let 
$D_f$ be the derivation of $\F$ defined by $D_f(x)=0$, $D_f(y)=[y,f]$.  Then
$$ma_{D_f(g)}=-arit(ma_f)\cdot ma_g.\eqno(3.3.3)$$}
\noindent {\bf Proof.}  We have $D_{f+g}=D_f+D_g$, so we may assume that
$f=C_{a_1}\cdots C_{a_r}$ is a monomial in the $C_i$.  Furthermore, a
derivation of $\F$ is defined by its action on the generators $C_i$, 
so we may take $g=C_m={\rm ad}(x)^{m-1}(y)$.  
Let $F_0=[y,f]$, and for $i\ge 1$, let 
$F_i={\rm ad}(x)^i([y,f])$. In particular, we have
$$D_f(g)=[x,[x,\cdots,[x,[y,f]]\cdots]={\rm ad}(x)^{m-1}([y,f])
=F_{m-1}.$$
Then by Lemma 3.3.1, since all the $F_i$ are in depth $r+1$, we have
$$ma_{F_i}=-(u_1+\ldots+u_{r+1})ma_{F_{i-1}}\ \ {\rm for}\ \
i>0,$$
so
$$ma_{F_i}=(-1)^i(u_1+\ldots+u_{r+1})^ima_{F_0},$$
so the left-hand side of (3.3.3) is equal to
$$\eqalign{ma_{D_f(g)}&=ma_{F_{m-1}}\cr
&=(-1)^{m-1}(u_1+\ldots+u_{r+1})^{m-1}ma_{F_0}\cr
&=(-1)^{m-1}(u_1+\ldots+u_{r+1})^{m-1}ma_{[y,f]}\cr
&=(-1)^{m+r+n+1}(u_1+\cdots+u_{r+1})^{m-1}\iota_C(C_1C_{a_1}\cdots C_{a_r}-C_{a_1}\cdots 
C_{a_r}C_1)\cr
&=(-1)^{m+r+n+1}(u_1+\cdots+u_{r+1})^{m-1}  \bigl(u_2^{a_1-1}\cdots u_{r+1}^{a_r-1}-
u_1^{a_1-1}\cdots u_r^{a_r-1}\bigr),}\eqno(3.3.4)$$
since $ma_{[y,f]}=(-1)^{r+n+2}\iota_C([y,f])$.
Now consider the right-hand side of (3.3.3).  By (3.2.5), we have
$$ma_f(u_1,\ldots,u_r)=(-1)^{r+n}\iota_C(f)=
(-1)^{r+n}u_1^{a_1-1}\cdots u_r^{a_r-1},$$
where $n=a_1+\cdots +a_r$, and 
$$ma_g(u_1)=(-1)^{m-1}u_1^{m-1}.$$
Since $ma_g$ has value zero on any word of length greater than 1, the
defining formula for $arit(A)$ simplifies when $A=ma_f$, $B=ma_g$ to
$$\eqalign{\bigl(arit(ma_f)\cdot ma_g\bigr)(u_1,\ldots,u_r))&=
ma_g(u_1)ma_f(u_2,\ldots,u_{r+1})\cr
&+ma_g(u_1+\cdots+u_{r+1}) ma_f(u_1,\cdots,u_r)\cr
&-ma_g(u_1+\cdots+u_{r+1}) ma_f(u_2,\ldots,u_{r+1})\cr
&-ma_g(u_1)ma_f(u_2,\ldots,u_{r+1})}$$
$$\qquad\qquad\qquad\qquad\qquad
=ma_g(u_1+\cdots+u_{r+1}) \Bigl(ma_f(u_1,\cdots,u_r)
-ma_f(u_2,\ldots,u_{r+1})\Bigr)$$
$$\qquad\qquad\qquad\qquad\qquad\ \ \ 
=(-1)^{m+r+n}(u_1+\cdots+u_{r+1})^{m-1}\bigl(u_2^{a_1-1}\cdots
u_{r+1}^{a_r-1}-u_1^{a_1-1}\cdots u_r^{a_r-1}\bigr).$$
This proves (3.3.3).\hfill{$\square$}
\vskip .3cm
\noindent {\bf Corollary 3.3.4.} {\it Let $f\in\F_n$ be of homogeneous
depth $r$ and $g\in \F_m$ of homogeneous depth $s$.  Then
$$ma_{\{f,g\}}=ari(ma_f,ma_g).\eqno(3.3.5)$$}
\noindent {\bf Proof.} Recall that $\{f,g\}=D_f(g)-D_g(f)+fg-gf.$
By (3.3.3) and (3.2.13), we then have
$$\eqalign{ma_{\{f,g\}}&=-arit(ma_f)\cdot ma_g+arit(ma_g)\cdot ma_f+mu(ma_f,ma_g)-mu(ma_g,ma_f)\cr
&=arit(ma_g)\cdot ma_f-arit(ma_f)\cdot ma_g+lu(ma_f,ma_g)\cr
&=ari(ma_f,ma_g)}$$
by (2.2.9). This concludes the proof.\hfill{$\square$}
\vskip .8cm
\noindent {\bf \S 3.4. The $ma$ map from $\ds$ to \ARI}
\vskip .3cm
In this section we relate the special Lie subspaces $\mt$, $\ls$ and $\ds$
of $\F$ to some of the special subspaces of $\ARI$ defined in \S 2.5. The
proofs are based on the following explicit comparison of double shuffle 
properties of polynomials in $\F$ with symmetry properties on moulds.
\vskip .3cm
\noindent {\bf Lemma 3.4.1.} {\it Let $f\in \F_n$.  Then
\vskip .1cm
(i) $f$ satisfies shuffle in $x,y$ if and only if $ma_f\in \ARI^{pol}_{al}$;
\vskip .1cm
(ii) $f_Y$ satisfies shuffle in the $y_i$ if and only if $mi_f\in \overline{\ARI}^{pol}_{al}$;
\vskip .1cm
(iii) $f_Y$ satisfies stuffle in the $y_i$ if and only if $mi_f\in \overline{\ARI}^{pol}_{il}$;
\vskip .1cm
(iv) $f_Y$ satisfies stuffle in the $y_i$ in depth $1\le r<n$ if and only if 
$mi_f\in \overline{\ARI}^{pol}_{*il}$.}
\vskip .2cm
This Lemma is proved in the Appendix, \S A.5.
\vskip .3cm
\noindent {\bf Theorem 3.4.2.} {\it The isomorphism $ma:\F\buildrel\sim\over
\rightarrow \ARI^{pol}$ restricts to an isomorphism of Lie algebras
$$ma:\mt\buildrel\sim\over\rightarrow \ARI^{pol}_{al}.\eqno(3.4.1)$$}
\noindent {\bf Proof.} We first observe that $\mt\subset \F$ since by
definition, the underlying vector space of $\mt$ is the Lie 
algebra freely generated by the $C_i$, $i\ge 1$ (see \S 1.3).  Indeed, we have
${\rm Lie}[x,y]\cap \F=\mt$.

Since $ma$ is injective on 
$\F$, it is injective restricted to $\mt$. By \S 1.3 (1), a polynomial
$f\in\F$ satisfies shuffle if and only if $f\in {\rm Lie}[x,y]$, which shows
that every $f\in\mt$ satisfies shuffle. Then Lemma 3.4.1 (i) shows that
$ma_f\in \ARI^{pol}_{al}$. Conversely, if $A\in \ARI^{pol}_{al}$, then
since $ma:\F\rightarrow \ARI^{pol}$ is an isomorphism, there exists
$f\in\F$ such that $A=ma_f$, and then again by Lemma 3.4.1 (i), $f$
must satisfy shuffle, i.e. $f\in {\rm Lie}[x,y]\cap \F=\mt$.\hfill{$\square$}
\vskip .3cm
We can now proceed to the first main result of this section.
\vskip .3cm
\noindent {\bf Theorem 3.4.3.} {\it The map $f\mapsto ma_f$ yields 
a Lie algebra isomorphism
$$\ls\buildrel\sim\over\rightarrow \ARI^{pol}_{\underline{al}/
\underline{al}}.\eqno(3.4.2)$$}
\noindent {\bf Proof.} Thanks to (3.3.5), which shows that the Poisson
bracket on $\mt$ carries over to the $ari$-bracket, it suffices to show that
(3.4.2) is a vector space isomorphism.  Let $f\in \ls$; we
may assume that $f$ is homogeneous of degree $n$. Recall that the definition
of $\ls$ is that $f$ must satisfy shuffle in $x,y$ and $\pi_y(f)$ must satisfy 
shuffle in the $y_i$.  Since $f$ is a Lie polynomial, we have
$f_Y=ret_X(\pi_Y(f))=(-1)^{n-1}\pi_y(f)$, so $f_Y$ satisfies shuffle in the
 $y_i$ if and only if $\pi_y(f)$ (rewritten in the $y_i$) does.  But 
by Lemma 3.4.1 (ii), $f_Y$ satisfies the shuffle in the $y_i$ if and only if 
$mi_f\in\overline{\ARI}_{al}$. Thus the image of $\ls$ under the injective map $f\mapsto ma_f$
lies in $\ARI^{pol}_{al/al}$.  Recall that by the definition of $\ls$ (see
\S 1.4), the even degree depth 1 polynomials $ad(x)^{2i+1}(y)$ are excluded 
from $\ls$; thus the image of $\ls$ lies in 
$\ARI^{pol}_{\underline{al}/\underline{al}}$. 

Conversely, if $A\in\ARI^{pol}_{\underline{al}/\underline{al}}$, then
since $ma:\F\rightarrow \ARI^{pol}$ is an isomorphism, there exists a unique
$f\in\F$ such that $ma_f=A$, and then by Lemma 3.4.1, $f$ must satisfy
shuffle and $f_Y$ must satisfy shuffle in the $y_i$, and if
$ma_f$ is of depth $1$ then $f$ is of odd degree, so $f\in \ls$.\hfill{$\square$}
\vskip .3cm
From this we deduce the proofs of Theorem 1.4.1 (which then implies
Theorem 1.3.4), which is essentially no more than a translation back into
$\F$ of Theorem 2.5.6 stating that $\ARI_{al/al}$ is a Lie algebra under the
$ari$-bracket.
\vskip .3cm
\noindent {\bf Corollary 3.4.4.} {\it The weight $n$, depth $d$ space
$\ls_n^d$ is zero if $n\not\equiv d$ mod 2; thus in particular the
graded quotient $\ds_n^d/\ds_n^{d+1}$ which lies inside it is zero
if $n\not\equiv d$ mod 2.}
\vskip .2cm
\noindent {\bf Proof.} Using the translation into moulds (3.4.2),
the statement is equivalent to the fact that if $A\in\ARI^{pol}_{\underline{al}/
\underline{al}}$ is a homogeneous polynomial mould $A(u_1,\ldots,u_d)$ of
of odd degree $n-d$, then $A=0$. But this follows immediately from Lemma 2.5.5 
which says that elements of $\ARI_{\underline{al},\underline{al}}$ are 
$neg$-invariant, i.e. $A(u_1,\ldots,u_d)=A(-u_1,\ldots,-u_d)$; indeed if
$A$ is homogeneous of odd degree, then $A$ must be zero.
\hfill{$\square$}
\vskip .3cm
This proof, or rather the proof of Lemma 2.5.5, is a perfect example of the 
real simplicity and magic of Ecalle's methods.
\vskip .3cm
Our next step is to prove the analogue of (3.4.2) for $\ds$.  We first need
a lemma that slightly rephrases the definition of $\ds$.
\vskip .3cm
\noindent {\bf Lemma 3.4.5.} {\it The Lie algebra $\ds$ is equal to the set of
$f\in {\rm Lie}[x,y]$ of degree $\ge 3$ such that $f_Y$, rewritten in the variables
$y_i$, satisfies all the stuffle relations (1.3.3) except for those where
both words in the pair $(u,v)$ are powers of $y$.}
\vskip .2cm
\noindent {\bf Proof.} Let the depth of a stuffle relation as in (1.3.3) be 
equal to the sum of the depths of the two words $(u,v)$. Let $f\in\ds$; we may
assume that $f$ is homogeneous of degree $n$. Suppose that $f_Y$ satisfies
all the stuffle relations of depths $<n$. Since $f$ is Lie, we have
$ret_X(f)=(-1)^{n-1}f$, so in particular 
$$f_Y=ret_X\bigl(\pi_Y(f)\bigr)=(-1)^{n-1}\pi_y(f);$$ 
thus $\pi_y(f)$ satisfies the same stuffle relations. Then [CS, Theorem 2] shows 
that there exists a unique constant, namely $a={{(-1)^{n-1}}\over{n}}(\pi_y(f)|
x^{n-1}y)$, 
such that $\pi_y(f)+ay^n$, rewritten in the $y_i$, satisfies all of the stuffle 
relations.  But the term $ay^n$ is equal to $f_{corr}$ as in (1.3.2), so this is 
equivalent to the original definition of $\ds$ given in \S 1.3.\hfill{$\square$}
\vskip .3cm
\noindent {\bf Theorem 3.4.4.} {\it The isomorphism $ma:\F\rightarrow\ARI^{pol}$
restricts to a Lie algebra isomorphism
$$\ds\buildrel\sim\over\rightarrow 
\ARI^{pol}_{\underline{al}*\underline{il}}\eqno(3.4.3)$$}
\noindent {\bf Proof.} We saw above that $f\mapsto ma_f$ maps $\ds$ injectively into $\ARI^{pol}_{al}$. 
Let $f\in\ds$, and assume that $f$ is homogeneous of degree $n$. Then as in the
proof of Lemma 3.4.5, $\pi_y(f)$ satisfies all the stuffle relations of depth 
$<n$, and $f_*=\pi_y(f)+ay^n$ satisfies all the stuffle relation, where
$a={{(-1)^{n-1}}\over{n}}(\pi_y(f)|x^{n-1}y)$. 

Now, let $mi_f=\iota_Y(f_Y)$ as in (3.2.5), and let $mi'_f=\iota_Y(f_*)$.
Then since $f_*$ satisfies the stuffle relations, by Lemma 3.4.1 (iii) we know that
$mi'_f$ is alternil.  But since (apart from the sign) $f_Y$ differs from $f_*$ 
only by the depth $n$ term $ay^n$, the two moulds $mi_f$ and $mi'_f$ differ
(up to sign) only by the depth $n$ component, which is a constant due to the
homogeneity of $f$, which in terms of moulds means that each 
$mi_f(v_1,\ldots,v_r)$ is a polynomial
of degree $n-r$. This means that it suffices to modify $mi_f$ by a constant in
depth $n$ to make it fully alternil, which is the definition of $\overline{\ARI}_{*il}$.
Thus $ma_f\in\ARI^{pol}_{al*il}$.  The surjectivity holds as before, since
surjectivity of $ma$ means that there exists a polynomial in $\F$ such that
$ma_f=A$ for any $A\in \ARI^{pol}_{al*il}$, and then by Lemma 3.4.1, $f$
must satisfy shuffle and $f_Y$ stuffle for depths $<n$; then using
Lemma 3.4.5 proves that $f\in\ds$.\hfill{$\square$}
\vskip .3cm
\noindent {\bf Example.}  We take the same example as in (3.2.11), and check that
$ma_f/mi_f$ is $al*il$ (i.e. $ma_f\in\ARI_{al}$ and $mi_f\in \overline{\ARI}_{*il}$).  Recall 
that 
$$\cases{ma_f(u_1)=u_1^2\cr
ma_f(u_1,u_2)=-u_1+u_2,}\ \ \ \ \ \ 
\cases{mi_f(v_1)=v_1^2\cr
mi_f(v_1,v_2)=v_1-2v_2.}$$
To show that $ma_f$ is alternal, the only condition to check
is that $ma_f(u_1,u_2)+ma_f(u_2,u_1)=0$, which is immediate.
To show that $mi_f$ is alternil, we only have to check the alternility relation
corresponding to the stuffle relation for depth $r=2$, given in (2.3.5):
$$(v_1-2v_2)+(v_2-2v_1)+{{1}\over{v_1-v_2}}v_1^2+{{1}\over{v_2-v_1}}v_2^2=
(-v_1-v_2)+(v_1+v_2)=0.$$
\vskip .8cm
\noindent {\bf \S 3.5. The group G\ARI and the twisted Magnus group.}
\vskip .5cm
In this section we establish the isomorphism between the twisted Magnus
group (defined below) and $G\ARI^{pol}_{as}$ which is the group analog of
Theorem 3.4.2. The proof is basically a corollary of Theorem 3.4.2 using
the exponential, but it is useful to recall the objects and definitions
that are the translations of $G\ARI^{pol}$ and its associated operators 
($ganit$, $garit$, $gari$ etc.) so as to clarify the fact that in this
familiar context they are in fact familiar operators, on the one hand, and
to emphasize the power of Ecalle's theory in extending from polynomial-valued
moulds to rational-valued moulds on the other.  We end the section by 
explaining the meaning of some of the main identities from \S 2.7 in the
twisted Magnus situation.
\vskip .2cm
\noindent {\bf Definition.} Let $f,g\in\mt$, and define $p(f,g)=fg-D_g(f)$
to be the {\it pre-Lie law} associated to $\mt$.  Obviously
$p(f,g)-p(g,f)=\{f,g\}$, and thanks to (3.3.3), we have
$$ma_{p(f,g)}=mu(ma_f,ma_g)+arit(ma_g)\cdot ma_f=preari(ma_f,ma_g).\eqno(3.5.1)$$
The expression $p(f,g)=fg-D_g(f)$ actually expresses the 
multiplication rule on the universal enveloping algebra ${\cal U}\mt$
for all $g\in \mt$, $f\in {\cal U}\mt$, not only when $f\in \mt$.

Define the twisted Magnus exponential on $\mt$ by 
$$exp^{\odot}(f)=1+f+\sum_{n\ge 2} {{1}\over{n!}} p(f^n),\eqno(3.5.2)$$
Then by (2.6.1) we have
$$ma_{exp^{\odot}(f)}=exp_{ari}(ma_f).\eqno(3.5.3)$$
where $p(f^n)=p(p(f^{n-1}),f)$, $p(f^3)=p(p(f,f),f)$ etc.

The {\it twisted Magnus group} $MT$ is the pro-unipotent group $exp^\odot(\mt)$.

By the Milnor-Moore theorem, we have an isomorphism of vector spaces
$${\cal U}\mt\simeq \F\eqno(3.5.4),$$
where both sides are Hopf algebras with the multiplication on the right-hand
ring being different than the usual concatenation, but the coproduct being
the restriction to $\F$ of the standard coproduct defined by
$$\Delta(C_i)=C_i\otimes 1+1\otimes C_i,\ \ i\ge 1.\eqno(3.5.5)$$
Indeed, the primitive elements of $\F$ for $\Delta$ are well-known to be the
Lie polynomials in the $C_i$, which form the underlying vector space ${\bf L}$ 
of $\mt$ (see \S 1.3). Since the ring $\F$ is a graded polynomial ring
(where the grading can be considered to the be degree in $x,y$ or else
the weight in the $C_i$ where each $C_i$ is of weight $i$) with
$\F_0=\Q$ and each graded part is finite-dimensional,
Milnor-Moore applies and yields the isomorphism (3.5.3).

As in the general case of Lie algebras, we have the inclusion of the
exponential group into the completion of the enveloping algebra, namely
$$exp^\odot(\mt)\subset \widehat{{\cal U}\mt}\simeq \widehat\F,\eqno(3.5.6)$$
where the right-hand ring is included (as vector spaces) in the power series
ring on $x$ and $y$. 

The group $exp^\odot(\mt)$ consists of the power series 
in $x,y$ that have constant term $1$ and no linear term in $x$, and are 
{\it group-like}, i.e. such that
$$\Delta(f)=f\otimes f.\eqno(3.5.7)$$
The expression for product of two elements of the subgroup $exp^\odot(\mt)$
is the {\it twisted Magnus multiplication law}
$$f(x,y)\odot g(x,y)=f(x,gyg^{-1})g(x,y).\eqno(3.5.8)$$
This multiplication corresponds to identifying $f\in exp^\odot(\mt)$ with
the endomorphism $R_f$ of $\Q\langle\langle x,y \rangle\rangle$ given by
$x\mapsto x$, $y\mapsto fyf^{-1}$. The twisted Magnus multiplication then
simply corresponds to anticomposition of endomorphisms; indeed, we have
$$R_g\circ R_f(y)=R_g(fyf^{-1})=f(x,gyg^{-1})gyg^{-1}f(x,gyg^{-1}),$$
so
$$R_g\circ R_f=R_{f(x,gyg^{-1})g}=R_{f\odot g}.\eqno(3.5.9)$$
We have
$$garit(ma_g)\cdot ma_f=ma_{R_g(f)}\eqno(3.5.10)$$
and
$$gari(ma_f,ma_g)=f\odot g.\eqno(3.5.11)$$
The group $MT$ is the set of all group-like power series in 
$\widehat\F$ with constant term 1, equipped with the twisted Magnus 
multiplication $\odot$ given in (3.5.8).  Let $\widehat\F_1$ denote the
set of all power series in $\widehat\F$ with constant term 1, equipped with
the multiplication $\odot$ of (3.5.8). Then (3.5.11) shows that $ma$ gives 
rise to an isomorphism
$$ma:\widehat\F_1\buildrel\sim\over\rightarrow \G\ARI^{pol}.\eqno(3.5.12)$$
Restricting this isomorphism to the subgroup of group-like power series
$MT=exp^\odot\mt$ yields an isomorphism
$$ma:MT\buildrel\sim\over\rightarrow \G\ARI^{pol}_{as},\eqno(3.5.13)$$
where $\G\ARI_{as}$ is the group of {\it symmetral} moulds, i.e. moulds
$A$ satisfying
$$\sum_{\w\in sh(\u,\v)} A(\w)=A(\u)A(\v).\eqno(3.5.14)$$
\vskip .5cm
With this background situation established, let us now explain one of
the identities from \S 2.7 in the power series situation.  We consider the
equality of automorphisms (2.8.6).

For $f,g,g'\in \widehat\F_1$, we define endomorphisms $X_{(g,g')}$, $R_f$ 
and $N_f$ of $\widehat\F_1$ as follows: each one sends $x\mapsto x$, and
$$\cases{X_{(g,g')}(y)=gyg'\cr
R_f(y)=fyf^{-1}\cr
N_f(y)=yf,}$$
i.e. $R_f=X_{(f,f^{-1})}$ and $N_f=X_{(1,f)}$.  
We have 
$$\cases{ma_{X_{(g,g')}(f)}=gaxit(ma_g,ma_{g'})\cdot ma_f\cr
ma_{R_g(f)}=garit(ma_g)\cdot ma_f\cr
ma_{N_g(f)}=ganit_{ma_g}\cdot ma_f,}\eqno(3.5.15)$$
where the second equality is (3.5.10) above and the others are analogous.
Just as $X_{(g,g')}$, $R_g$ and $N_g$ are automorphisms of the group
(under the usual multiplication) of power series with constant term $1$,
so $gaxit$, $garit$ and $ganit$ are automorphisms of $\G\ARI^{pol}$ equipped
with the multiplication $mu$.

We have
$$X_g\circ X_f=X_{X_{(g,g')}(f)},\eqno(3.5.16)$$
so if $f$ is such that $X_{(g,g')}(f)g=1$, then $ma_f=invgaxi(ma_g).$
Thus, the translation of the equality (2.8.6) back to the twisted Magnus
situation is given by
$$X_{(g,g')}\circ R_f=N_{g'g},\eqno(3.5.17)$$
where $g^{-1}=X_{(g,g')}(f)$, i.e. $ma_f=invgaxi(ma_g)$.
But it is easy to prove (3.5.17). Indeed, the automorphisms on both sides
fix $x$, so we only need to compare their images on $y$.  The RHS
yields $N_{g'g}(y)=yg'g$, and the LHS yields
$$\eqalign{X_{(g,g')}R_f(y)&=X_{(g,g')}(fyf^{-1})\cr
&=X_{(g,g')}(f)gyg'X_{(g,g')}(f^{-1})\cr
&=yg'X_{(g,g')}(f^{-1})\cr
&=yg'g,}$$
which proves that they are equal.

\vskip .8cm
\vfill\eject
\centerline{\bf Chapter 4}
\vskip .5cm
\centerline{\bf The mould pair $pal/pil$ and its properties}
\vskip .8cm
\noindent {\bf \S 4.1. Diffeomorphisms and the mould $pil$}
\vskip .5cm 
The passage from the space $DIFF_{\langle x\rangle}$ of diffeomorphisms
$f(x)=x(1+\sum_{r\ge 1} a_rx^r)$ to $\GV\ARI$ is one of Ecalle's key discoveries.
Given $f(x)$, he defines an associated mould $p_f$ in $\GV\ARI$, in fact giving two
equivalent definitions for $p_f$.  These stem from two functions associated 
to $f(x)$, namely the {\it infinitesimal dilator} $f_\#(x)$, defined by
$$f_\#(x)=x-{{f(x)}\over{f'(x)}}=\sum_{r\ge 1} \gamma_rx^{r+1},\eqno(4.1.1)$$
and the {\it infinitesimal generator} $f_*(x)$ defined by
$$f_*(x)=\sum_{r\ge 1}\epsilon_rx^{r+1}\eqno(4.1.2)$$
where the coefficients $\epsilon_r$ are determined by the identity
$$\Bigl(exp(f_*(x){{d}\over{dx}})\Bigr)\cdot x=f(x).$$

Let $re_1={{1}\over{v_1}}$, and for $r>1$ define the mould $re_r$ recursively
by $re_r=arit(re_{r-1})\cdot re_1$.  The mould $re_r$ is concentrated in
depth $r$, and it is easy to show by induction that it has explicit expression
$$re_r(v_1,\ldots,v_r)={{v_1+\cdots+v_r}\over{v_1(v_1-v_2)\cdots (v_{r-1}-v_r)
v_r}}.\eqno(4.1.3)$$
Let $lop_f$ denote the mould in $\overline{\ARI}$ defined by
$$lop_f(v_1,\ldots,v_r)=\epsilon_rre_r(v_1,\ldots,v_r)=\epsilon_r{{v_1+\cdots+v_r}
\over{v_1(v_1-v_2)\cdots (v_{r-1}-v_r)v_r}}\ \ \ {\rm for}\ \ r\ge 1.\eqno(4.1.4)$$
The first definition of the mould $p_f$ associated to $f(x)$ comes from
the infinitesimal generator of $f(x)$ and is given by
$$p_f=exp_{ari}(lop_f).\eqno(4.1.5)$$
By construction, the moulds $p_f$ associated to $f$ satisfy
$$p_{f\circ g}=gari(p_f,p_g).\eqno(4.1.6)$$
The second definition comes from the infinitesimal dilator, via the mould
$d_f\in \overline{\ARI}$ defined by
$$d_f(v_1,\ldots,v_r)=\gamma_rre_r(v_1,\ldots,v_r)=
\gamma_r{{v_1+\cdots+v_r}\over{v_1(v_1-v_2)\cdots (v_{r-1}-v_r)v_r}}\ \ \ {\rm for}\ \ r\ge 1;\eqno(4.1.7)$$
we define the mould $p_f$ recursively by setting $p_f(\emptyset)=1$ and
$$der\cdot p_f=preari(p_f,d_f),\eqno(4.1.8)$$
where $der$ is the operator on moulds such that 
$$\bigl(der\cdot A\bigr)(w_1,\ldots,w_r)=r\,A(w_1,\ldots,w_r).$$
Indeed, note that since $d_f(\emptyset)=0$, the depth $r$ term of $p_f$ 
can be deduced from the parts of $p_f$ up to depth $r-1$ via the right-hand 
side of (4.1.8).
\vskip .4cm
\noindent {\bf Proposition 4.1.1.} {\it The two definitions of $p_f$
are equivalent.}
\vskip .2cm
\noindent {\bf Proof.} The main fact is that if we apply the linearization
procedure, working in $k[[\epsilon]]/(\epsilon^2)$, then the 
linearized dilator $1+\epsilon\ f_\#(x)$ satisfies
the identity
$$\bigl(f\circ (1+\epsilon\ f_\#)\bigr)(x)=f(x)+
\epsilon\ \sum_{n\ge 1} na_nx^{n+1}.$$
Passing to the associated moulds by (4.1.6), using (2.8.4), the
left-hand side maps to
$$gari(p_f,p_{1+\epsilon\ f_\#})=p_f+\epsilon\ preari(p_f,p_{f_\#}).$$
We also see that the sum $\sum_{n\ge 1}na_nx^{n+1}$ maps to
$der\cdot p_f$ since each term is multiplied by its degree, so the
right-hand side altogether maps to 
$$p_f+\epsilon\ der\cdot p_f.$$
This shows that $p_f$ satisfies (4.1.8).\hfill{$\square$}
\vskip .4cm
\noindent {\bf Proposition 4.1.2.} {\it The moulds $p_f$ are symmetral.}
\vskip .2cm
\noindent {\bf Proof.} By Proposition 2.6.1, since $p_f=exp_{ari}(lop_f)$, 
it is enough to show that $lop_f$ is
alternal.  But $re_1$ is trivially alternal since it is concentrated in
depth 1. Assuming as an induction hypothesis that $re_{r-1}$ is alternal,
we see by Proposition 2.5.2 that $re_r=arit(re_{r-1})\cdot re_1$ is also
alternal, which proves that $lop_f$ is alternal.\hfill{$\square$}
\vskip .3cm
\noindent {\bf Definition.} Let $pil$ be the mould $p_f$ constructed
as above, where $f(x)=1-e^{-x}$, and let $dipil$ denote the mould $d_f$ for
this $f$.  In low depths, we have
$$\cases{pil(v_1)={{-1}\over{2v_1}}\cr
pil(v_1,v_2)={{1}\over{12}}{{2v_1-v_2}\over{v_1(v_1-v_2)v_2}}\cr
pil(v_1,v_2,v_3)={{-1}\over{24}}{{1}\over{(v_1-v_2)v_2v_3}}\cr
pil(v_1,v_2,v_3,v_4)={{1}\over{720}}{{6v_1v_3-10v_1v_4+v_2v_3+5v_2v_4-4v_3^2+v_3v_4}\over{v_1v_3v_4(v_1-v_2)(v_2-v_3)(v_3-v_4)}}.}$$
\vskip .2cm
\noindent {\bf Remarks.} Ecalle gives some very pretty results on
moulds associated to diffeomorphisms that we cite here without proof.  
\vskip .2cm
\noindent (1) A mould $A\in\GV\ARI$ lies in the image of $DIFF_{\langle x
\rangle}$ if and only if there exist constants $c_r$, $r\ge 1$ such that
$$mu\bigl(anti\cdot swap(A),swap(A)\bigr)=c_r{{1}\over{u_1\cdots u_r}},
\eqno(4.1.9)$$
and if this is the case, then $A=p_f$ where $f(x)=x+\sum_{r\ge 1}{{c_r}
\over{r+1}}x^{r+1}$.
\vskip .2cm
\noindent (2) If a mould $A\in\G\ARI$ is symmetral, then
$mu(anti\cdot A,A)$ is also symmetral. Therefore, setting $A=swap(p_f)$,
it is a necessary condition for the bisymmetrality of $p_f$ that
$mu(anti\cdot A,A)$ be symmetral, i.e. that the mould defined by the
right-hand side of (4.1.9) be symmetral. One can show directly that
the only mould of this form which is symmetral is the one where
$c_r=(-1)^r/r!$, i.e. $mu(anti\cdot A,A)=expmu({\cal O})$ where ${\cal O}$ is
the mould concentrated in depth $1$ defined by ${\cal O}(u_1)=1/u_1$.
Thus, since we can get the diffeomorphism $f$ back from the $c_r$ by
setting $a_r=c_r/(r+1)=(-1)^r/(r+1)!$, we find that the only diffeomorphism $f$
for which $p_f$ could be bisymmetral is 
$$f(x)=x+\sum_{r\ge 1} {{(-1)^r}\over{(r+1)!}}x^{r+1}=1-e^{-x}.$$
The next two sections will be devoted to giving Ecalle's direct proof, not
relying on this property, that $pil$ is indeed bisymmetral.
\vskip .8cm
\noindent {\bf \S 4.2. Two definitions of the mould $pal$}
\vskip .5cm
The mould pair $pal/pil$ is undoubtedly one of Ecalle's most beautiful 
and powerful discoveries.  In this chapter we give the most recent
definition that Ecalle has given for the mould $pal$ (cf. [Eupolars]), and
then give the complete proof that $pal=swap(pil)$.
\vskip .3cm
\noindent {\bf Definition 4.2.1.} Let $dur$ be the mould operator 
defined by $dur\cdot G(\emptyset)=0$ and for $r\ge 1$,
$$dur\cdot G(u_1,\ldots,u_r)=(u_1+\cdots+u_r)\,
G(u_1,\ldots,u_r).\eqno(4.2.1)$$
Let $du$ be the mould operator on $G\ARI$ defined by 
$$duG=mu(invmu(G),dur\cdot G).\eqno(4.2.2)$$
Inversely, if $duG$ is a given mould in $\ARI$, then the mould $G\in G\ARI$
satisfying (4.2.1) can be recovered depth by depth from $duG$
starting with $G(\emptyset)=1$, then using the formula
$$dur\cdot G=mu(G,duG).\eqno(4.2.3)$$
\'Ecalle calls the mould $duG$ the {\it $mu$-dilator} of $G$.
\vskip .1cm
\noindent {\bf Definition 4.2.2.} 
Let $dupal\in \ARI$ be the mould defined explicitly as follows:
$dupal(\emptyset)=0$ and for each $r\ge 1$, 
$$\eqalign{dupal(u_1,\ldots,u_r)&={{B_r}\over{r!}}\sum_{i=0}^{r-1} (-1)^i
\Bigl({{r-1}\atop{i}}\Bigr){{1}\over{u_1\cdots\hat u_{r-i}\cdots u_r}}\cr
&={{B_r}\over{r!}}{{1}\over{u_1\cdots u_r}}
\bigl(\sum_{i=0}^{r-1} (-1)^i\Bigl({{r-1}\atop{i}}\Bigr)u_{i+1}\bigr).}\eqno(4.2.4)$$
Note in particular that $dupal(u_1,\ldots,u_r)=0$ for all odd $r>1$.
The mould $pal\in \G\ARI$ is defined by $pal(\emptyset)=1$ and then, 
recursively depth by depth as in (4.2.3), by the formula
$$dur\cdot pal=mu(pal,dupal),\eqno(4.2.5)$$
where $dur$ is as in (4.2.1a).
Up to depth 4, we have
$$\cases{pal(u_1)=-{{1}\over{2u_1}}\cr
pal(u_1,u_2)={{1}\over{12}}{{u_1+2u_2}\over{u_1u_2(u_1+u_2)}}\cr
pal(u_1,u_2,u_3)={{-1}\over{24}}{{1}\over{u_1(u_1+u_2)u_3}}\cr
pal(u_1,u_2,u_3,u_4)=-{{1}\over{720}}{{u_1^2-2u_1u_2-2u_1u_3+4u_1u_4
-3u_2^2-7u_2u_3-6u_2u_4}\over{u_1u_2u_3u_4(u_1+u_2)(u_1+u_2+u_3+u_4)}},}$$
\vskip .3cm
\noindent {\bf Theorem 4.2.3.} {\it We have $pal=swap(pil)$.}
\vskip .2cm
\noindent {\bf Proof.} We need two preliminary results.
\vskip .3cm
\noindent {\bf Lemma 4.2.4.} {\it The derivations $dur$ and $der$ commute, and 
for any mould $B\in \ARI$, $dur$ commutes with $amit(B)$, $anit(B)$, 
$arit(B)$ and $irat(B)$.}
\vskip .2cm
\noindent {\bf Proof.} The commutation of $der$ and $dur$ is obvious since
$der\cdot dur$ and $dur\cdot der$ both come down to multiplying the mould $A$
by $r(u_1+\cdots+u_r)$ in depth $r$.  The commutation of $dur$ with 
$arit(B)$ and $irat(B)$ follow immediately from the commutation with
$amit(B)$ and $anit(B)$ since $arit(B)=amit(B)-anit(B)$ by (2.2.4) and
$irat(B)=amit(B)-anit(push(B))$ by (2.4.11).
Looking at the definition of $amit(B)$ in (2.2.1), we see that 
$$amit(B)\cdot dur\cdot A(\w)=\sum_{{{\w=\a\b\c}\atop{\b,\c\ne\emptyset}}}
(dur\cdot A)(\a\lceil\c)B(\b).$$
But if $\a=(u_1,\ldots,u_i)$, $\b=(u_{i+1},\ldots,u_{i+k})$ and
$\c=(u_{i+k+1},\ldots,u_r)$, we have
$$\a\lceil\c=\bigl(u_1,\ldots,u_i,u_{i+1}+\cdots+u_{i+k+1},u_{i+k+2}, \ldots,u_r),\eqno(4.2.6)$$
we see that $(dur\cdot A)(\a\lceil\c)=(u_1+\cdots+u_r)A(u_1,\ldots,u_r)$,
so the same factor $(u_1+\cdots+u_r)$ occurs in every term of the sum
over $\w=\a\b\c$ and therefore can be taken outside the sum, leaving
exactly $dur\cdot amit(B)\cdot A$.  The exact same argument holds for 
$anit(B)$ (defined in (2.2.2)), with $\a\rceil\c$ instead of $\a\lceil\c$.
This concludes the proof.\hfill{$\square$}
\vskip .3cm
\noindent {\bf Definition 4.2.5.} Let $dipil\in \overline{\ARI}$ be 
the mould $d_f$ of the previous section, with
$f(x)=1-e^{-x}$.  Explicitly,
$$dipil(v_1,\ldots,v_r)={{-1}\over{(r+1)!}}re_r(v_1,\ldots,v_r)=
{{-1}\over{(r+1)!}}
{{v_1+\cdots+v_r}\over{v_1(v_1-v_2)\cdots (v_{r-1}-v_r)v_r}},\eqno(4.2.7)$$ 
and by (4.1.8), we have
$$der\cdot pil=preari(pil,dipil).\eqno(4.2.8)$$
\vskip .2cm
\noindent {\bf Proposition 4.2.6.} {\it Set $dapal=swap(dipil)$.  Then
$$der\cdot dupal=dur\cdot dapal+irat(dapal)\cdot dupal-lu(dapal,dupal).
\eqno(4.2.9)$$}
The detailed proof of this identity is given in the Appendix, \S A.6.
\vskip .3cm
We can now complete the proof of Theorem 4.2.3.
We first apply the $swap$ to (4.2.8), obtaining
$$\eqalign{der\cdot swap(pil)&=swap\bigl(preari(pil,dipil)\bigr)\cr
&=swap\bigl(preari(swap(swap(pil)),dapal)\bigr)\cr
&=preira\bigl(swap(pil),dapal)\bigr).}\eqno(4.2.10)$$
Given that $swap(pil)(\emptyset)=1$, (4.2.10) can actually be used as a recursive
depth-by-depth definition for $swap(pil)$; i.e. we have two equivalent ways to compute 
$swap(pil)$, either by swapping the terms of $pil$ or by (4.2.10).
Therefore, if $pal$ is the mould defined in (4.2.5), to show that $pal=swap(pil)$,
it suffices to prove that $pal$ satisfies (4.2.10), i.e. that
$$der\cdot pal=preira(pal,dapal).\eqno(4.2.11)$$
Set
$$\eqalign{A&=der\cdot pal-preira(pal,dapal)\cr
&=der\cdot pal-irat(dapal)\cdot pal-mu(pal,dapal).}\eqno(4.2.12)$$
We apply $der$ to the left hand side of (4.2.5).  
Using the fact that $irat(dapal)$ is a $mu$-derivation, we have
$$\eqalign{der\cdot &dur\cdot pal=der\cdot mu(pal,dupal)\ \ {\rm by\ (4.2.5)}\cr
&=mu(der\cdot pal,dupal)+mu(pal,der\cdot dupal)\cr
&=mu(der\cdot pal,dupal)+mu(pal,irat(dapal)\cdot dupal)\cr
&\ \ \ +mu(pal,dur\cdot dapal)-mu(pal,dapal,dupal)+mu(pal,dupal,dapal)
\ \ {\rm by\ (4.2.9)}\cr
&=mu(der\cdot pal,dupal)+irat(dapal)\cdot mu(pal,dupal)
-mu(irat(dapal)\cdot pal,dupal)\cr
&\ \ +mu(pal,dur\cdot dapal)-mu(pal,dapal,dupal)+mu(pal,dupal,dapal)\cr
&=mu(der\cdot pal,dupal)-mu(irat(dapal)\cdot pal,dupal)
-mu(pal,dapal,dupal)\cr
&\ \ +irat(dapal)\cdot mu(pal,dupal)
+mu(pal,dur\cdot dapal)
+mu(pal,dupal,dapal)\cr
&=mu(A,dupal)+irat(dapal)\cdot mu(pal,dupal)+
mu(pal,dur\cdot dapal) +mu(pal,dupal,dapal)\cr
&=mu(A,dupal)+irat(dapal)\cdot dur\cdot pal+mu(pal,dur\cdot dapal)+mu(pal,dupal,dapal)\cr
&=mu(A,dapal)+irat(dapal)\cdot dur\cdot pal+mu(pal,dur\cdot dapal)+mu(dur\cdot pal,dapal)
\cr
&=mu(A,dupal)+irat(dapal)\cdot dur\cdot pal+dur\cdot mu(pal,dapal)\cr
&=mu(A,dupal)+dur\cdot irat(dapal)\cdot pal+dur\cdot mu(pal,dapal)\ \
{\rm by\ Lemma\ 4.2.2.} }$$
By Lemma 4.2.4, we also have $der\cdot dur\cdot pal=dur\cdot der\cdot pal$,
and the equality of $der\cdot dur\cdot pal$ with the last line above can
thus be rewritten as
$$dur\cdot der\cdot pal-
dur\cdot irat(dapal)\cdot pal-dur\cdot mu(pal,dapal)=mu(A,dupal),$$
i.e.
$$dur\cdot A=mu(A,dupal).\eqno(4.2.13)$$
Now, although this looks like the defining equation (4.2.5) for $pal$, in fact
the defining equation (4.2.12) for $A$ shows that $A(\emptyset)=0$. But it is
easy to show that if a mould $A$ satisfies $A(\emptyset)=0$ and (4.2.13), then
$A$ is identically $0$. Indeed, suppose by induction that $A(u_1,\ldots,
u_i)=0$ for $0\le i<r$.  Then
$$\eqalign{(u_1+\cdots+u_r)A(u_1,\ldots,u_r)&=\sum_{i=0}^r A(u_1,\ldots,u_i)dupal(
u_{i+1},\ldots,u_r)\cr
&=A(u_1,\ldots,u_r)dupal(\emptyset)\cr
&=0,}$$
so $A(u_1,\ldots,u_r)=0$.  Thus the expression (4.2.12) is equal to $0$,
proving the desired identity (4.2.11).
This concludes the proof of Theorem 4.2.3.\hfill{$\square$}
\vskip .8cm
\noindent {\bf \S 4.3. Symmetrality of $pal$}
\vskip .5cm
Let
$$Paj(r_1,\ldots,r_s)={{1}\over{r_1(r_1+r_2)\cdots (r_1+\cdots+r_s)}}.$$
\vskip .3cm
\noindent {\bf Lemma 4.3.1.} {\it The mould $Paj$ is symmetral.}
\vskip .2cm
\noindent {\bf Proof.} 
Let us proceed by induction on the length of the shuffles
$sh(\u,\v)$, i.e. the total length of the two words $\u$ and $\v$.  When
$\u$ and $\v$ are both of length $1$, i.e. $\u=(r_1)$, $\u=(r_2)$, we have
$$\sum_{\w\in sh(\u,\v)} Paj(\w)=Paj(r_1,r_2)+Paj(r_2,r_1)=
{{1}\over{r_1(r_1+r_2)}}+{{1}\over{r_2(r_1+r_2)}}={{1}\over{r_1r_2}},$$
so $Paj$ is symmetral in length 2.  Assume it is symmetral up to length $s-1$,
and let $\u=(r_1,\ldots,r_l)$, $\v=(r_{l+1},\ldots,r_s)$ be two words of
total length $s$.  We use the recursive definition 
$$sh(\u,\v)=sh(\u',\v)\cdot r_l+sh(\u,\v')\cdot r_s,$$
where $\u'=(r_1,\ldots,r_{l-1})$ and $\v'=(r_{l+1},\ldots,r_{s-1})$.  Letting
$R=\sum_{i=1}^s r_i$, we have
$$\eqalign{\sum_{\w\in sh(\u,\v)} Paj(\w)&=
\sum_{\w\in sh(\u',\v)} Paj(\w,r_l)+\sum_{\x\in sh(\u,\v')} Paj(\x,r_s)\cr
&=\sum_{\w\in  sh(\u',\v)} Paj(w_1,\ldots,w_{s-1},r_l)+
\sum_{\x\in  sh(\u,\v')} Paj(x_1,\ldots,x_{s-1},r_s)\cr
&=\sum_{\w\in  sh(\u',\v)} {{1}\over{w_1(w_1+w_2)\cdots(w_1+\cdots+w_{s-1})R}}\cr
&\ \ \ \ \ \ \ \ +\sum_{\x\in  sh(\u,\v')} {{1}\over{x_1(x_1+x_2)\cdots(x_1+\cdots+x_{s-1})R}}\cr
&={{1}\over{R}}
\sum_{\w\in  sh(\u',\v)} {{1}\over{w_1(w_1+w_2)\cdots(w_1+\cdots+w_{s-1})}}\cr
&\ \ \ \ \ \ \ \ +{{1}\over{R}}\sum_{\x\in  sh(\u,\v')} {{1}\over{x_1(x_1+x_2)\cdots(x_1+\cdots+x_{s-1})}}\cr
&={{1}\over{R}}\sum_{\w\in sh(\u',\v)} Paj(\w)+{{1}\over{R}}
\sum_{\x\in sh(\u,\v')} Paj(\x)\cr
&={{1}\over{R}}Paj(\u')Paj(\v)+{{1}\over{R}}Paj(\u)Paj(\v')\ \ 
{\rm by\ the\ induction\ hypothesis}\cr}$$
$$\eqalign{
&={{1}\over{R}}{{1}\over{r_1(r_1+r_2)\cdots(r_1+\cdots+r_{l-1})}}
{{1}\over{r_{l+1}(r_{l+1}+r_{l+2})\cdots (r_{l+1}+\cdots +r_s)}}\cr
&\ \ \ \ +{{1}\over{R}}{{1}\over{{r_1(r_1+r_2)\cdots(r_1+\cdots+r_l})}}
{{1}\over{r_{l+1}(r_{l+1}+r_{l+2})\cdots (r_{l+1}+\cdots +r_{s-1})}}\cr
&=\Bigl({{1}\over{R}}\Bigr)\Bigl((r_1+\cdots+r_l)+(r_{l+1}+\cdots+r_s)\Bigr)
Paj(\u)Paj(\v)\cr
&=Paj(\u)Paj(\v).}$$
This proves that $Paj$ is symmetral.\hfill{$\square$}
\vskip .5cm
\noindent {\bf Lemma 4.3.2.} {\it Let $S$ be a mould such that
$S(\emptyset)=1$. Then the defining formula
$$dur\cdot S=mu(S,duS)\eqno(4.3.1)$$
is equivalent to the inversion formula
$$S({\bf u})=1+\sum_{{{\bf u_1\cdots u_s=u}\atop{\u_i\ne\emptyset}}} Paj(|{\bf u_1}|,
\ldots |{\bf u_s}|)\,duS({\bf u_1})\cdots duS({\bf u_s}),\eqno(4.3.2)$$
where if $\u=(r_1,\ldots,r_l)$ then $|\u|=r_1+\cdots +r_l$.}
\vskip .2cm
\noindent {\bf Proof.} 
We prove the equivalence of (4.3.1) and (4.3.2) by induction on the length
of $\u$.  When $\u=\emptyset$, the constant term 1 on the right-hand side
of (4.3.2) ensures equality . For $\u=(u_1)$, we have 
$$S(u_1)=Paj(u_1)duS(u_1)={{1}\over{u_1}}duS(u_1)$$ 
from (4.3.2), and from (4.3.1) we have
$$u_1S(u_1)=S(\emptyset)duS(u_1)=duS(u_1)$$ 
so they are equivalent.  This settles the base case.
Now assume the induction hypothesis that (4.3.1) and (4.3.2) give the same
formula for $S(u_1,\ldots,u_i)$ for $i<r$.
From (4.3.1), and using the induction hypothesis on each term in $S$, we have
$$\eqalign{(u_1+\cdots+u_r)&S(u_1,\ldots,u_r)=
\sum_{i=0}^{r-1} S(u_1,\ldots,u_i)duS(u_{i+1},\ldots,u_r)\cr
&=\sum_{i=0}^{r-1} \sum_{(u_1,\ldots,u_i)=\u_1\cdots\u_s} 
Paj(|\u_1|,\cdots,|\u_s|)duS(\u_1)\cdots duS(\u_s)duS(u_{i+1},\ldots,u_r),
}$$
so writing $\u_{s+1}=(u_{i+1},\ldots,u_r)$ in each term and dividing both
sides by $R=(u_1+\cdots+u_r)$, we find
$$\eqalign{S(u_1,\ldots,u_r)&=
\sum_{1\le |\u_{s+1}|\le r}\sum_{\u=\u_1\cdots\u_s\u_{s+1}} 
{{1}\over{R}}Paj(|\u_1|,\cdots,|\u_s|)duS(\u_1)\cdots duS(\u_s)duS(\u_{s+1})\cr
&=\sum_{\u=\u_1\cdots\u_s\u_{s+1}} 
{{|\u_1|+\cdots+|\u_{s+1}|}\over{R}}Paj(|\u_1|,\cdots,|\u_s|,|\u_{s+1}|)duS(\u_1)\cdots duS(\u_{s+1})\cr
&=\sum_{\u=\u_1\cdots\u_s\u_{s+1}} 
Paj(|\u_1|,\cdots,|\u_{s+1}|)duS(\u_1)\cdots duS(\u_{s+1}).}$$
This proves that (4.3.1) is equivalent to (4.3.2).\hfill{$\square$}
\vskip .5cm
\noindent {\bf Proposition 4.3.3.} {\it Let $S$ and $duS$ be two moulds
related as in (4.3.1).  If $duS$ is alternal, then $S$ is symmetral.}
\vskip .2cm
\noindent {\bf Proof.} We will use the equivalent formula (4.3.2) for $S$.
Indeed, by formula (2.1.1) for mould composition, we see
that (4.3.2) is equivalent to the statement that the definition
$dur\cdot S=mu(S,duS)$ is equivalent to $S=1+Paj\circ duS$.  Assume
that $duS$ is alternal.  From Lemma 2.6.2, we know that for any alternal mould 
$A$ and symmetral mould $B$, the composition $B\circ A$ is symmetral, 
and from Lemma 4.3.1 we know that $Paj$ is symmetral. This concludes the
proof.\hfill{$\square$}
\vskip .3cm
\noindent {\bf Theorem 4.3.4.} {\it The mould $pal$ is symmetral.}
\vskip .2cm
\noindent {\bf Proof.} Thanks to Proposition
4.3.3, it is enough to show that $dupal$ is alternal.  But this reduces in fact
to an easy exercise, namely showing that the only linear alternal moulds
$a_1u_1+\cdots+a_ru_r$ are, up to scalar multiple, the binomial moulds
$\sum_{i=1}^r (-1)^i\Bigl({{r-1}\atop{i-1}}\Bigr) u_i.$\hfill{$\square$}

\vskip .8cm
\noindent {\bf \S 4.4. The identity $crash(pal)=pac$}
\vskip .5cm
Let $pac$ be the mould defined by 
$$pac(u_1,\ldots,u_r)={{1}\over{u_1\cdots u_r}}\eqno(4.4.1)$$
and let $pic$ be defined by
$$pic(v_1,\ldots,v_r)={{1}\over{v_1\cdots v_r}}.\eqno(4.4.2)$$
In this section we show two identities (4.4.3) and (4.4.8) that are essential
to the proof of the second fundamental identity (4.5.2) stated and proved
below in \S 4.5.
\vskip .3cm
\noindent {\bf Lemma 4.4.1.} {\it We have
$$crash(pal):=mu(push\cdot swap\cdot invmu\cdot invpil,swap\cdot invpil)=pac.
\eqno(4.4.3)$$}
\vskip .1cm
\noindent {\bf Proof.} Since $pil$ is symmetral, we have 
$$mu\bigl(pari\cdot anti(pil),pil\bigr))=1,\eqno(4.4.4)$$
and it's easy to see by the homogeneous degrees of $pil$ that
$$anti\cdot neg(pil)=pari\cdot anti(pil),\eqno(4.4.5)$$
so we find that 
$$anti\cdot neg(pil)=invmu(pil).\eqno(4.4.6)$$
Now, because of (4.4.6), we find that $pil\in G\ARI\cap GAWI$ (see [E,p. 44]
for definition of $GAWI$), and thus the $gari$ and $gawi$ inverses are the
same, so it makes sense to write $invpil\in G\ARI\cap GAWI$.  This
means that for $pil$ and $invpil$ we have 
$$\cases{push\cdot swap\cdot invmu\cdot swap\cdot swap(pil)=anti\cdot swap(pil)\cr
push\cdot swap\cdot invmu\cdot swap\cdot swap(invpil)=anti\cdot swap(invpil).
}\eqno(4.4.7)$$
Thus the LHS of (4.4.3) is equal to
$$crash(pal)=mu\bigl(anti\cdot swap(invpil),swap(invpil)\bigr),$$
which is nothing other than $gepar(invpil),$ so we can use \S 4.1.3 for
$f(x)=-log(1-x)$ which shows that
$$gepar(invpil)=pic,$$
proving (4.4.3). \hfill{$\square$}
\vskip .3cm
\noindent {\bf Lemma 4.4.2.} {\it 
We have 
$$ganit_{pic}\cdot invpil=swap\cdot invpal.\eqno(4.4.8)$$}
\vskip .1cm
\noindent {\bf Proof.} From (2.9.17) applied to $A=1$, $B=pal$, we have
$$swap\cdot invgari\cdot swap\cdot pal=swap\cdot invpil=
ganit_{crash\cdot pal}(invpal).\eqno(4.4.9)$$
Using (2.9.12), from (4.4.3) we also know that
$$ganit_{pac}\cdot invpal=swap\cdot invpil.$$
We need to use the elementary result
$$invgani(pac)=pari\cdot anti\cdot paj,\eqno(4.4.10)$$
where 
$$paj(u_1,\ldots,u_r)={{1}\over{(u_1(u_1+u_2)(u_1+u_2+u_3)\cdots (u_1+\cdots
+u_r)}}.$$
This gives
$$invpal=ganit_{pari\cdot anti\cdot paj}\cdot swap\cdot invpil,$$
so
$$swap\cdot invpal=swap\cdot ganit_{pari\cdot anti\cdot paj}\cdot swap\cdot invpil.$$
It remains only to prove that the following two automorphisms of G\ARI are 
equal:
$$ganit_{pic}=swap\cdot ganit_{pari\cdot anti\cdot paj}\cdot swap.\eqno(4.4.11)$$

Now, every mould $C$ in the $v_i$ such that $C(v_1,\ldots,v_r)$ is actually a 
rational function $B$ of the variables $v_2-v_1,\ldots,v_r-v_1$ satisfies
the identity $C=ganit_B(Y)$, 
by the calculation
$$\eqalign{ganit_B(Y)(v_1,\ldots,v_r)&=\sum_{b_1c_1\cdots b_sc_s}
Y(b_1\cdots b_s)B(\lfloor c_1)\cdots B(\lfloor c_2)\cr
&=\sum_{b_1=(v_1),c_1=(v_2,\ldots,v_r)} Y(v_1)B(v_2-v_1,\ldots,v_r-v_1)\cr
&=B(v_2-v_1,\ldots,v_r-v_1)\cr
&=C(v_1,\ldots,v_r).}\eqno(4.4.12)$$

Let us write $swap(Y)=Y$ a little abusively, since although the values 
in depths 0 and 1 are still 1, $swap(Y)$ is considered a mould in the $u_i$.
We start to compute the right-hand side of (4.4.11) explicitly as
$$ganit_{pari\cdot anti\cdot paj}\cdot Y(u_1,\ldots,u_r)=
{{(-1)^{r-1}}\over{u_r(u_{r-1}+u_r)\cdots (u_2+\cdots u_r)}}$$
(with $ganit_{pari\cdot anti\cdot paj}\cdot Y(\emptyset)=1$,
$ganit_{pari\cdot anti\cdot paj}\cdot Y(u_1)=1$).
Swapping this, we obtain for the RHS of (4.4.11):
$$swap\cdot ganit_{pari\cdot anti\cdot paj}\cdot Y(u_1,\ldots,u_r)=
{{1}\over{(v_2-v_1)(v_3-v_1)\cdots (v_r-v_1)}}.$$
Letting 
$$C(v_1,\ldots,v_r)={{1}\over{(v_2-v_1)(v_3-v_1)\cdots (v_r-v_1)}},$$
we see by 4.4.12)) that $C=ganit_B(Y)$ where
$$B(v_1,\ldots,v_r)={{1}\over{v_1\cdots v_r}},\eqno(4.4.13)$$
i.e. $B=pic$.\hfill{$\square$}
\vskip .5cm
Note that we have not shown that $crash(pil)=pic$, although
it seems to be true. However, the above result is enough for our purposes,
together with the important result stated by \'Ecalle concerning 
the automorphism $ganit_{pic}$ given in the corollary to Theorem
4.4.3 below.
\vskip .3cm
\noindent {\bf Proposition 4.4.3} {\it Let $A,B\in \overline{\ARI}$ be such
that $A=ganit_{pic}\cdot B$.  Then $A$
satisfies the stuffle relations if and only if $B$ satisfies the shuffle 
relations, i.e. $A_{r,s}(v_1,\ldots,v_r)=0$ for all pairs $(r,s)$, where
$A_{r,s}$ defined as in (2.3.3), if and only if $B$ is alternal.}
\vskip .2cm
\noindent {\bf Proof.} The full and complex proof of this fundamental
statement has been worked out by N.~Komiyama in [K], 
Theorem 3.24.\hfill{$\square$}
\vskip .5cm
\noindent {\bf Corollary.} {\it Let $A=ganit_{pic}\cdot B$.  Then $A$ 
satisfies the stuffle relations, i.e.~$A_{r,s}(v_1,\ldots,v_r)=0$ for
all $(r,s)$, if and only if $B$ satisfies the shuffle relations.}
\vskip .5cm
\noindent {\bf \S 4.5. Ecalle's second fundamental identity}
\vskip .2cm
In this section we use Ecalle's first fundamental identity (2.9.4) and the 
results of \S 4.4 to prove another formula that is one of the main tools in his 
theory, namely the second fundamental identity, given in Theorem 4.5.2.
It will be deduced from an initial version given in the following proposition.
\vskip .3cm
\noindent {\bf Proposition 4.5.1.}  We have
$$swap\cdot fragari(swap\cdot A,pal)=ganit_{pic}\cdot
fragari(A,pil).\eqno(4.5.1)$$
\vskip .2cm
\noindent {\bf Proof.} Applying the fundamental identity (2.9.4) to 
$A=swap\cdot M$ and $B=pal$ and using Lemma 4.4.1 yields
$$\eqalign{swap\cdot fragari(M,swap\cdot pal)&=
ganit_{crash\cdot pal}\cdot fragari(swap\cdot M,pal)\cr
&=ganit_{pac}\cdot fragari(swap\cdot M,pal).}$$
Thus by (4.4.10) we have
$$\eqalign{ganit_{invgani\cdot pac}\cdot swap\cdot fragari(M,pil)&=
ganit_{pari\cdot anti\cdot paj}\cdot swap\cdot fragari(M,pil)\cr
&=fragari(swap\cdot M,pal).}$$
Applying swap to both sides and (4.4.11), we have
$$\eqalign{swap\cdot ganit_{pari\cdot anti\cdot paj}\cdot swap\cdot 
fragari(M,pil) &=
ganit_{pic}\cdot fragari(M,pil)\cr
& = swap\cdot fragari(swap\cdot M,pal),}$$
which proves the desired (4.5.1).\hfill{$\square$}
\vskip .3cm
\noindent {\bf Theorem 4.5.2.} {\it 
For every push-invariant mould $M$, we have {\rm Ecalle's second fundamental
identity}:
$$swap\cdot Ad_{ari}(pal)\cdot M=ganit_{pic}\cdot Ad_{ari}(pil)\cdot swap(M).
\eqno(4.5.2)$$}
\vskip .1cm
\noindent {\bf Proof.} We use the defining identity
$$Ad_{ari}(A)\cdot B=fragari\bigl(preari(A,B),A\bigr)\eqno(4.5.3)$$
and equation (2.4.10) given by
$$swap\bigl(preari(swap\cdot A,swap\cdot B)\bigr)=axit(B,-push(B)) \cdot A+
mu(A,B).\eqno(4.5.4)$$
Using this for $A=pal$ and $B=M$, we find in particular that
$$\eqalign{preari(pil,swap\cdot M)&=swap\bigl(axit(M,-push(M))\cdot pal
+mu(pal,M)\bigr)\cr
&=swap\bigl(arit(M)\cdot pal+mu(pal,M)\bigr) \ \ \hbox{because $M$ is push-inv}\cr
&=swap\cdot preari(pal,M).}\eqno(4.5.5)$$
Using (2.8.6) for $A=pal$, $B=M$, we have 
$$\eqalign{swap\cdot Ad_{ari}(pal)\cdot M&=
swap\cdot fragari\bigl(preari(pal,M),pal\bigr)\cr
&=swap\cdot fragari\bigl(swap\bigl(swap\cdot preari(pal,M)\bigr),pal\bigr)\cr
&=ganit_{pic}\cdot fragari(swap\cdot preari(pal,M),pil)
\ \ {\rm by\ (4.5.1)}\cr
&=ganit_{pic}\cdot fragari(preari(pil,swap\cdot M),pil)
\ \ {\rm by\ (2.8.6)}\cr
&=ganit_{pic}\cdot Ad_{ari}(pil)\cdot swap\cdot M,}$$
proving (4.5.2).\hfill{$\square$}
\vfill\eject
\noindent {\bf \S 4.6. Double shuffle is a Lie algebra}
\vskip .5cm
Recall that by Theorem 3.4.4, the double shuffle Lie algebra $\ds$
is isomorphic to $\ARI^{pol}_{\underline{al}*\underline{il}}$. In
this section we give Ecalle's proof that the latter is a Lie algebra
for the $ari$-bracket, thus giving a complete different proof of
Racinet's well-known theorem 1.3.1.  Our proof comes directly from the
paper [SS], and was indicated to us in a personal communication from Ecalle.
\vskip .3cm
\noindent {\bf Theorem 4.6.1.} {\it The action of the operator $Ad_{ari}(pal)$
on the Lie subalgebra $\ARI_{\underline{al}*\underline{al}}\subset \ARI$ yields a
Lie isomorphism of subspaces
$$Ad_{ari}(pal):\ARI_{\underline{al}*\underline{al}} 
\buildrel\sim\over\rightarrow \ARI_{\underline{al}*\underline{il}}.\eqno(4.6.1)$$
Thus in particular $\ARI_{\underline{al}*\underline{il}}$ forms a Lie algebra
under the $ari$-bracket.}
\vskip .2cm
\noindent {\bf Proof.} 
Let $A\in \ARI$ be an even function in depth 1. Note first that $Ad_{ari}(pal)$ 
preserves the depth 1 component of moulds in $\ARI$, so $Ad_{ari}(pal)\cdot A$
is also even in depth 1.

We first consider the case where $A\in \ARI_{\underline{al}/\underline{al}}$,
i.e. $swap(A)$ is alternal without addition of a constant correction.
By Proposition 2.6.1, $\G\ARI_{as}=exp_{ari}(\ARI_{al})$, so in particular
$\G\ARI_{as}$ acts by the adjoint action on $\ARI_{al}$, and therefore
since $pal$ is symmetral by Theorem 4.3.4, the mould $Ad_{ari}(pal)\cdot A$ 
is alternal.  By Lemma 2.5.5, $A$ is push-invariant, so we can apply Ecalle's 
second fundamental identity (4.5.2) and find that 
$$swap\bigl(Ad_{ari}(pal)\cdot A\bigr)=ganit_{pic}\cdot
\bigl(Ad_{ari}(pil)\cdot swap(A)\bigr).\eqno(4.6.2)$$
Since $A\in \ARI_{\underline{al}/\underline{al}}$, $swap(A)$ is alternal, and 
thus again by Proposition 2.6.1, $Ad_{ari}(pil)\cdot swap(A)$ is again
alternal; thus $ganit_{pic}\cdot Ad_{ari}(pil)\cdot swap(A)$ is alternil, 
and finally by
(4.6.2), $swap\bigl(Ad_{ari}(pal)\cdot A\bigr)$ is alternil, which proves that
$Ad_{ari}(pal)\cdot A\in \ARI_{\underline{al}/\underline{il}}$ as desired.
\vskip .2cm
We now consider the general case where $A\in\ARI_{\underline{al}*\underline{al}}$.  Let $C$ be the 
constant-valued mould such that $swap(A)+C$ is alternal. We will need 
the following result to deal with the constant mould $C$.
\vskip .3cm
\noindent {\bf Lemma 4.6.2.} [B, Corollary 4.43] {\it If $C$ is a constant-valued mould, then
$$ganit_{pic}\cdot Ad_{ari}(pil)\cdot C=C.\eqno(4.6.3)$$}
\noindent {\bf Proof.} 
We apply the fundamental identity (4.5.2) in the case where
$A=swap(A)=C$ is a constant-valued mould, obtaining
$$swap\bigl(Ad_{ari}(pal)\cdot C\bigr)=ganit_{pic}\cdot \bigl(Ad_{ari}(pil)\cdot C\bigr).$$
So it is enough to show that the left-hand side of this is equal to $C$,
i.e. that $Ad_{ari}(pal)\cdot C=C$.  
Directly from the definitions, we see that if $A\in\ARI$, 
then $arit(C)\cdot A=0$ and $arit(A)\cdot C=lu(C,A)$.  Thus
$$ari(A,C)=lu(A,C)+arit(A)\cdot C-arit(C)\cdot A=0.\eqno(4.6.4)$$
Now, by (2.8.5) we see that $Ad_{ari}(pal)\cdot C$ is a linear combination of 
iterated $ari$-brackets of $logari(pal)$ with $C$, but
since $pal\in\G\ARI$, $logari(pal)\in \ARI$, so (4.6.4) shows that 
$ari(logari(pal),C)=0$, i.e. all the terms in (2.8.5) are 0, which 
concludes the proof.\hfill{$\square$}
\vskip .3cm
Returning to the case $A\in\ARI_{\underline{al}*\underline{al}}$, we again have that
$Ad_{ari}(pal)\cdot A$ is alternal, so to conclude the proof of the theorem
it remains only to show 
that its swap is alternil up to addition of a constant mould, and
we will show that this constant mould is exactly $C$. 
As before, since $swap(A)+C\in\overline{\ARI}$ is alternal, the mould
$$Ad_{ari}(pil)\cdot \bigl(swap(A)+C\bigr)=Ad_{ari}(pil)\cdot swap(A)+
Ad_{ari}(pil)\cdot C$$ 
is also alternal.  Thus applying $ganit_{pic}$ to it yields 
the alternil mould 
$$ganit_{pic}\cdot Ad_{ari}(pil)\cdot swap(A) + ganit_{pic}\cdot 
Ad_{ari}(pil)\cdot C.$$
By Lemma 4.6.2, this is equal to
$$ganit_{pic}\cdot Ad_{ari}(pil)\cdot swap(A) + C,\eqno(4.6.5)$$
which is thus alternil. Now, since $A$ is push-invariant by Lemma 2.5.5,
we can apply (4.5.2) and find that (4.6.5) is equal to
$$swap\bigl(Ad_{ari}(pal)\cdot A\bigr)+C,$$
which is thus also alternil.  Therefore $swap\bigl(Ad_{ari}(pal)\cdot A\bigr)$
is alternil up to a constant, which precisely means that $Ad_{ari}(pal)\cdot A
\in \ARI_{\underline{al}*\underline{il}}$ as claimed. Since $Ad_{ari}(pal)$ is 
invertible (with inverse $Ad_{ari}(invgari\cdot pal)$), we can use all of
these arguments in the other direction to show that
$Ad_{ari}(invgari\cdot pal)$ maps $\ARI_{\underline{al}*\underline{il}}$ to 
$\ARI_{\underline{al}*\underline{al}}$. Thus (4.6.1) is a Lie
algebra isomorphism.  \hfill{$\square$}
\vskip .5cm
\noindent {\bf \S 4.7. The $\Delta$-denominator}
\vskip .5cm
Let $\Delta$ be the mould operator defined on moulds in the $u_i$ by
$$\Delta(A)(u_1,\ldots,u_r)=(u_1+\cdots+u_r)u_1\cdots u_r\,A(u_1,\ldots,u_r),
\eqno(4.7.1)$$
and on moulds in the $v_i$ by its swapped version
$$\Delta(A)(v_1,\ldots,v_r)=v_1(v_1-v_2)\cdots (v_{r-1}-v_r)v_r\,A(v_1,\ldots,v_r).\eqno(4.7.2)$$
Let $\ARI^\Delta$ (resp.~$\overline{\ARI}^\Delta$) denote the space of 
rational-valued moulds $P$ in the $u_i$ such that $\Delta(P)\in \ARI^{pol}$, 
i.e.~such that the denominator of the rational function $P(u_1,\ldots,u_r)$ 
is ``at worst'' $(u_1+\cdots+u_r)u_1\cdots u_r$, and similarly let 
$\overline{\ARI}^\Delta=swap(\ARI^\Delta)$ denote the space of moulds in the
$v_i$ with denominator ``at worst'' 
$v_1(v_1-v_2)\cdots (v_{r-1}-v_r)v_r$. In general we indicate moulds having
the property that $\Delta(A)\in \ARI^{pol}$ with the superscript $\Delta$
(in either the $u_i$ or the $v_i$),
writing for example $\ARI_{al}^\Delta$ for the space of 
moulds in $\ARI^\Delta$ which are also alternal. 

The statements and proofs in this section are mostly drawn from
S.~Baumard's Ph.D. thesis; we reproduce them here since the thesis was not 
published. 
\vskip .3cm
\noindent {\bf Theorem 4.7.1.} [Baumard, Lemma 4.40] {\it The spaces
$\overline{\ARI}_{al}^\Delta$, $\overline{\ARI}^\Delta_{*circneut}$ and
$\ARI_{\underline{al}/\underline{al}}^\Delta$ are all 
closed under the $ari$-bracket.}
\vskip .3cm
\noindent {\bf Proof.} We first need the following useful lemma.
Recall that a mould $A\in \overline{\ARI}$ is said to be 
*circ-neutral if it is circ-neutral up to addition of a constant mould
(see \S 2.6.1 for circ-neutrality).
\vskip .3cm
\noindent {\bf Lemma 4.7.2.} {\it (i)} [Baumard, Lemma 4.39] {\it Let $M\in 
\overline{\ARI}_{al}$ and let $\Delta(M)$ denote the image of $M$ under 
the $\Delta$ operator as in (4.7.2).  Then for all $r>1$, 
$\Delta(M)$ satisfies the identity
$$\Delta(M)(0,v_2,\ldots,v_r)=\Delta(M)(v_2,\ldots,v_r,0).\eqno(4.7.3)$$
(ii) Let $M\in \overline{\ARI}_{*circneut}$. Then $\Delta(M)$
again satisfies (4.7.3).}
\vskip .2cm
\noindent {\bf Proof.} (i) Let $r>1$. We obtain (4.7.3) from the first 
alternality relation on $M$, which we write as
$$\eqalign{0&=\sum_{i=1}^r M(v_2,\ldots,v_i,v_1,v_{i+1},\ldots,v_r)\cr
&=\sum_{i=2}^{r-1} {{\Delta(M)(v_2,\ldots,v_i,v_1,v_{i+1},\ldots,v_r)}\over
{v_2(v_2-v_3)\cdots(v_i-v_1)(v_1-v_{i+1})\cdots (v_{r-1}-v_r)v_r}}\cr
&\ \ \ \ \ +{{\Delta(M)(v_1,\ldots,v_r)}\over{v_1(v_1-v_2)\cdots(v_{r-1}-v_r)v_r}}
+{{\Delta(M)(v_2,\ldots,v_r,v_1)}\over{v_2(v_2-v_3)\cdots(v_{r-1}-v_r)(v_r-v_1)v_1}}.}$$
Multiplying the right-hand side by $v_1$ and then setting $v_1=0$ kills
all the terms in the sum (since $v_1$ does not appear in any of
the denominators of those terms), leaving only
$$0={{\Delta(M)(0,v_2,\ldots,v_r)}\over{(-v_2)(v_2-v_3)\cdots(v_{r-1}-v_r)v_r}}
+{{\Delta(M)(v_2,\ldots,v_r,0)}\over{v_2(v_2-v_3)\cdots(v_{r-1}-v_r)v_r}}.$$
(ii) Let $M\in \overline{\ARI}_{*circneut}$ and $r>1$. The *circ-neutrality of 
$M$ means that there exists a constant $M_0$ such that
$$\eqalign{0&=M(v_1,\ldots,v_r)+M(v_2,\ldots,v_r,v_1)+\cdots+M(v_r,v_1,\ldots,v_{r-1})+rM_0\cr
&={{\Delta(M)(v_1,\ldots,v_r)}\over{v_1(v_1-v_2)\cdots(v_{r-1}-v_r)v_r}}
+{{\Delta(M)(v_2,\ldots,v_r,v_1)}\over{v_2(v_2-v_3)\cdots (v_r-v_1)v_1}}+\cdots
\cr
&\qquad \qquad \ \  +{{\Delta(M)(v_r,v_1,\ldots,v_{r-1})}\over{v_r(v_r-v_1)\cdots(v_{r-2}-v_{r-1})v_{r-1}}}+rM_0.}$$
Again multiplying this identity by $v_1$ and setting $v_1=0$ makes all but the
first two terms disappear, and these become
$${{\Delta(M)(0,v_2,\ldots,v_r)}\over{-v_2(v_2-v_3)\cdots(v_{r-1}-v_r)v_r}}
+{{\Delta(M)(v_2,\ldots,v_r,0)}\over{v_2(v_2-v_3)\cdots(v_{r-1}-v_r)v_r}}.$$
Since the two terms have the same denominator but opposite signs, this
is equivalent to (4.7.3), which concludes the proof.\hfill{$\square$}
\vskip .3cm
\noindent {\bf Proof of Theorem 4.7.1.} We first note that the statement for
$\ARI^\Delta_{\underline{al}/\underline{al}}$ follows easily from the statement 
for $\overline{\ARI}^\Delta_{al}$. Indeed, let $A,B\in 
\ARI^\Delta_{\underline{al}/\underline{al}}$. Since $swap(A)$ and $swap(B)$
lie in $\overline{\ARI}^\Delta_{al}$, under the assumption that 
$\overline{\ARI}^\Delta_{al}$ is closed under the $ari$-bracket, 
we have $ari(swap(A),swap(B))\in \overline{\ARI}^\Delta_{al}$.
Since $A$ and $B$ are push-invariant by 
Lemma 2.5.5 (this is where we use the evenness property in depth 1,
i.e.~the assumption that $A,B\in \ARI_{\underline{al}/\underline{al}}$ rather
than just $\ARI_{al/al}$), we know by Lemma 2.4.1 that 
$$ari(swap(A),swap(B)\bigr)=swap\cdot ari(A,B),$$
and thus $swap\bigl(ari(A,B)\bigr)\in \overline{\ARI}_{al}^\Delta$, 
which means that $ari(A,B)\in \ARI_{al/al}^\Delta$. Since $ari(A,B)$ has no 
depth 1 part it lies in $\ARI_{\underline{al}/\underline{al}}^\Delta$. 
\vskip .2cm
We now prove simultaneously that $\overline{\ARI}^\Delta_{al}$ and
$\overline{\ARI}^\Delta_{*circneut}$ are closed under the $ari$-bracket,
in two steps.  Let $A,B$ lie in either one of the two spaces.
\vskip .3cm
\noindent {\it Proof that $ari(A,B)\in \overline{\ARI}^\Delta$.}
The proof of this fact is identical for the two spaces, because it 
does not use the actual conditions
of alternality or *circ-neutrality but only the identity (4.7.3), which
holds for moulds $M$ in both spaces by Lemma 4.7.2.
We use the proof given in Baumard's thesis (\S 4.3.4). Let $A,B\in 
\overline{\ARI}_{al}^\Delta$. Since everything is additive, we may assume that 
$A$ is concentrated in a single depth $r$ and $B$ in a single depth $s$.
Since we know that $\overline{\ARI}_{al}$ is closed under the $ari$-bracket
(cf.~Proposition 2.5.2), we only need to ensure that 
$ari(A,B)\in\overline{\ARI}^\Delta$. For this, we study what poles
can occur in the separate terms $arit(A)\cdot B$, 
$arit(B)\cdot A$, $mu(A,B)$ and $mu(B,A)$ of $ari(A,B)$. Taking the
definition of $arit$ given in (2.2.5) and reducing it to the case where
$A$ is concentrated in depth $r$ and $B$ in depth $s$, we write it as
$$\eqalign{(arit&(A)\cdot B)(v_1,\ldots,v_{r+s})
=\sum_{0\le i<s} B(v_1,\ldots,v_i,v_{i+r+1},\ldots,v_{r+s})A(v_{i+1}\cdots v_{i+r}\rfloor)\cr
&\ \ \ -\sum_{0<i\le s} B(v_1,\ldots,v_i,v_{i+r+1},\ldots,v_{r+s})A(\lfloor v_{i+1}\cdots v_{i+r})\cr
&=B(v_{r+1}\cdots v_{r+s})A(v_1-v_{r+1},\ldots,v_r-v_{r+1})+\sum_{i=1}^{s-1}B(v_1,\ldots,v_i,v_{i+r+1},\ldots,v_{r+s})\cdot\cr
&\qquad\qquad \ \Biggl(A(v_{i+1}-v_{i+r+1},\ldots,v_{i+r}-v_{i+r+1})-A(v_{i+1}-v_i,\ldots,v_{i+r}-v_i\Biggr)\cr
&\qquad \qquad\qquad -B(v_1,\ldots,v_s)A(v_{s+1}-v_s,\ldots,v_{r+s}-v_s).}$$
\vskip .2cm
We rewrite this as 
$$\eqalign{\Delta&\bigl(arit(A)\cdot B\bigr)(v_1,\ldots,v_{r+s})
=\sum_{i=1}^{s-1} \Delta(S_i) \cr
&+{{v_1}\over{v_{r+1}(v_1-v_{r+1})}}\Delta(B)(v_{r+1},\ldots,v_{r+s})
\Delta(A)(v_1-v_{r+1},\ldots,v_r-v_{r+1})\cr
&\ \ +{{v_{r+s}}\over{v_s(v_{r+s}-v_s)}}\Delta(B)(v_1,\ldots,v_s)\Delta(A)
(v_{s+1}-v_s,\ldots,v_{r+s}-v_s),}\eqno(4.7.4)$$
where
$$\eqalign{\Delta(S_i)&={{(v_i-v_{i+1})\cdots (v_{i+r}-v_{i+r+1})}\over
{v_i-v_{i+r+1}}}\Delta(B)(v_1,\ldots,v_i,v_{i+r+1},\ldots,v_{r+s})\cdot\cr
&\ \ \Biggl(A(v_{i+1}-v_{i+r+1},\ldots,v_{i+r}-v_{i+r+1})-A(v_{i+1}-v_i,\ldots,
v_{i+r}-v_i)\Biggr)\cr
&={{(v_i-v_{i_1})\cdots (v_{i+r}-v_{i+r+1})}\over{v_i-v_{i+r+1}}}
B(v_1,\ldots,v_i,v_{i+r+1},\ldots,v_{r+s})\cdot\cr
&\ \ \ \ \ \Biggl({{\Delta(A)(v_{i+1}-v_{i+r+1},\ldots,v_{i+r}-v_{i+r+1})}\over{(v_{i+1}-v_{i+r+1})(v_{i+1}-v_{i+2})\cdots(v_{i+r-1}-v_{i+r})(v_{i+r}-v_{i+r+1})}}\cr
&\qquad\qquad \ -{{\Delta(A)(v_{i+1}-v_i,\ldots,v_{i+r}-v_i)}\over{(v_{i+1}-v_i)(v_{i+1}-v_{i+2})\cdots(v_{i+r-1}-v_{i+r})(v_{i+r}-v_i)}}\Biggr)}$$
$$\eqalign{&={{1}\over{v_i-v_{i+r+1}}}\Delta(B)(v_1,\ldots,v_i,
v_{i+r+1},\ldots,v_{r+s})\cdot \cr
&\ \ \ \ \ \Biggl({{v_i-v_{i+1}}\over{v_{i+1}-v_{i+r+1}}}\Delta(A)(v_{i+1}-v_{i+r+1},\ldots,v_{i+r}-v_{i+r+1}) \cr
&\ \ \ \ \ \ \ \ \ \ +{{v_{i+r}-v_{i+r+1}}\over{v_{i+r}-v_i}}\Delta(A)(v_{i+1}-v_i,\ldots,v_{i+r}-v_i)\Biggr).}\eqno(4.7.5)$$
The expression (4.7.5) shows that there are only three possible types
of poles in $\Delta\bigl(arit(A)\cdot B\bigr)(v_1,\ldots,v_{r+s})$:
\vskip .1cm\noindent
(i) the poles of the form ${{1}\over{v_i-v_{i+r+1}}}$;
\vskip .1cm\noindent
(ii) the poles of the form ${{1}\over{v_i-v_{i+r}}}$.
\vskip .1cm\noindent
(iii) the poles ${{1}\over{v_{r+1}}}$ and ${{1}\over{v_s}}$, which only appear 
in one term;
\vskip .3cm
\noindent {\it Poles of type (i).} The pole ${{1}\over{v_i-v_{i+r+1}}}$ appears
uniquely as a factor of the term $\Delta(S_i)$. We show that it is in fact 
compensated by the sum of two terms in $\Delta(A)$ appearing in 
$\Delta(S_i)$, i.e. that $v_i-v_{i+r+1}$ divides the sum
$${{v_i-v_{i+1}}\over{v_{i+1}-v_{i+r+1}}}\Delta(A)(v_{i+1}-v_{i+r+1},\ldots,v_{i+r}-v_{i+r+1})$$
$$+{{v_{i+r}-v_{i+r+1}}\over{v_{i+r}-v_i}}\Delta(A)(v_{i+1}-v_i,\ldots,v_{i+r}-v_i).\eqno(4.7.6)$$
To see this, we write $x=v_i=v_{i+r+1}$ and substitute this into (4.7.6),
obtaining
$${{x-v_{i+1}}\over{v_{i+1}-x}}\Delta(A)(v_{i+1}-x,\ldots,v_{i+r}-x)
+{{v_{i+r}-x}\over{v_{i+r}-x}}\Delta(A)(v_{i+1}-x,\ldots,v_{i+r}-x)$$
$$=-\Delta(A)(v_{i+1}-x,\ldots,v_{i+r}-x)
+\Delta(A)(v_{i+1}-x,\ldots,v_{i+r}-x)=0.$$
Thus there are no poles of type (i) in $\Delta\bigl(arit(A)\cdot B\bigr)$.
\vskip .2cm
\noindent {\it Poles of type (ii).} We consider the three cases
$i=1$, $2\le i\le s-1$ and $i=s$ separately. When $i=1$, the pole
${{1}\over{v_1-v_{r+1}}}$ is multiplied by 
$${{v_1}\over{v_{r+1}}}\Delta(B)(v_{r+1},v_{r+2},\ldots,v_{r+2})
\Delta(A)(v_1-v_{r+1},\ldots,v_r-v_{r+1})\qquad\qquad\qquad\qquad\qquad$$
$$-{{1}\over{v_1-v_{r+2}}}\Delta(B)(v_1,v_{r+2},\ldots,v_{r+s})\cdot
(v_{r+1}-v_{r+2})\Delta(A)(v_2-v_1,\ldots,v_{r+1}-v_1).$$
Setting $v_1=v_{r+1}=x$, this becomes
$$\Delta(B)(x,v_{r+2},\ldots,v_{r+2})
\Delta(A)(0,v_2-x,\ldots,v_r-x)\qquad\qquad\qquad\qquad\qquad$$
$$-{{1}\over{x-v_{r+2}}}\Delta(B)(x,v_{r+2},\ldots,v_{r+s})\cdot
(x-v_{r+2})\Delta(A)(v_2-x,\ldots,v_r-x,0)$$
$$=\Delta(B)(x,v_{r+2},\ldots,v_{r+2})\Biggl(
\Delta(A)(0,v_2-x,\ldots,v_r-x) -\Delta(A)(v_2-x,\ldots,v_r-x,0)\Biggr)$$
which is equal to $0$ by Lemma 4.7.2, so there are no poles of
type (ii) when $i=1$.  
When $i=s$, the pole ${{1}\over{v_s-v_{r+s}}}$ is multiplied by
$$-{{v_{r+s}}\over{v_s}}\Delta(B)(v_1,\ldots,v_s)\Delta(A)(v_{s+1}-v_s,
\ldots,v_{r+s}-v_s)$$
$$+{{1}\over{v_{s-1}-v_{r+s}}}\Delta(B)(v_1,\ldots,v_{s-1},v_{r+s})\cdot
(v_{s-1}-v_s)\Delta(A)(v_s-v_{r+s},\ldots,v_{r+s-1}-v_{r+s}).$$
Setting $v_s=v_{r+s}=x$ in this expression, we find
$$-\Delta(B)(v_1,\ldots,v_{s-1},x)\Delta(A)(v_{s+1}-x,,v_{s+2}-x,
\ldots,v_{r+s-1}-x,0)$$
$$+{{1}\over{v_{s-1}-x}}\Delta(B)(v_1,\ldots,v_{s-1},x)\cdot
(v_{s-1}-x)\Delta(A)(0,v_{s+1}-x,\ldots,v_{r+s-1}-x)$$
which is again equal to zero by Lemma 4.7.2, so there are no poles of
type (ii) with $i=s$. Finally, for $2\le i\le s-1$, the pole
${{1}\over{v_i-v_{i+r}}}$ comes from the two terms $\Delta(S_{i-1})$
and $\Delta(S_i)$; putting the factors from these two terms together,
the pole appears in front of the expression
$${{v_{i-1}-v_i}\over{v_{i-1}-v_{i+r}}}\Delta(B)(v_1,\ldots,v_{i-1},v_{i+r},
v_{i+r+1}, \ldots,v_{r+s})\Delta(A)(v_i-v_{i+r},v_{i+1}-v_{i+r},\ldots,v_{i-1+r}-v_{i+r})$$
$$-{{v_{i+r}-v_{i+r+1}}\over{v_i-v_{i+r+1}}}\Delta(B)(v_1,\ldots,v_{i-1},v_i,
v_{i+r+1},\ldots,v_{r+s}) \Delta(A)(v_{i+1}-v_i,\ldots,v_{i+r-1}-v_i,v_{i+r}-v_i).$$
Setting $v_i=v_{i+r}=x$, this reduces to
$$\Delta(B)(v_1,\ldots,v_{i-1},x,v_{i+r+1},
\ldots,v_{r+s})\Delta(A)(0,v_{i+1}-x,\ldots,v_{i-1+r}-x)\qquad\qquad$$
$$-\Delta(B)(v_1,\ldots,v_{i-1},x,
v_{i+r+1},\ldots,v_{r+s}) \Delta(A)(v_{i+1}-x,\ldots,v_{i+r-1}-x,0),$$
which is once again equal to zero thanks to Lemma 4.7.2. Thus we have shown
that $\Delta(arit(B)\cdot A)$ has no poles of type (ii).
\vskip .2cm
\noindent {\it Poles of type (iii).} It remains to consider the potential
poles from the terms ${{1}\over{v_{r+1}}}$ and ${{1}\over{v_s}}$.
These arise from the terms
$${{v_1}\over{v_{r+1}(v_1-v_{r+1})}}\Delta(B)(v_{r+1},\ldots,v_{r+s})\Delta(A)(v_1-v_{r+1},\ldots,
v_r-v_{r+1})$$
$$\ \ \ +{{v_{r+s}}\over{v_s(v_{r+s}-v_s)}}\Delta(B)(v_1,\ldots,v_s)
\Delta(A)(v_{s+1}-v_s,\ldots,v_{r+s}-v_s).\eqno(4.7.7)$$
In fact, these poles are real poles in $\Delta(arit(A)\cdot B)$; thus
symmetrically, there are real poles at ${{1}\over{v_r}}$ and 
${{1}\over{v_{s+1}}}$ in $\Delta(arit(B)\cdot A)$.  Since
$$ari(A,B)=arit(B)\cdot A-arit(A)\cdot B+lu(A,B),\eqno(4.7.8)$$
to show that $\Delta\bigl(ari(A,B)\bigr)$ has no poles, we show that poles 
at $v_r$, $v_s$, $v_{r+1}$ and $v_{s+1}$ in $\Delta\bigl(arit(B)\cdot A
-arit(A)\cdot B\bigr)$ are cancelled out by poles at the same places in 
$\Delta\bigl(lu(A,B)\bigr)$, and also that $\Delta(lu(A,B)\bigr)$ has no 
other poles. The expression for $\Delta\bigl(lu(A,B)\bigr)$ is given by
$$\eqalign{\Delta\bigl(&lu(A,B)\bigr)(v_1,\ldots,v_{r+s})=v_1(v_1-v_2)\cdots 
(v_{r-1}-v_r)(v_r-v_{r+1})\cdots (v_{r+s-1}-v_{r+s})v_{r+s}\cdot\cr
&\qquad \Biggl({{\Delta(A)(v_1,\ldots,v_r)}\over{v_1(v_1-v_2)\cdots (v_{r-1}-v_r)v_r}}{{\Delta(B)(v_{r+1},\ldots,v_{r+s})}\over{v_{r+1}(v_{r+1}-v_{r+2})\cdots(v_{r+s-1}-v_{r+s})v_{r+s}}}\cr
&\qquad\qquad -{{\Delta(B)(v_1,\ldots,v_s)}\over{v_1(v_1-v_2)\cdots (v_{s-1}-v_s)v_s}}
{{\Delta(A)(v_{s+1},\ldots,v_{s+r})}\over{v_{s+1}(v_{s+1}-v_{s+2})\cdots v_{r+s-1}-v_{r+s})v_{r+s}}}\Biggr)\cr
&={{v_r-v_{r+1}}\over{v_rv_{r+1}}}\Delta(A)(v_1,\ldots,v_r)\Delta(B)(v_{r+1},\ldots,v_{r+s})\cr
&\qquad\qquad \ \ -{{v_s-v_{s+1}}\over{v_sv_{s+1}}}\Delta(B)(v_1,\ldots,v_s)\Delta(A)(v_{s+1},\ldots,v_{r+s}).}\eqno(4.7.9)$$
This shows that the only poles of $\Delta\bigl(lu(A,B)\bigr)$ are indeed at
$v_r$, $v_{r+1}$, $v_s$ and $v_{s+1}$. It remains only to show that these
poles cancel out with the poles at the same places appearing in
$\Delta\bigl(arit(B)\cdot A-arit(A)\cdot B\bigr)$. Let us show this
first for the pole at $v_s=0$.  For this, we multiply (4.7.7) and (4.7.9) 
by $v_s$, set $v_s=0$ in the results, 
and compare them. From (the second line of) (4.7.7) we obtain the residue
$$\Delta(B)(v_1,\ldots,v_s)\Delta(A)(v_{s+1},\ldots,v_{r+s})$$
at $v_s=0$, and from (4.7.9) we obtain exactly the same expression (also from
the second line). Noting that $arit(A)\cdot B$ appears in 
(4.7.8) with a negative sign, this means that the pole at ${{1}\over{v_s}}$
cancels out between $-\Delta(arit(A)\cdot B)$ and $\Delta\bigl(lu(A,B)\bigr)$.
The pole ${{1}\over{v_r}}$ also cancels out in the same way, thanks to the
symmetry between $A$ and $B$. Let us check that the pole at 
${{1}\over{v_{r+1}}}$ also cancels out. Again we multiply (4.7.7)
and (4.7.9) by $v_{r+1}$ and then set $v_{r+1}$ equal to zero. In (4.7.7)
there remains
$$\Delta(B)(v_{r+1},\ldots,v_{r+s})\Delta(A)(v_1,\ldots,v_r),$$
and in (4.7.9) exactly the same expression. Thus the pole at $v_{r+1}$
cancels out, and again by the symmetry between $A$ and $B$, so does the
pole at $v_{s+1}$. So $\Delta\bigl(ari(A,B)\bigr)$ actually has no poles,
which proves that $ari(A,B)\in \overline{\ARI}^\Delta$.
\vskip .3cm
\noindent {\it Proof that $\overline{\ARI}^\Delta_{al}$ and
$\overline{\ARI}^\Delta_{*circneut}$ are closed.}
In order to complete the proof that $\overline{\ARI}^\Delta_{al}$ is closed
under the $ari$-bracket, we use the fact that $\overline{\ARI}_{al}$
is closed under the $ari$-bracket (cf.~Proposition 2.5.2); thus
$ari(A,B)\in \overline{\ARI}_{al}$, and since we proved in the first step
that $ari(A,B)\in \overline{\ARI}^\Delta$, we find that
$$ari(A,B)\in \overline{\ARI}_{al}\cap \overline{\ARI}^\Delta=\overline{\ARI}^\Delta_{al}$$
as desired.

Finally, to complete the proof that $\overline{\ARI}^\Delta_{*circneut}$ is
closed under the $ari$-bracket, we need to check that
$\overline{\ARI}_{*circneut}$ is closed under the $ari$-bracket.
Let $A,B\in \overline{\ARI}_{*circneut}$,
and let $A_0$ and $B_0$ denote the constant moulds such that
$A+A_0$ and $B+B_0$ are circ-neutral. By Proposition 2.6.1,
$$ari\bigl(A+A_0,B+B_0\bigr)\in \overline{\ARI}_{circneut}.$$
Since constant moulds are invariant under the swap, this is equal to
$$ari\bigl(swap(A),swap(B)\bigr)
+ari\bigl(A_0,swap(B)\bigr)
+ari\bigl(swap(A),B_0\bigr)
+ari\bigl(A_0,B_0\bigr).$$
But the $ari$-bracket of a constant mould with any mould $A$ is zero
(see (4.6.4)), so we have
$$ari\bigl(swap(A+A_0),swap(B+B_0)\bigr)=ari\bigl(swap(A),swap(B)\bigr)$$
and therefore
$$ari\bigl(swap(A),swap(B)\bigr)\in \overline{\ARI}_{circneut}
\subset \overline{\ARI}_{*circneut}.$$
This shows that $\overline{\ARI}_{*circneut}$ is closed under the $ari$-bracket.
Thus $A,B\in \overline{\ARI}^\Delta_{*circneut}$, we saw above that
$ari(A,B)\in \overline{\ARI}^\Delta$, and we now see that $ari(A,B)\in
\overline{\ARI}_{*circneut}$, so
$$ari(A,B)\in \overline{\ARI}^\Delta\cap \overline{\ARI}_{*circneut}=
\overline{\ARI}^\Delta_{*circneut},$$
completing the proof of Theorem 4.7.1.  \hfill{$\square$}
\vfill\eject
\centerline{\bf Chapter 5} 
\vskip .5cm
\centerline{\bf Elliptic mould theory}
\vskip .6cm
The sections of this section relate all the previous results on moulds
and the double shuffle Lie algebra to the elliptic situation. We first
study the properties of the action of the adjoint operator $Ad_{ari}(invpal)$.
Then we use it to define the elliptic double shuffle Lie algebra, and
investigate the elliptic double shuffle relations and the astonishing
connection between the elliptic double shuffle Lie algebra and the
associated graded of the usual double shuffle Lie algebra.
\vskip .6cm
\noindent {\bf \S 5.1. The operator $Ad_{ari}(invpal)$ and the denominator $\Delta$}
\vskip .3cm
The main goal of this section is to use \'Ecalle's second fundamental 
identity to prove that applying the operator $Ad_{ari}(invpal)$ to
double shuffle moulds (i.e.~moulds in $\ARI^{pol}_{\underline{al}*\underline{il}}$) leads to denominators controlled by $\Delta$.
This result is again drawn from Baumard's thesis.
\vskip .3cm
\noindent {\bf Theorem 5.1.1.} [Baumard, Th\'eor\`eme 4.35] {\it Let $N\in \ARI^{pol}_{\underline{al}*\underline{il}}$ be a double shuffle mould in the $u_i$.  Then 
$$Ad_{ari}(invpal)\cdot N\in \ARI_{\underline{al}*\underline{al}}^\Delta.
\eqno(5.1.1)$$}\par
\noindent Before proving the theorem, we give two useful lemmas.
Recall the definition of the mould $pic$ given in (4.4.2).
\vskip .2cm
\noindent {\bf Lemma 5.1.2.} [Baumard, Lemme 4.37] {\it Let $poc$ be the mould
defined by $poc(\emptyset)=1$ and
$$poc(v_1,\ldots,v_r)={{1}\over{v_1(v_1-v_2)\cdots(v_{r-1}-v_r)}}$$
for $r\ge 1$.  Then $ganit_{pic}\circ ganit_{poc}=id$.}
\vskip .2cm
\noindent {\bf Proof.} The statement is equivalent to showing that
$mu(pic,ganit_{pic}\cdot poc)=1$ (where $1$ denotes the mould that takes 
value $1$ in depth $0$ and value $0$ in depths $r>0$). 
By definition, the inverse of $pic$ is the mould $1-V$ where $V$ is
defined by $V(\emptyset)=0$, $V(v_1)=1/v_1$ and $V(v_1,\ldots,v_r)=0$
for $r>1$. Direct calculation shows that $(ganit_{pic}\cdot poc)(\emptyset)=0$
and $(ganit_{pic}\cdot poc)(v_1)=1/v_1$, so it remains only to show that
$ganit_{pic}\cdot poc$ is zero in depths $r>1$.
Recall that $ganit_{pic}\cdot poc$ is defined by the formula
$$ganit_{pic}\cdot poc(v_1,\ldots,v_r) = 
\sum_{v_1\cdots v_r={\bf b}_1{\bf c}_1\cdots {\bf b}_s{\bf c}_s}
poc({\bf b}_1\cdots {\bf b}_s) pic(\lfloor{\bf c}_1)\cdots pic(\lfloor{\bf c}_s)\eqno(5.1.2)$$
where the sum runs over the set of decompositions of $v_1\cdots v_r$ into
$2s$ words of which only the last one ${\bf c}_s$ may be empty.

Let us define a map $\sigma$ from the set of decompositions 
${\bf b}_1{\bf c}_1\cdots {\bf b}_s$, i.e.~those having empty 
final part ${\bf c}_s$, to the set of decompositions in which ${\bf c}_s$
is non-empty.
\vskip .1cm
\noindent {\it Case 1}: if the final part ${\bf b}_s$ has length 1, we define
$$\sigma({\bf b}_1{\bf c}_1\cdots {\bf b}_{s-1}{\bf c}_{s-1}{\bf b}_s)={\bf b}_1{\bf c}_1\cdots
{\bf b}_{s-1}{\bf c}'_{s-1}$$
where ${\bf c}'_{s-1}={\bf c}_{s-1}{\bf b}_s$, i.e. we join up the single
letter ${\bf b}_s$ to the previous term ${\bf c}_{s-1}$.
\vskip .1cm
\noindent {\it Case 2}: if the final part ${\bf b}_s$ has length $>1$, we
break up ${\bf b}_s$ into two pieces ${\bf b}'_s{\bf c}'_s$ where ${\bf c}_s$ 
consists only of the final letter of ${\bf b}_s$, and set
$$\sigma({\bf b}_1{\bf c}_1\cdots {\bf c}_{s-1}{\bf b}_s)={\bf b}_1{\bf c}_1
\cdots {\bf c}_{s-1}{\bf b}'_s{\bf c}'_s.$$
\vskip .1cm
Since every decomposition with non-empty final part ${\bf c}_s$ is the
image under $\sigma$ of a unique decomposition with empty final part
(those having single-letter ${\bf c}_s$ coming from case 2 and those
with longer ${\bf c}_s$ from case 1), we see that $\sigma$ is a bijection
which pairs up terms of the two types. We will show that each pair of
terms cancels out in the sum (5.1.2).
For this, let us first compute the two corresponding terms in case 1,
where ${\bf b}_s$ consists of the single letter $v_r$. The corresponding
terms in the sum (5.1.2) are given by
$$poc({\bf b}_1\cdots {\bf b}_{s-1}{\bf b}_s)pic(\lfloor {\bf c}_1)
\cdots pic(\lfloor {\bf c}_{s-2})pic(\lfloor {\bf c}_{s-1}),$$
and the term corresponding to the decomposition $sigma({\bf b}_1{\bf c}_1\cdots
{\bf b}_s)$ is given by
$$poc({\bf b}_1\cdots {\bf b}_{s-1})pic(\lfloor {\bf c}_1)
\cdots pic(\lfloor {\bf c}_{s-2})pic(\lfloor{\bf c}_{s-1}v_r).$$
Comparing these two terms we see that letting $v_k$ denote the final letter of 
${\bf b}_{s-1}$, and writing ${\bf b}_s=v_r$, we have
$$poc({\bf b}_1\cdots {\bf b}_{s-1}{\bf b}_s)=poc({\bf b}_1\cdots {\bf b}_{s-1}){{1}\over{(v_k-v_r)}}$$
and
$$pic(\lfloor{\bf c}_{s-1}v_r)=pic(\lfloor{\bf c}_{s-1}){{1}\over{(v_r-v_k)}},$$
so they cancel out. Similarly, since in case 2 we have 
${\bf b}_s={\bf b}'_sv_r$ and ${\bf c}_s=v_r$, the two terms for a pair 
are given by
$$poc({\bf b}_1\cdots {\bf b}_{s-1}{\bf b}'_sv_r)pic(\lfloor{\bf c}_1)\cdots
pic(\lfloor{\bf c}_{s-1})$$
and
$$poc({\bf b}_1\cdots {\bf b}_{s-1}{\bf b}'_s)pic(\lfloor{\bf c}_1)\cdots
pic(\lfloor{\bf c}_{s-1})pic(\lfloor v_r),$$
but since $v_{r-1}$ is the last letter of ${\bf b}'_s$, we have
$$poc({\bf b}_1\cdots {\bf b}'_sv_r)=
poc({\bf b}_1\cdots {\bf b}'_s){{1}\over{v_{r-1}-v_r}}$$
and
$$pic(\lfloor{\bf c}_1)\cdots pic(\lfloor{\bf c}_{s-1})pic(\lfloor v_r)=
pic(\lfloor{\bf c}_1)\cdots pic(\lfloor{\bf c}_{s-1}){{1}\over{v_r-v_{r-1}}},$$
so again these two terms cancel in the sum (5.1.2), proving that it is equal
to zero for $r>1$. This completes the proof of Lemma 5.1.2\hfill{$\square$}
\vskip .3cm
\noindent {\bf Lemma 5.1.3.} [Baumard, Lemme 4.38] {\it Let $A\in \overline{\ARI}^{pol}$. Then
$$swap\cdot ganit_{poc}\cdot A\in \ARI^\Delta.\eqno(5.1.3)$$}
\noindent Proof. The explicit expression for $ganit$ in (2.8.2)
shows that the only denominators that can occur
in $ganit_{poc}\cdot A$ come from the factors
$$poc(\lfloor {\bf b}_1)\cdots poc(\lfloor {\bf b}_s)$$
for all decompositions $d_{\bf v}={\bf a}_1{\bf b}_1\cdots
{\bf a}_s{\bf b}_s$ of ${\bf v}=(v_1,\ldots,v_r)$ into chunks, where if
the chunk ${\bf b}_i$ is given by $(v_k,\ldots,v_{k+l})$, then
$$\lfloor {\bf b}_i=(v_k-v_{k-1},v_{k+1}-v_{k-1},\ldots,v_{k+l}-v_{k-1})$$
(note that ${\bf a}_1\ne\emptyset$ and therefore $k>1$).  By the definition of $poc$, the only factors that can appear are
$(v_l-v_{l-1})$ where $v_l$ is a letter in one of ${\bf b}_i$, and these
factors appear in each term with multiplicity one.  Since
the sum ranges over all possible decompositions, the only letter of
${\bf v}$ that never belongs to any ${\bf b}_i$ is $v_1$, so the
factor $(v_r-v_1)$ never appears but all the other
factors $(v_i-v_{i-1})$ for $1<i\le r$ do appear.  
Thus $(v_1-v_2)(v_2-v_3)\cdots (v_{r-1}-v_r)$
is a common denominator for all the terms in the sum defining
$ganit_{poc}\cdot A$.  The swap of this common denominator is
$u_2\cdots u_r$, so this term is a common denominator for
$swap\cdot ganit_{poc}\cdot A\in \ARI$ and thus $\Delta\bigl(swap\cdot ganit_{poc}\cdot A\bigr)\in \ARI^{pol}$, proving the result.\hfill{$\square$}
\vskip .3cm
\noindent {\bf Proof of Theorem 5.1.1.} Let $M$ be a mould in the $u_i$ which 
is $push$-invariant and let $N=Ad_{ari}(pal)\cdot M$, 
i.e.~$M=Ad_{ari}(invpal)\cdot N$.  Then \'Ecalle's second fundamental identity 
(4.5.2) can be rewritten in terms of $N$ as follows:
$$swap\cdot Ad_{ari}(invpil)\cdot ganit_{poc}\cdot swap(N)=Ad_{ari}(invpal)\cdot N.\eqno(5.1.4)$$
We saw in Theorem 4.6.1 that if $N$ is as in the statement of the theorem,
then $M\in \ARI_{\underline{al}*\underline{al}}$, and therefore by Lemma
2.5.5, $M$ is push-invariant, so (4.7.7) holds. It remains only to prove
that the denominators of $M$ are controlled by $\Delta$.

Applying Lemma 5.1.3 with $A=swap(N)\in \overline{\ARI}^{pol}$ shows that
$ganit_{poc}\cdot swap(N)\in\overline{\ARI}^\Delta$.
By Theorem 4.7.1, the space $\overline{\ARI}_{al}^\Delta$ is closed
under the $ari$-bracket. Let us show that this space is preserved by
the operator $Ad_{ari}(invpil)$ (Corollaire 4.41 of Baumard's thesis).
Let $f(x)=1-e^{-x}$, and recall the sequence of moulds $re_r$ for 
$r\ge 1$ defined in (4.1.3) and the mould $lop_f$ defined in (4.1.4).
By (4.1.5) and the definition just following Prop. 4.1.2, we have
$$pil=exp_{ari}(lop_f),\ \ {\rm so}\ \ invpil=exp_{ari}(-lop_f).\eqno(5.1.5)$$
Since $re_1\in \overline{\ARI}_{al}^\Delta$ and this space
is closed under the $ari$-bracket by Theorem 4.7.1, all the moulds 
$re_r$ lie in this space and therefore $\pm lop_f\in \overline{\ARI}_{al}^\Delta$.
By definition, we have the equality of operators 
$$Ad_{ari}(invpil)=exp\bigl(ad_{ari}(-lop_f)\bigr),$$
where $ad_{ari}(P)$ is the Lie adjoint operator, i.e.~$ad_{ari}(P)\cdot Q=ari(P,Q)$.  Thus we can write
$$Ad_{ari}(invpil)\cdot ganit_{poc}\cdot swap(N)=\sum_{n=0}^\infty {{(-1)^n}\over{n!}} ad_{ari}(lop_f)^n\cdot \bigl(ganit_{poc}\cdot swap(N)\bigr).\eqno(5.1.6)$$
Since $ganit_{poc}\cdot swap(N)$ and $lop_f$ are both in 
$\overline{\ARI}^\Delta$, the fact that $\overline{\ARI}^\Delta$ is closed
under the $ari$-bracket shows that each term in this sum lies
in $\overline{\ARI}^\Delta$. Thus
$$Ad_{ari}(invpil)\cdot ganit_{poc}\cdot swap(N)\in \overline{\ARI}^\Delta,$$
and taking the swap of this mould and using (5.1.4) then proves (5.1.1),
completing the proof of Theorem 5.1.1.\hfill{$\square$}
\vskip .5cm
\noindent {\bf \S 5.2. $\Delta$ as a Lie algebra isomorphism, and the
$Dari$-bracket}
\vskip .3cm
The goal of this section is to study the transport of the $ari$-bracket
by the linear isomorphism $\Delta$. We define a new Lie bracket $Dari$ on
the vector space $\ARI$ by
$$Dari(A,B)=\Delta\bigl(ari(\Delta^{-1}(A),\Delta^{-1}(B)\bigr),\eqno(5.2.1)$$
so that writing $\ARI_{ari}$ for the space $\ARI$ equipped with the
$ari$-bracket and $\ARI_{Dari}$ for the space equipped with the $Dari$-bracket,
$\Delta$ gives a Lie algebra isomorphism
$$\Delta:\ARI_{ari}\buildrel\sim\over\rightarrow \ARI_{Dari}.\eqno(5.2.2)$$
In this section we give several properties of the Lie bracket $Dari$;
in particular, like $ari$, $Dari$ can be interpreted as a bracket of 
a certain type of derivation. The results in this section are all drawn from
[S2].
\vskip .3cm
\noindent {\bf Proposition 5.2.1.} {\it For every $A\in \ARI$, let $Darit(A)$
denote the operator on $\ARI$ defined by
$$Darit(A)=-dar\circ \Bigl(arit\bigl(\Delta^{-1}(P)\bigr)-ad_{ari}\bigl(\Delta^{-1}(P)\bigr)\Bigr)\circ dar^{-1}.\eqno(5.2.3)$$
Then $Darit(A)$ is a derivation of $\ARI_{lu}$, and 
$$Dari(A,B)=Darit(A)\cdot B-Darit(B)\cdot A.\eqno(5.2.4)$$}
\vskip .1cm
\noindent {\bf Proof.} It is clear that $Darit(A)$ is a derivation of 
$\ARI_{lu}$ since both $arit(A)$ and $ad_{ari}(A)$ are. 
$$\eqalign{Darit(A)\cdot &B-Darit(B)\cdot A
=-\bigl(dar\circ arit(\Delta^{-1}A)\circ dar^{-1}\bigr)\cdot B+
\bigl(dar\circ ad(\Delta^{-1}A)\circ dar^{-1}\bigr)\cdot B\cr
&\qquad\qquad +\bigl(dar\circ arit(\Delta^{-1}B)\circ dar^{-1}\bigr)\cdot A
-\bigl(dar\circ ad(\Delta^{-1}B)\circ dar^{-1}\bigr)\cdot A\cr
&=-\bigl(\Delta\circ arit(\Delta^{-1}A)\circ \Delta^{-1}\bigr)\cdot B
+\bigl(\Delta\circ arit(\Delta^{-1}B)\circ \Delta^{-1}\bigr)\cdot A\cr
&\qquad\qquad +\bigl(dar\circ ad(\Delta^{-1}A)\circ dar^{-1}\bigr)\cdot B
-\bigl(dar\circ ad(\Delta^{-1}B)\circ dar^{-1}\bigr)\cdot A\cr
&=-\bigl(\Delta\circ arit(\Delta^{-1}A)\circ \Delta^{-1}\bigr)\cdot B
+\bigl(\Delta\circ arit(\Delta^{-1}B)\circ \Delta^{-1}\bigr)\cdot A\cr
&\qquad\qquad +dar\bigl([\Delta^{-1}(A),dar^{-1}B]\bigr)-
dar\bigl([\Delta^{-1}(A),dar^{-1}A]\bigr)\cr
&=\Delta\Bigl(-arit(\Delta^{-1}A\cdot \Delta^{-1}B+arit(\Delta^{-1}B)\cdot
\Delta^{-1}A\cr
&\qquad\qquad +dur^{-1}\bigl([\Delta^{-1}A,dar^{-1}B]+[dar^{-1}A,\Delta^{-1}B]\bigr)\Bigr)\cr
&=\Delta\Bigl(-arit(\Delta^{-1}A\cdot \Delta^{-1}B+arit(\Delta^{-1}B)\cdot
\Delta^{-1}A\cr
&\qquad\qquad +dur^{-1}\bigl([\Delta^{-1}A,dur\Delta^{-1}B]+[dur\Delta^{-1}A,\Delta^{-1}B]\bigr)\Bigr)\cr
&=\Delta\Bigl(-arit(\Delta^{-1}A\cdot \Delta^{-1}B+arit(\Delta^{-1}B)\cdot
\Delta^{-1}A\cr
&\qquad\qquad +dur^{-1}dur\bigl([\Delta^{-1}A,\Delta^{-1}B]\bigr)\Bigr)\cr
&=\Delta\Bigl(-arit(\Delta^{-1}A\cdot \Delta^{-1}B+arit(\Delta^{-1}B)\cdot
\Delta^{-1}A+[\Delta^{-1}A,\Delta^{-1}B]\Bigr)\cr
&=\Delta\bigl(ari(\Delta^{-1}A,\Delta^{-1}B)\bigr)\cr
&=Dari(A,B).}$$
This completes the proof.\hfill{$\square$}
\vskip .3cm
Recall from Definition 1.3.2 that $L\subset {\rm Lie}[a,b]$ denotes the
(degree-completed) Lie subalgebra ${\rm Lie}[C_1,C_2,\ldots]$, 
where $C_i=ad(a)^{i-1}(b)$ (with variables $a,b$ instead of $x,y$), so 
that $ma:L\rightarrow \ARI^{pol}_{al}$ is an isomorphism. 
\vskip .8cm
\noindent {\bf \S 5.3. Adding the mould $a$ to $\ARI$}
\vskip .3cm
Recall the mould operators 
$$\cases{dur(Q)(u_1,\ldots,u_r)=(u_1+\cdots+u_r)Q(u_1,\ldots,u_r)\cr
dar(Q)(u_1,\ldots,u_r)=u_1\cdots u_r\,Q(u_1,\ldots,u_r)\cr
\Delta(Q)=dur(dar(Q)).}\eqno(5.3.1)$$
In this section we consider the free Lie algebra on two 
non-commutative variables $a$ and $b$, which we differentiate from
${\rm Lie}[x,y]$ by considering ${\rm Lie}[x,y]$ as the Lie algebra
of the fundamental group of the thrice-punctured sphere, and
${\rm Lie}[a,b]$ as the fundamental group of the once-punctured
torus. All the results presented in this section are drawn from [S2].

Let $\ARI^a$ denote the vector space spanned by $\ARI$ and by one
further mould, denoted $a$, which takes value $a$ in depth $0$ 
and $0$ in all other depths. We extend the $lu$ bracket to
$\ARI^a$ by setting 
$$lu(P,a)=dur(P)\eqno(5.3.2)$$
for all $P\in \ARI_{lu}$.
We write $\ARI^a_{lu}$ for the vector space $\ARI^a$
viewed as a Lie algebra under the $lu$-bracket. Note in particular
that if $P=ma(p)$ for a Lie series $p$ in the Lie subalgebra of
${\rm Lie}[a,b]$ generated by $ad(a)^{i-1}(b)$ for $i\ge 1$, then
from Lemma 3.3.1, we have
$$ma([p,a])=dur(P),\eqno(5.3.3)$$
so adding the mould $a$ to $\ARI$ gives us an injective linear morphism
$${\rm Lie}[a,b]\hookrightarrow \ARI^a$$
from the (degree-completed) Lie algebra on independent variables
$a,b$ to $\ARI$,
mapping the variable $a$ to the mould denoted by $a$. 
We thus obtain a Lie algebra isomorphism
$${\rm Lie}[a,b]\buildrel\sim\over\rightarrow (\ARI^a_{lu})^{pol}_{al}.\eqno(5.3.4)$$
Note that by (5.3.2), all moulds in $\ARI^a$ are of the form $ca+P$ for
$P\in \ARI$ and a scalar $c$ in the base field, which here we take to be 
${\Bbb Q}$. 
\vskip .3cm
The following proposition shows how to extend all the
derivations on $\ARI_{lu}$ that we need to the space $\ARI^a_{lu}$.
\vskip .3cm
\noindent {\bf Proposition 5.3.1.} {\it (i) The automorphism $dar$ 
extends to $a$
taking the value $dar(a)=a$;
\vskip .1cm
\noindent (ii) The derivation $dur$ extends to $a$ taking the value $dur(a)=0$;
\vskip .1cm
\noindent (iii) For all $P\in \ARI$, the derivation
$arit(P)$ of $\ARI_{lu}$ extends to $a$, taking the value $arit(P)\cdot a=0$.
\vskip .1cm
\noindent (iv) For all $P\in \ARI$, the derivation $Darit(P)$
of $\ARI_{lu}$ extends to $a$, with $Darit(P)\cdot a=P$. Furthermore,
$Darit(P)\cdot B_1=0$.}
\vskip .2cm
\noindent {\bf Proof.} Since $dar$ is an automorphism, to check (5.3.2) 
we write
$$lu(dar(Q),dar(a))=lu(dar(Q),a)=dur\bigl(dar(Q)\bigr).$$ 
But it is obvious from their definitions that $dur$ and $dar$ commute, 
so this is indeed equal to $dar\bigl(dur(Q)\bigr)$. This proves (i).
We check (5.3.2) for (ii) similarly. Because $dur(a)=0$ and
$dur$ is a derivation, we have
$$dur\bigl(lu(Q,a)\bigr)=lu\bigl(dur(Q),a\bigr)=dur\bigl(dur(Q)\bigr).$$
For (iii), we have
$$arit(P)\cdot lu(Q,a)=lu\bigl(arit(P)\cdot Q,a\bigr)=dur\bigl(arit(P)\cdot Q)\bigr).$$
But as pointed out by Ecalle [E2] (cf. [S, Lemma 4.2.4] for details), 
$arit(P)$ commutes with $dur$ for all $P$, which proves the result. 
\vskip .1cm\noindent
For (iv), the calculation to check that (5.3.2) is respected is a little
more complicated.  Let $Q\in \ARI$. Again using the
commutation of $arit(P)$ with $dur$, as well as that of $dar$ and $dur$,
we compute
$$\eqalign{Darit&(P)\cdot lu(Q,a)=lu\bigl(Darit(P)(Q),a\bigr)+
lu\bigl(Q,Darit(P)(a)\bigr)\cr
&=dur\bigl(Darit(P)\cdot Q\bigr)+lu(Q,P)\cr
&=-dur\Biggl(dar\Bigl(arit\bigl(\Delta^{-1}(P)\bigr)\cdot dar^{-1}(Q)-lu\bigl(\Delta^{-1}(P),dar^{-1}(Q)\bigr)\Bigr)\Biggr)+lu(Q,P)\cr
&=-dur\Biggl(dar\Bigl(arit\bigl(\Delta^{-1}(P)\bigr)\cdot dar^{-1}(Q)\Bigr)\Biggr)-dur\Bigl(lu\bigl(Q,dur^{-1}(P)\bigr)\Bigr)+lu(Q,P)\cr
&=-dar\Biggl(dur\Bigl(arit\bigl(\Delta^{-1}(P)\bigr)\cdot dar^{-1}(Q)\Bigr)\Biggr)-lu\bigl(lu(Q,N),a\bigr)+lu\bigl(Q,lu(N,a)\bigr)\cr
&\qquad\qquad\qquad\qquad\qquad\qquad{\rm with\ } N=dur^{-1}P,\ {\rm i.e.,}\ P=lu(N,a)\cr}$$
$$\eqalign{\qquad\qquad 
&=-dar\Bigl(arit\bigl(\Delta^{-1}(P)\bigr)\cdot dur\,dar^{-1}(Q)\Bigr)
-lu\bigl(lu(Q,a),N\bigr)\ \ {\rm by\ Jacobi}\cr 
&=-dar\Bigl(arit\bigl(\Delta^{-1}(P)\bigr)\cdot dar^{-1}\,dur(Q)\Bigr)
-lu\bigl(dur(Q),dur^{-1}P\bigr)\cr
&=-dar\Bigl(arit\bigl(\Delta^{-1}(P)\bigr)\cdot dar^{-1}\,dur(Q)\Bigr)
-dar\Bigl(lu\bigl(dar^{-1}dur(Q),dar^{-1}dur^{-1}(P)\bigr)\Bigr)\cr
&=-dar\Bigl(arit\bigl(\Delta^{-1}(P)\bigr)\cdot dar^{-1}\,dur(Q)\Bigr)
+dar\Bigl(lu\bigl(\Delta^{-1}(P),dar^{-1}dur(Q)\bigr)\Bigr)\cr
&=Darit(P)\cdot dur(Q).}$$
This proves the first statement of (iv).  For the second statement, we note that
$dar^{-1}(B_1)=B$.  Set $R=\Delta^{-1}(P)$, and we compute
$$\eqalign{Darit(P)\cdot B_1&=
-dar\bigl(arit(R)\cdot B\bigr)+dar\bigl([R,B]\bigr)\cr
&=-u_1\cdots u_r \bigl(R(u_1,\ldots,u_{r-1})-R(u_2,\ldots,u_r)\bigr)\cr
&\qquad -u_1\cdots u_r\bigl(R(u_2,\ldots,u_r)-R(u_1,\ldots,u_{r-1})\bigr)\cr
&=0.}$$
This concludes the proof of Proposition 5.3.1.\hfill{$\square$}
\vskip .3cm
\noindent {\bf Definition 5.3.2.} For any mould $P\in \ARI$, we define the
{\it partner} $P'$ of $P$ by
$$P'(u_1,\ldots,u_r)={{1}\over{u_1+\cdots+u_r}}\Bigl(P(u_2,\ldots,u_{r-1},u_r)-P(u_2,\ldots,u_r)\Bigr).\eqno(5.3.5)$$
Observe that we have the equality
$$P'(u_1,\ldots,u_r)={{1}\over{u_1+\cdots+u_r}}\Bigl(P(u_2,\ldots,u_{r-1},-u_1-\cdots-u_{r-1})-P(u_2,\ldots,u_r)\Bigr)\eqno(5.3.6)$$
if and only if $P$ is push-invariant, in which case $P'\in \ARI^{pol}$.

For any Lie series $p\in L$, let $p'$ be the (not necessarily
Lie) power series associated to $p$ by the formula
$$p'=\sum_{i\ge 0} {{(-1)^{i-1}}\over{i!}}a^ib\,\partial^i_a(p_a)\eqno(5.3.7)$$
where we write $p=p_aa+p_bb$ and $\partial_a$ denotes the derivation
of ${\rm Lie}[a,b]$ defined by $\partial_a(a)=1$, $\partial_a(b)=0$.
We also call $p'$ the {\it partner} of $p$.
\vskip .3cm
\noindent {\bf Lemma 5.3.3.} {\it Let $B=ma(b)$; it is the mould 
concentrated in depth 1 given by $B(u_1)=1$.  Then the derivation 
$Darit(P)$ of $\ARI^a_{lu}$ associated
to any mould $P\in \ARI$ by (5.2.3) satisfies
$$Darit(P)\cdot a=P,\ \ \ Darit(P)\cdot B=P',\eqno(5.3.8)$$
where $P'$ denotes the partner of $P$ defined in (5.3.5). }
\vskip .2cm
\noindent {\bf Proof.}
Let us compute the mould $Darit(P)\cdot B$ using (5.2.3). First
we set $\tilde B=dar^{-1}B$; it is the mould concentrated in depth 1
defined by $\tilde B(u_1)=1/u_1$. By (2.2.5) we have
$$\eqalign{\Bigl(arit\bigl(\Delta^{-1}&(P)\bigr)\cdot \tilde B\Bigr)
(u_1,\ldots,u_r)
=\tilde B(u_1+\cdots +u_r)\Bigl(\Delta^{-1}(P)(u_1,\ldots,u_{r-1})-\Delta^{-1}(P)(u_2,\ldots,u_r)\Bigr)\cr
&=\tilde B(u_1+\cdots +u_r)\Bigl({{P(u_1,\ldots,u_{r-1})}\over{u_1\cdots u_{r-1}(u_1+\cdots+u_{r-1})}}-{{P(u_2,\ldots,u_r)}\over
{u_2\cdots u_r(u_2+\cdots+u_r)}}\Bigr)\cr
&={{1}\over{(u_1+\cdots +u_r)}}\Bigl({{P(u_1,\ldots,u_{r-1})}\over{u_1\cdots u_{r-1}(u_1+\cdots+u_{r-1})}}-{{P(u_2,\ldots,u_r)}\over
{u_2\cdots u_r(u_2+\cdots+u_r)}}\Bigr).}\eqno(5.3.9)$$
We also have
$$\eqalign{\Bigl(-ad_{ari}&\bigl(\Delta^{-1}(P)\bigr)\cdot \tilde B\Bigr)
(u_1,\ldots,u_r)=\tilde B(u_1)\Delta^{-1}(P)(u_2,\ldots,u_r)
-\Delta^{-1}(P)(u_1,\ldots,u_{r-1})\tilde B(u_r)\cr
&=\tilde B(u_1){{P(u_2,\ldots,u_r)}\over{u_2\cdots u_r(u_2+\cdots+u_r)}}
-{{P(u_1,\ldots,u_{r-1})}\over{u_1\cdots u_{r-1}(u_1+\cdots+u_{r-1})}}\tilde B(u_r)\cr
&={{P(u_2,\ldots,u_r)}\over{u_1\cdots u_r(u_2+\cdots +u_r)}}
-{{P(u_1,\ldots,u_{r-1})}\over{u_1\cdots u_r(u_1+\cdots+u_{r-1})}}.}
\eqno(5.3.10)$$
Applying $dar$ to the sum of (5.3.9) and (5.3.10) yields
$$\bigl(Darit(P)\cdot B\bigr)(u_1,\ldots,u_r)=
{{1}\over{(u_1+\cdots +u_r)}}\Bigl({{u_rP(u_1,\ldots,u_{r-1})}\over{(u_1+\cdots+u_{r-1})}}-{{u_1P(u_2,\ldots,u_r)}\over
{(u_2+\cdots+u_r)}}\Bigr)$$
$$+{{P(u_2,\ldots,u_r)}\over{(u_2+\cdots +u_r)}}
-{{P(u_1,\ldots,u_{r-1})}\over{(u_1+\cdots+u_{r-1})}}$$
$$={{P(u_2,\ldots,u_r)}\over{(u_2+\cdots +u_r)}}
\Bigl(1-{{u_1}\over{u_1+\cdots+u_r}}\Bigr)
-{{P(u_1,\ldots,u_{r-1})}\over{(u_1+\cdots+u_{r-1})}}
\Bigl(1-{{u_r}\over{u_1+\cdots+u_r}}\Bigr)$$
$$={{P(u_2,\ldots,u_r)}\over{(u_2+\cdots +u_r)}}
\Bigl({{u_2+\cdots+u_r}\over{u_1+\cdots+u_r}}\Bigr)
-{{P(u_1,\ldots,u_{r-1})}\over{(u_1+\cdots+u_{r-1})}}
\Bigl({{u_1+\cdots+u_{r-1}}\over{u_1+\cdots+u_r}}\Bigr)$$
$${{1}\over{u_1+\cdots+u_r}}\Bigl(P(u_2,\ldots,u_r)-P(u_1,\ldots,
u_{r-1})\Bigr).\eqno(5.3.11)$$
This is equal to the definition of the partner $P'$ of $P$ by (5.3.5).
\hfill{$\square$}
\vskip .3cm
\noindent {\bf Theorem 5.3.4.} {\it Let $p$ be a Lie series in 
${\rm Lie}[a,b]$ with no linear term, let
$p'$ be its partner as in (5.3.6), 
let $P=ma(p)$ and let $P'$ denote the partner of $P$ as in (5.3.5).
Let $E_p$ denote the derivation of ${\rm Lie}[a,b]$ defined by 
$E_p(a)=p$, $E_p(b)=p'$. Then 
\vskip .2cm\noindent
(i) $E_p([a,b])=0$ if and only if $p$ is push-invariant. Furthermore if
$p$ is push-invariant with no linear term, then $p'$ is the unique 
Lie series with no linear term such that $E_p([a,b])=0$.
\vskip .2cm\noindent
\noindent (ii) The derivation $Darit(P)$ 
restricted to $ma\bigl({\rm Lie}[a,b]\bigr)=(\ARI^a_{lu})^{pol}_{al}$ 
is the mould version of the derivation $E_p$, meaning that
for all $\in (\ARI^a_{lu})^{pol}_{al}$ we have
$$ma\bigl(E_p(q)\bigr)=Darit(P)\cdot ma(q),\eqno(5.3.12)$$
if and only if $p$ is push-invariant.
\vskip .2cm\noindent
(iii) We have the equality 
$$P'=ma(p')\eqno(5.3.13)$$ 
if and only if $p$ is push-invariant.
Furthemore if $P=ma(p)$ is a polynomial, alternal and push-invariant
mould then the partner $P'$ of $P$ defined in (5.3.5) is also
alternal and polynomial.
\vskip .1cm\noindent
(iv) If $p,q$ are two push-invariant Lie
series in $L$ and $p'$ and $q'$ denote their partners as
defined in (5.3.6), and $E_p$ and $E_q$ the associated derivations
of ${\rm Lie}[a,b]$, then setting $P=ma(p)$ and $Q=ma(Q)$, the bracket of 
$E_p$ and $E_q$ is related to the $Dari$-bracket of $P$ and $Q$ by
$$Dari(P,Q)=ma\bigl([E_p,E_q](a)\bigr).\eqno(5.3.14)$$}
\vskip .3cm
\noindent {\bf Proof.} (i) This is proven in [S1]; more precisely it 
is the equivalence between parts (ii) and (iv) of Theorem 2.1 there. 
\vskip .2cm\noindent
(ii) By Proposition 5.3.1 (iv), we have 
$$P=ma(p)=ma\bigl(E_p(a)\bigr)=Darit(P)\cdot a.\eqno(5.3.15)$$
By Proposition 5.3.1 (iv) we have
$Darit(P)\cdot ma([a,b])=0$, and by (i) above we have
$ma\bigl(E_p([a,b])\bigr)=0$ if and only if $p$ is push-invariant.
Thus $Darit(P)$ agrees with $E_p$ on $a$ and $[a,b]$ if and only
if $p$ is push-invariant. But derivations of ${\rm Lie}[a,b]$ 
which annihilate $[a,b]$ are determined by their value on $a$,
so $Darit(P)$ must agree with the mould version of $E_p$ on all
of $ma\bigl({\rm Lie}[a,b]\bigr)$ as expressed in (5.3.12).
\vskip .2cm\noindent
(iii) Let $B=ma(b)$. We saw in Lemma 5.3.3 that
$Darit(P)\cdot B=P'$. By (ii), $Darit(P)$ coincides with $E_p$ on
$ma\bigl({\rm Lie}[a,b]\bigr)$ if and only if $p$ is push-invariant, 
in which case we have
$$Darit(P)\cdot B=P'=ma\bigl(E_p(b)\bigr)=ma(p'),$$
proving (5.3.13). For the second statement, suppose $P=ma(p)$ is
polynomial, alternal and push-invariant.  Then by (i) there is
a unique derivation $E_p$ of ${\rm Lie}[a,b]$ mapping $a\mapsto p$
and annihilating $[a,b]$, so $p'=E_p(b)\in {\rm Lie}[a,b]$.
Thus since $P'=ma(p')$ by (5.3.13), $P'$ is alternal and polynomial.
\vskip .2cm\noindent
(iv) By (5.2.4), we have
$$Dari(P,Q)=Darit(P)\cdot Q-Darit(Q)\cdot P,$$
which by Proposition 5.3.1 (iv) we can write as
$$Dari(P,Q)=Darit(P)\cdot Darit(Q)\cdot a-Darit(Q)\cdot Darit(P)\cdot a
=[Darit(P),Darit(Q)]\cdot a.$$
But since $Darit(P)$ agrees with $E_p$ and $Darit(Q)$ agrees with
$E_q$ on $ma\bigl({\rm Lie}[a,b]\bigr)$, we have
$$ma\bigl([E_p,E_q](a)\bigr)=[Darit(P),Darit(Q)]\cdot a,$$
proving (iv). This concludes the proof of Theorem 5.3.4.\hfill{$\square$}
\vskip .5cm
\noindent {\bf \S 5.4. Closed subspaces of $\ARI_{Dari}$}
\vskip .5cm
The Lie morphism $\Delta$ and the $Dari$-bracket turn out to be useful in
proving results on $\ARI_{ari}$.
\vskip .3cm
\noindent {\bf Proposition 5.4.1.} {\it The space $\ARI^{pol}_{al+push}$ is
closed under the $Dari$-bracket, and the space $\ARI^\Delta_{al+push}$ is
closed under the $ari$-bracket.}
\vskip .3cm
\noindent {\bf Proof.} Let $p,q$ be push-invariant Lie series in $L$,
let $p'$ and $q'$ denote their partners as in (5.3.6), and let $E_p$ and
$E_q$ be the associated derivations of ${\rm Lie}[a,b]$. Then $E_p$ and
$E_q$ annihilate $[a,b]$, so the
bracket $[E_p,E_q]$ also annihilates $[a,b]$. Thus if we set $r=[E_p,E_q](a)$
and $r'=[E_p,E_q](b)$, then by Theorem 5.3.4 (i), $r$ is push-invariant
and $r'$ is its partner.  
Writing $L^{push}$ for the push-invariant and consider the injective map
$L^{push}\rightarrow {\rm Der}{\rm Lie}[a,b]$ defined by $p\mapsto E_p$.
Under this map, we can pull back the Lie bracket of derivations to a
Lie bracket on $L^{push}$, denoted $\langle .,.\rangle$, satisfying
$\langle p,q\rangle=r$.  In terms of moulds, since letting $P=ma(p)$ and
$Q=ma(Q)$ we have $ma([E_p,E_q](a))=Dari(P,Q)$ by (5.3.14), 
the fact that $L^{push}$ is preserved by $\langle .,.\rangle$ implies that
the $Dari$-bracket preserves $\ARI^{pol}_{al+push}$. Now, the map
$\Delta:\ARI\rightarrow \ARI$
trivially preserves alternality and push-invariance, as does its
inverse $\Delta^{-1}$, so we have
$$\Delta^{-1}(\ARI^{pol}_{al+push})=\ARI^\Delta_{al+push}.$$
The map $\Delta^{-1}$ pulls the $Dari$-bracket back to the $ari$-bracket,
so since $\ARI^{pol}_{al+push}$ is closed under the $Dari$-bracket,
$\ARI^\Delta_{al+push}$ is closed under the $ari$-bracket.\hfill{$\square$}
\vskip .3cm
\noindent {\bf Proposition 5.4.2.} {\it The space $\ARI_{al+push*circneut}^\Delta$ 
of alternal push-invariant moulds in $\ARI^\Delta$ whose swap is circ-neutral 
up to addition of a constant mould is closed under the $ari$-bracket.}
\vskip .3cm
\noindent {\bf Proof.} Proposition 2.6.1 showed that $\overline{\ARI}_{circneut}$ is closed under the $ari$-bracket, and it was shown at the end of the
proof of Theorem 4.7.1 that $\overline{\ARI}_{*circneut}$ is also closed
under the $ari$-bracket.  Let $A,B\in \ARI^{pol}_{al+push*circneut}$.
Then by Proposition 5.4.1, $ari(A,B)\in \ARI^\Delta_{al+push}$. 
But $swap(A)$ and $swap(B)$ lie in $\overline{\ARI}_{*circneut}$, which
is closed under the $ari$-bracket of moulds in $\overline{\ARI}$,
so $ari\bigl(swap(A),swap(B)\bigr)\in \overline{\ARI}_{*circneut}$. 
Since $A$ and $B$ are push-invariant, by (2.5.9) we have
$$swap\cdot ari(A,B)=ari\bigl(swap(A),swap(B)\bigr)$$
so $swap\cdot ari(A,B)\in \overline{\ARI}_{*circneut}$, and therefore
$ari(A,B)\in \ARI^\Delta_{al+push*circneut}$ as desired.\hfill{$\square$}
\vskip .5cm
\vskip .5cm
\noindent {\bf \S 5.5. The real function of the moulds $pal$ and $invpal$}
\vskip .3cm
Let $G\ARI$ denotes the set of all moulds with constant term 1, which
can be equipped with the multiplication law $gari$ (resp.~$Dgari$) 
corresponding to the Campbell-Hausdorff law on $\ARI_{ari}$ 
(resp.~on $\ARI_{Dari}$). 
We have exponential maps 
$$exp_{ari}:\ARI_{ari}\rightarrow G\ARI_{gari},\ \ \ \ exp_{Dari}:\ARI_{Dari}\rightarrow G\ARI_{Dgari}$$ 
(with inverses $log_{ari}$ and $log_{Dari}$); the map $exp_{ari}$
was defined in (2.7.1), and $exp_{Dari}$ is given by
$$exp_{Dari}(A)=1+\sum_{n\ge 1}\,Darit(A)^{n-1}(A).\eqno(5.5.1)$$
There is a unique group isomorphism $\Delta^*$ making the diagram
$$\xymatrix{
G\ARI_{gari}\ar[r]^{\Delta^*}&G\ARI_{Dgari}\\
\ARI_{ari}\ar[u]^{exp_{ari}}\ar[r]^{\Delta}&\ARI_{Dari}\ar[u]_{exp_{Dari}}.}\eqno(5.5.2)$$
commute. For all $G\in G\ARI_{Dgari}$, 
define the automorphism $Dgarit(G)$ of $G\ARI_{Dgari}$ by
$$Dgarit(G):=exp\bigl(Darit(A)\bigr)\eqno(5.5.3)$$
where $A=log_{Dari}(G)\in \ARI_{Dari}$. 
\vskip .3cm
Let $Ber_b$ denote the Bernoulli function defined by
$$Ber_b=ad(b)/\bigl(exp(ad(b))-1\bigr)=\sum_{r\ge 0} {{B_r}\over{r!}}ad(b)^r,\eqno(5.5.4)$$
and set
$$t_{01}=Ber_b(-a),\ \ t_{02}=Ber_{-b}(a),\ \ t_{12}=[a,b].\eqno(5.5.5)$$
Then 
$$t_{01}+t_{02}+t_{12}=0.$$
We write 
$$T=ma\bigl({\rm Lie}[t_{01},t_{12}]\bigr)\subset ma\bigl({\rm Lie}[a,b]\bigr)\subset \ARI_{lu}.\eqno(5.5.6)$$ 

We now give the result that really explains the important role
of the mould $pal$ throughout \'Ecalle's theory. Ideally, one would
like to have an automorphism of ${\rm Lie}[a,b]$ mapping 
$a\mapsto t_{02}$ and fixing $[a,b]$. But no Lie series exists
such that mapping $a\mapsto t_{02}$ and $b$ to that Lie series would
fix $[a,b]$. However, extended to moulds, i.e.~accepting that the image
of $b$ is a mould with denominators, there is such an isomorphism,
as stated in the next theorem.
\vskip .3cm
\noindent {\bf Theorem 5.5.1.} {\it Let $\Delta^*$ be the map
in diagram (5.5.2). Then
$$\Delta^*(invpal)=1-a+ma(t_{02}).$$}\par
\noindent Before proving this theorem, we give several preliminary 
results.
\vskip .3cm
\noindent {\bf Proposition 5.5.2.} {\it Let $G\in G\ARI_{Dgari}$. Then
$$Dgarit(G)\cdot a=a-1+G,\ \ Dgarit(G)\cdot ma([a,b])=ma([a,b]).\eqno(5.5.7)$$}\par
\noindent {\bf Proof.} Let $A=log_{Dari}(G)$, so that
$$Dgarit(G)=exp\bigl(Darit(A)\bigr)
=id + \sum_{n\ge 1}{{1}\over{n!}} Darit(A)^n.\eqno(5.5.8)$$
Applying (5.5.3) to the mould $B_1=ma([a,b])$, we see that $Dgarit(G)$
fixes $B_1$ since $Darit(A)$ annihilates $B_1$ (cf.~Prop. 5.3.2 (iv)).
Applying (5.5.3) to $a$, we find 
$$\eqalign{Dgarit(G)\cdot 
&=a+Darit(A)\cdot a+{{1}\over{2}}\,Darit(A)^2\cdot a+\cdots\cr
&=a+A+{{1}\over{2}}\,Darit(A)\cdot A+\cdots\cr
&=a-1+exp_{Dari}(A)\ \ \ {\rm by\ (5.5.1)}\cr
&=a-1+G.  }$$
This concludes the proof.\hfill{$\square$}
\vskip .1cm
\vskip .4cm
\noindent {\bf Lemma 5.5.3.} {\it Let $A\in \ARI$. The derivation
$-arit(A)+ad(A)$ extends from $\ARI_{lu}$ to $\ARI^a_{lu}$ taking the
value $lu(P,a)$ on $a$, and 
the automorphism ${\cal A}=exp\bigl(-arit(P)+ad(P)\bigr)$ satisfies
$${\cal A}\cdot a=R^{-1}aR$$
where $R=exp_{ari}(-A)$.}
\vskip .2cm
\noindent Proof. We know that the derivation $arit(A)$ extends to 
$a$ taking the value $0$ by Proposition 5.3.1 (iii), it suffices to 
check that $ad(A)$ extends to $a$ via $ad(A)\cdot a = lu(A,a)$, 
i.e., that this action respects (5.3.2).
Indeed, for all $P\in \ARI$ we have
$$ad(A)\cdot lu(P,a)=lu\bigl(ad(A)\cdot P,a\bigr)+lu\bigl(P,ad(A)\cdot a
\bigr)= lu\bigl(lu(A,P),a\bigr)+lu\bigl(P,lu(A,a)\bigr)$$
$$= lu\bigl(A,lu(P,a)\bigr)=ad(A)\cdot dur(P).$$
This proves the first statement.
Now, for a real parameter $t\in [0,1]$, let $R_t=exp_{ari}(-tA)$,
and let ${\cal A}_t$ denote the automorphism of $(\ARI^a_{lu})^{pol}$ defined by
$${\cal A}_t(a)=R_t^{-1}aR_t, \ \ \ {\cal A}_t(B_1)=B_1,$$
so that in particular ${\cal A}_1(a)=R^{-1}aR$. Let $D=log({\cal A}_1)$; we only need to prove that 
$D=-arit(A)+ad(A)$ on $(\ARI^a_{lu})^{pol}$.  We compute $D(a)$
and $D(b)$ by the linearization formula
$$D(a)={{d}\over{dt}}|_{t=0}\bigl({\cal A}_t(a)\bigr)\ \ \ {\rm and}\ \ \
D(b)={{d}\over{dt}}|_{t=0}\bigl({\cal A}_t(b)\bigr).$$
The second equality yields $D(b)=0$.  Let us compute $D(a)$.  Using
$R_0=1$ and ${{d}\over{dt}}|_{t=0}R_t=-A$, we find
$$\eqalign{D(a)&={{d}\over{dt}}|_{t=0}\Bigl({\cal A}_t(a)\Bigr)\cr
&={{d}\over{dt}}|_{t=0}\Bigl(R_t^{-1}aR_t\Bigr)\cr
&=\Bigl(-R_t^{-1}{{d}\over{dt}}(R_t)R_t^{-1}aR_t+R_t^{-1}a{{d}\over{dt}} (R_t)
\Bigr)|_{t=0}\cr
&=Aa-aA.}$$
Thus $D(a)=lu(A,a)=\bigl(-arit(A)+ad(A)\bigr)\cdot a$ and $D(b)=0=
\bigl(-arit(A)+ad(A)\bigr)\cdot b$, so $D=-arit(A)+ad(A)$, 
which concludes the proof of the lemma.\hfill{$\square$}
\vskip .4cm
Let $G\in G\ARI$, and recall the definition of the $mu$-dilator
$duG$ given in (4.2.2). The equivalent formula (4.2.3) 
$$dur\cdot G=mu(G,duG)$$
By (5.3.2) we have $dur\cdot G=lu(G,a)=mu(G,a)-mu(a,B)$, 
so this equality can be expressed as
$$mu(G,a)-mu(a,B)=mu(G,duG)$$
this means that $[G,a]=Ga-aG=G\,duG$, 
which multiplying by $G^{-1}$, gives us the useful formulation
$$G^{-1}aG=a-duG.\eqno(5.5.9)$$
\vskip .4cm
\noindent {\bf Proposition 5.5.4.} {\it The isomorphism
$$\Delta^*:G\ARI_{gari}\rightarrow G\ARI_{Dgari}$$
in diagram (5.5.2) is explicitly given by the formula
$$\Delta^*(G)=1-dar\bigl(du\,inv_{gari}(G)\bigr).\eqno(5.5.10)$$\par}
\noindent Proof.  Let $G\in G\ARI$, and set $A=log_{ari}(G)$ and
$R=exp_{ari}(-A)$.  By (5.2.3), we have
$$exp\Bigl(Darit\bigl(\Delta(A)\bigr)\Bigr)=
dar\circ exp\bigl(-arit(A)+ad(A)\bigr)\circ dar^{-1}.\eqno(5.5.11)$$
We have $dar(a)=a$ by Lemma 5.3.3 (i), and $dar$ is
an automorphism of $\ARI^a_{lu}$; in particular $du$ commutes with $dar$.  
Thus we have
$$\eqalign{exp\Bigl(Darit\bigl(\Delta(A)\bigr)\Bigr)\cdot a&=dar\circ 
exp\bigl(-arit(A)+ad(A)\bigr)\cdot a\cr
&=dar(R^{-1}\,a\,R)\ \ \ {\rm by\ Lemma\ 5.4.3}\cr
&=dar(R)^{-1}\,a\,dar(R)\cr
&=a-du\bigl(dar(R)\bigr)\ \ {\rm by\ (5.5.9)}\cr
&=a-dar\bigl(duR\bigr).}\eqno(5.5.12)$$
Now, using $A=log_{ari}(G)$, we compute 
$$\eqalign{\Delta^*(G)&=1-a+Dgarit\bigl(\Delta^*(G)\bigr)\cdot a
\ \ {\rm by\ (5.5.7)}\cr
&=1-a+Dgarit\Bigl(exp_{Dari}\bigl(\Delta(log_{ari}(G))\bigr)\Bigr)\cdot a
\ \ {\rm by\ diagram\ (5.5.2)}\cr
&=1-a+Dgarit\Bigl(exp_{Dari}\bigl(\Delta(A)\bigr)\Bigr)\cdot a\cr
&=1-a+exp\Bigl(Darit\bigl(\Delta(A)\bigr)\Bigr)\cdot a
\ \ {\rm by\ (5.5.3)}\cr
&=1-dar\bigl(du\,exp_{ari}(-A)\bigr)
\ {\rm by\ (5.5.12)}\cr
&=1-dar\bigl(du\,inv_{gari}(G)\bigr).}\eqno(5.5.13)$$
This proves the proposition.  \hfill{$\square$}
\vskip .3cm
\noindent {\bf Proof of Theorem 5.5.1.} The proof of the theorem
follows easily from the preliminary results together with the
definition of $dupal$ given in (4.2.4). Indeed, from (4.2.4) we
see immediately that for $r\ge 1$ we have
$$dar\cdot dupal(u_1,\ldots,u_r)={{B_r}\over{r!}}ma\bigl(ad(b)^r(-a)\bigr),$$
therefore $dar\cdot dupal$ agrees with $ma(t_{02})$ for $r\ge 2$,
but needs a sign correction for $r=1$, and is equal to $0$ for
$r=0$: more precisely we have
$$dar\cdot dupal = ma(t_{01}+t_{12})+a.\eqno(5.5.14)$$
Now, since $invpal=inv_{gari}(pal)$, we have $du\,inv_{gari}(invpal)=
dupal$ and so by (5.5.14) together with (5.5.10) applied to $G=invpal$,
we have $\Delta^*(invpal)(\emptyset)=1$ and
$$\eqalign{\Delta^*(invpal)(u_1,\ldots,u_r)
&=1-dar\cdot dupal(u_1,\ldots,u_r)\cr
&=1-ma(t_{01}+t_{12})-a\cr
&=1-a+ma(t_{02}),}\eqno(5.5.15)$$
which completes the proof of Theorem 5.5.1.\hfill{$\square$}

\vskip .5cm
\noindent {\bf \S 5.6. The real meaning of the operator $\Delta\circ 
Ad_{ari}(invpal)$}
\vskip .3cm
We complete this section with a final theorem concerning the nature 
of the operator $\Delta\circ Ad_{ari}(invpal)$ acting on a double 
shuffle element.
Recall that a polynomial mould is homogeneous of degree $n$ if
$F(u_1,\ldots,u_r)$ is a homogeneous polynomial of degree $n-r$
for all $r\ge 1$ (in particular $F$ is zero in depths greater than $n$).
The material in the following theorem all comes from the original
source [S2].
\vskip .3cm
\noindent {\bf Theorem 5.6.1.} {\it Let $F\in 
\ARI^{pol}_{\underline{al}*\underline{il}}=ma(\ds)$ be a double shuffle
mould of homogeneous degree $n$.  Set 
$$A=Ad_{ari}(invpal)(F),\ \ \ C=\Delta(A).\eqno(5.6.1)$$  
Set $B=ma(b)$ and $B_1=ma([a,b])$. Then 
\vskip .1cm\noindent
(i) $C$ is an alternal, polynomial, push-invariant mould.  Let $C'$
denote its partner (which is alternal and polynomial by Theorem 5.3.4 
(iii)). Thus the derivation $Darit(C)$ of $\ARI^a_{lu}$ restricts to a
derivation of $(\ARI^a_{lu})^{pol}_{al}=ma({\rm Lie}[a,b])$ such that
$$Darit(C)\cdot a=C,\ \ Darit(C)\cdot B=C',\ \ Darit(C)\cdot B_1=0.\eqno(5.6.2)$$
Let $c\in {\rm Lie}[a,b]$ be the Lie series such that $C=ma(c)$, so that
$c$ is also push-invariant, and let $D$ be the derivation of 
${\rm Lie}[a,b]$ defined by $D(a)=c$, $D([a,b])=0$. Then $Darit(C)$ 
is the mould version of $D$, i.e.~we have
$$ma\bigl(D(p)\bigr)=Darit(C)\cdot ma(p)\eqno(5.6.3)$$
for all $p\in {\rm Lie}[a,b]$.
\vskip .1cm\noindent
(ii) The derivation $Darit(C)$ restricts to a derivation of the Lie 
subalgebra $T$ of (5.5.6), given by
$$Darit(C)\cdot ma(t_{02})=ma([f(t_{02},-t_{12}),t_{02}]),\ \ 
Darit(C)\cdot ma(t_{12})=0.\eqno(5.6.4)$$
Equivalently, we have
$$D(t_{02})=f(t_{02},-t_{12}),t_{02}],\ \ D(t_{12})=0.\eqno(5.6.5)$$ 

\vskip .1cm\noindent
(iii) $Darit(C)$ is the unique derivation of $ma\bigl({\rm Lie}[a,b]\bigr)
=(\ARI^a_{lu})^{pol}_{al}$ which extends the derivation action on $T$ 
given in (5.6.4). Equivalently, $D$ is the unique derivation of
${\rm Lie}[a,b]$ extending the derivation on ${\rm Lie}[t_{02},
t_{12}]$ given in (5.6.5).
\vskip .2cm\noindent
(iv) For all $r\ge 1$, we have
$$\cases{C(u_1,\ldots,u_r)=0&if $r\not\equiv n$ mod 2\cr
C(u_1,\ldots,u_r)\ \hbox{is of degree }n+1&if $r\equiv n$ mod 2.}\eqno(5.6.6)$$
Equivalently, $D(a)$ has only terms of odd degree in $a,b$.}
\vskip .2cm
\noindent {\bf Proof.} (i) By Theorem 4.6.1, 
$Ad_{ari}(invpal)$ maps $\ARI_{\underline{al}*\underline{il}}$ to 
$\ARI_{\underline{al}*\underline{al}}$. By Lemma 2.5.5, $A$ is 
push-invariant. Thus $C$ is also push-invariant, since $\Delta$ 
respects push-invariance.  The fact that $C=\Delta(A)$ is a polynomial 
mould is shown in Theorem 5.1.1.  The fact that $Darit(C)$ acts as 
in (5.6.2) on $a$ and $B_1$ follows from Proposition 5.3.1 (iv).
To see that $Darit(C)\cdot B=C'$ follows from Lemma 5.3.3. To
show (5.6.3), it is enough to show that (5.6.3) holds for $a$
and $[a,b]$, with the value on $[a,b]$ being equal to $0$, since
in this case the value on $a$ determines the derivation. We
do have $D([a,b])=0=Darit(C)\cdot B_1$, and we also have
$$ma\bigl(D(a)\bigr)=ma(c)=C=Darit(C)\cdot a,$$
so the two derivations agree on all of $ma\bigl({\rm Lie}[a,b]\bigr)$,
completing the proof of (i).
\vskip .3cm
\noindent (ii)
By Lemma 5.3.3, $Darit(C)\cdot ma(t_{12})=0$. Let us compute
the action of $Darit(C)$ on $ma(t_{02})$.
By the nature of $Ad_{ari}(invpal)$ as an adjoint operator, we have 
$$\eqalign{Da&rit(C)=Darit\bigl(\Delta\circ Ad_{ari}(invpal)(F)\bigr)\cr
&=Dgarit\bigl(\Delta^*(invpal)\bigr)\circ Darit\bigl(\Delta(F)\bigr)\circ
Dgarit\bigl(\Delta^*(invpal)\bigr)^{-1}.}\eqno(5.6.7)$$
By Theorem 5.5.1, we have
$$\Delta^*(invpal)=1-a+ma(t_{02})$$
and by Proposition 5.5.2, for all moulds $G\in G\ARI$ we have 
$$Dgarit(G)\cdot a=a-1+G,$$
so taking $G=\Delta^*(invpal)$, we see that
$$Dgarit\bigl(\Delta^*(invpal)\bigr)\cdot a=ma(t_{02}).\eqno(5.6.8)$$
We use this to compute $Darit(C)\cdot ma(t_{02})$ using the
RHS of (5.6.7). By (5.6.8), the right-most operator
$Dgarit\bigl(\Delta^*(invpal)\bigr)^{-1}$ of (5.6.7) 
maps $ma(t_{02})\mapsto a$.  Next we compute the effect of the
middle operator of (5.6.7), $Darit\bigl(\Delta(F)\bigr)$,
on $a$. For this we recall that $\Delta=dar\circ dur$ (see (5.3.1)).
Let $f\in\ds$ be such that $F=ma(f)$. Then the effect of the 
$dar$-operator is expressed on $f$ by
$$dar(F)=ma\bigl(f(a,[b,a])\bigr)\eqno(5.6.9)$$
and we already saw that
$$dur(F)=lu(F,a)=ma([f,a]).\eqno(5.6.10)$$
Therefore since $Darit(P)\cdot a=P$ for all $P\in \ARI$ by (5.3.8), we have
$$Darit\bigl(\Delta(F)\bigr)\cdot a=\Delta(F)=dur\bigl(dar(F)\bigr)=
ma\bigl([f(a,[b,a]),a]\bigr).\eqno(5.6.11)$$
Finally, recalling that by Proposition 5.5.2 we have
$$Dgarit\bigl(\Delta^*(invpal)\bigr)\cdot ma([a,b])=ma([a,b])=ma(t_{12}),$$
we can apply the leftmost operator 
$Dgarit\bigl(\Delta^*(invpal)\bigr)$ of the RHS of (5.6.7) to (5.6.11) 
to obtain 
$$Dgarit\bigl(\Delta^*(invpal)\bigr)\cdot ma\bigl(f(a,[b,a]),a]\bigr)
=ma\bigl([f(t_{02},-t_{12}),t_{02}]\bigr).\eqno(5.6.12)$$
Thus altogether we have
$$Darit(C)\cdot ma(t_{02})
=ma\bigl([f(t_{02},-t_{12}),t_{02}]\bigr),$$
proving (5.6.4). The equality (5.6.5) follows immediately
from the agreement of $D$ and $Darit(C)$ proved in (i).
\vskip .3cm
\noindent (iii) We now show that there is a unique extension of
the derivation of (5.6.1) to all of $ma\bigl({\rm Lie}[a,b]\bigr)$.
We don't use mould theory for this part, so we can consider the
derivation $D$ defined by
$$D(t_{02})=[f(t_{02},-t_{12}),t_{02}],\ \ D(t_{12})=0\eqno(5.6.13)$$
and show that it has a unique extension to all of ${\rm Lie}[a,b]$
(cf. Lemma 2.1.2 of [S2]). In fact, knowing $D(t_{02})$, together
with the fact that $D(b)$ is necessarily the partner of $D(a)$
(as in (5.3.6)) because $D([a,b])=0$, allows us to recover 
$D(a)$ recursively, proceeding weight by weight. The minimal
weight term of $D(t_{0,2})$ is $[f^d(a,[b,a]),a]$, where
$d$ denotes the depth of $f$ and $f^d$ the minimal-depth part of $f$.
So the minimal weight of $D(t_{0,2})$ is equal to $n+d+1$, and
since $a$ is the lowest-weight part of $t_{02}$, this term
comes from $D(a)$.

Let $w=n+d+1$ denote the minimal weight. For all $m\ge w$, 
Let $t=D(t_{02})$, and for all $m\ge w$, let $t_m$ denote the
weight $m$ part of $t$, i.e.~$t=\sum_{m\ge w} t_m$. 
The recursive procedure to compute $D(a)$ runs as follows.  We 
first write out 
$$\eqalign{t&=D\bigl(Ber_{-b}(a)\bigr)\cr
&=D\bigl(a+{{1}\over{2}}[b,a]+
{{1}\over{12}}[b,[b,a]]-{{1}\over{720}}[b,[b,[b,[b,a]]]]+\cdots\bigr)\cr
&=D(a)+{{1}\over{2}}[D(b),a]+{{1}\over{2}}[b,D(a)]+{{1}\over{12}}[D(b),[b,a]]-{{1}\over{720}}[D(b),[b,[b,[b,a]]]\cr
&\qquad\qquad\qquad\qquad
-{{1}\over{720}}[b,[D(b),[b,[b,a]]]]
-{{1}\over{720}}[b,[b,[D(b),[b,a]]]]+
\cdots.}\eqno(5.6.14)$$
We construct $D(a)$ by solving (5.6.14) in successive weights starting 
with $w=n+d+1$.
We start by setting $D(a)_w=t_w$ since $D(a)$ is the only term in (2.1.9)
which can contribute to the lowest weight part $t_w$.
Let $D(b)_w$ be the partner of $D(a)_w$ as in (5.3.6), so as to ensure
that $D_w$ annihilates $[a,b]$.  We then continue to solve the
successive weight parts of (5.6.14) for $D(a)$ in terms of $t$ and the
previously determined lower weight
parts of $D(a)$ and $D(b)$.  For instance
the next few steps after weight $w$ are given by
$$\eqalign{D(a)_{w+1}&=t_{w+1}-{{1}\over{2}}[D(b)_w,a]-{{1}\over{2}}[b,D(a)_w],\cr
D(a)_{w+2}&=t_{w+2}-{{1}\over{2}}[D(b)_{w+1},a]-{{1}\over{2}}[b,D(a)_{w+1}]-{{1}\over{12}}[D(b)_w,[b,a]]\cr
&\qquad\qquad -{{1}\over{12}}[b,[D(b)_w,a]]-{{1}\over{2}}[b,[b,D(a)_w]],\cr
D(a)_{w+3}&=t_{w+3}-{{1}\over{2}}[D(b)_{w+2},a]-{{1}\over{2}}[b,D(a)_{w+2}]
 -{{1}\over{12}}[D(b)_{w+1},[b,a]]\cr
&\qquad\qquad -{{1}\over{12}}[b,[D(b)_{w+1},a]]
 -{{1}\over{12}}[b,[b,D(a)_{w+1}]]\cr
D(a)_{w+4}&=t_{w+4}-{{1}\over{2}}[D(b)_{w+3},a]-{{1}\over{2}}[b,D(a)_{w+3}]
-{{1}\over{12}}[D(b)_{w+2},[b,a]]\cr
&\qquad\qquad -{{1}\over{12}}[b,[D(b)_{w+2},a]]
-{{1}\over{12}}[b,[b,D(a)_{w+2}]]
+{{1}\over{720}} [D(b)_w,[b,[b,[b,a]]]\cr
&\qquad\qquad+{{1}\over{720}}[b,[D(b)_w,[b,[b,a]]]] +{{1}\over{720}}[b,[b,[D(b)_w,[b,a]]]]\cr
&\qquad\qquad+{{1}\over{720}}[b,[b,[b,[D(b)_w,a]]]] +{{1}\over{720}}[b,[b,[b,[b,D(a)_w]]]]\ldots}\eqno(5.6.15)$$
In this way we construct the unique Lie series $D(a)$ and its partner
$D(b)$ such that the derivation $D$ of ${\rm Lie}[a,b]$ extends
the derivation $D$ on ${\rm Lie}[t_{02},t_{12}]$ given in (5.6.13).
This construction shows that the derivation $D$ of ${\rm Lie}[a,b]$
extending (5.6.13) is unique, and since $Darit(C)$ does exactly
this, $Darit(C)$ must be the mould version of $D$, satisfying
$$Darit(C)\cdot a=ma\bigl(D(a)\bigr),
\ \ Darit(C)\cdot B=ma\bigl(D(b)\bigr).$$

\noindent (iv) We start by showing that 
$C(u_1,\ldots,u_r)$ is a polynomial of
degree $n+1$ in every depth $r\ge 1$. Observe that
the Lie series $t=D(t_{02})=[f(t_{02},-t_{12}),t_{02}]$ has constant $a$-degree
equal to $n+1$, since $f$ is assumed to be homogeneous of degree $n$
and $t_{02}$ and $t_{12}$ are both of degree $1$ in $a$. From
the weight-by-weight computation in (5.6.15), we note that
in every weight $m$, $D(a)_m$ is a Lie polynomial of constant
$a$-degree $n+1$ at every step, since the $a$-degree of the partner
$D(b)_m$ is one less than that of $D(a)_m$ at every weight $m$.
The part of the Lie series $D(a)$ of depth (=$b$-degree) $r$ corresponds
to the depth $r$ part of the mould $ma\bigl(D(a)\bigr)=C$, i.e.~to
$C(u_1,\ldots,u_r)$, and the $a$-degree corresponds to the degree
of the polynomial $C(u_1,\ldots,u_r)$, which is thus always equal to
$n+1$.  

Since $A\in (\ARI_{\underline{al}*\underline{al}}$, we know from 
Lemma 2.5.5 that $A$ is neg-invariant, i.e.
$$A(-u_1,\ldots,-u_r)=A(u_1,\ldots,u_r).$$
But we have $C=\Delta(A)$, i.e.
$$C(u_1,\ldots,u_r)=u_1\ldots,u_r(u_1+\cdots+u_r)\,A(u_1,\ldots,u_r),$$
so 
$$\eqalign{C(-u_1,\ldots,-u_r)&=(-1)^{r+1}u_1\ldots,u_r(u_1+\cdots+u_r)\,A(-u_1,\ldots,-u_r)\cr
&=(-1)^{r+1}u_1\cdots u_r(u_1+\cdots+u_r)A(u_1,\ldots,u_r)\cr
&=(-1)^{r+1}\,C(u_1,\ldots,u_r)\cr
&=(-1)^{n+1}\,C(u_1,\ldots,u_r),}$$
where the last equality holds because
$C(u_1,\ldots,u_r)$ is a polynomial of degree $n+1$. Therefore
if $r\not\equiv n$ mod 2, $C(u_1,\ldots,u_r)$ must be equal to zero.
In Lie algebra terms, the property (5.6.6) on $C=ma\bigl(D(a)\bigr)$
translates to $D(a)$ as saying that each term of $D(a)$ has 
$a$-degree $n+1$ and $b$-degree $\equiv n$ mod 2, which implies that
the total degree of every single term of $D(a)$ is odd.
This completes the proof of (iv), and thus of Theorem 5.6.1.\hfill{$\square$}
\vskip .3cm
\noindent {\bf Theorem 5.5.2.} {\it Let $f\in \ds$ and let the
moulds $A$ and $C$ and
the derivation $D$ be as in Theorem 5.6.1. Then the derivation $D$ of 
${\rm Lie}[a,b]$ acts on ${\rm Lie}[t_{01},t_{02}]$ by
$$\cases{D(t_{01})=[f(t_{01},-t_{12}),t_{01}]\bigr)\cr
D(t_{02})=[f(t_{02},-t_{12}),t_{02}]\bigr)\cr
D(t_{12})=0.}\eqno(5.6.16)$$}
\noindent Proof. For this, let $D$ be the
derivation on ${\rm Lie}[a,b]$ constructed in (iii) of Theorem 5.6.1
corresponding
to the mould derivation $Darit(C)$, and let $\iota$ denote the 
involutive automorphism of ${\rm Lie}[a,b]$ defined by 
$$\iota(a)=-a,\ \ \ \iota(b)=-b.$$
We claim that $D$ commutes with $\iota$ on ${\rm Lie}[a,b]$.
To check this, we consider the derivation $D'=\iota\circ D\circ \iota$
of ${\rm Lie}[a,b]$, and compare $D'$ with $D$ on $a$ and $[a,b]$. 
On $a$, we find that
$$D'(a)=(\iota\circ D\circ \iota)(a)=
\iota\circ D(-a)=-\iota\bigl(D(a)\bigr).$$
But since $D(a)$ has only odd-degree terms by (iv), we have
$\iota\bigl(D(a)\bigr)=-D(a)$ and therefore $D'(a)=D(a)$.
On $[a,b]$, since $\iota([a,b])=[a,b]$, we have
$$D'([a,b])=\iota\circ D([a,b])=0=D([a,b]).$$
Therefore $D$ and $D'$ agree on $a$ and $[a,b]$, and since a derivation
annihilating $[a,b]$ is uniquely determined by its value on $a$,
we have $D'=D$, proving that $D$ commutes with $\iota$. Now,
to prove (5.6.7), we simply observe that $t_{01}=\iota(t_{02})$, so
$$D(t_{01})=D\bigl(\iota(t_{02})\bigr)=\iota\bigl(D(t_{02})\bigr)
=\iota\bigl([f(t_{02},-t_{12}),t_{02}]\bigr)
=[f(t_{01},-t_{12}),t_{01}].$$
This concludes the proof.\hfill{$\square$}
\vskip .3cm

\vfill\eject
\centerline{\bf APPENDIX}
\vskip .7cm
\noindent {\bf \S A.1. Proof of Proposition 2.2.1.}
\vskip .5cm
Let $A\in \B\ARI$.  We prove that $amit(A)$ is a derivation for $mu$.  
The case for $anit(B)$ is analogous and we leave it as an exercise. It follows immediately 
from (2.2.3) and (2.2.4) that $axit(B)$ and $arit(B)$ are derivations.

For $amit$, we need to prove the identity
$$amit(A)\cdot mu(B,C)=mu\bigl(amit(A)\cdot B\,,\,C\bigr)+mu\bigl(B,amit(A)\cdot C\bigr).$$
Since $A,B,C$ all lie in B\ARI and therefore $0$-valued on the emptyset, we 
can remove $\b\ne\emptyset$ from the definition of $amit$; we have
$$\eqalign{amit(A)\cdot mu(B,C)&=\sum_{{{\w=\a\b\c}\atop{\c\ne\emptyset}}} mu(B,C)(\a\lceil\c)A(\b\rfloor)\cr
&=\sum_{{{\w=\a\b\c}\atop{\c\ne\emptyset}}} \sum_{\d_1\d_2=\a\lceil\c}B(\d_1)C(\d_2)A(\b\rfloor)\cr 
&=\sum_{{{\w=\a\b\c}\atop{\c\ne\emptyset}}} \sum_{\a_1\a_2=\a}B(\a_1)C(\a_2\lceil\c)A(\b\rfloor)
+\sum_{{{\w=\a\b\c}\atop{\c\ne\emptyset}}} \sum_{{{\c_1\c_2=\lceil\c}\atop{\c_1\ne\emptyset}}}B(\a\c_1)C(\c_2)A(\b\rfloor)\cr
&=\sum_{{{\w=\a_1\a_2\b\c}\atop{\c\ne\emptyset}}}B(\a_1)C(\a_2\lceil\c)A(\b\rfloor)
+\sum_{{{\w=\a\b\c_1\c_2}\atop{\c_1\ne\emptyset}}} B(\a\lceil\c_1)C(\c_2)A(\b\rfloor)\cr
&=\sum_{{{\w=\a_1\d}\atop{\d\ne\emptyset}}} B(\a_1)\sum_{{{\d=\a_2\b\c}\atop{\c\ne\emptyset}}}C(\a_2\lceil\c)A(\b\rfloor)
+\sum_{{{\w=\d\c_2}\atop{\d\ne\emptyset}}} \sum_{{{\d=\a\b\c_1}\atop{\c_1\ne\emptyset}}} B(\a\lceil\c_1)A(\b\rfloor)C(\c_2)\cr
&=\sum_{{{\w=\a_1\d}\atop{\d\ne\emptyset}}} B(\a_1)\bigl(amit(A)\cdot C)(\d)
+\sum_{{{\w=\d\c_2}\atop{\d\ne\emptyset}}} \bigl(amit(A)\cdot B\bigr)(\d)C(\c_2).  }$$
Noting that for $A,B,C\in \ARI$ we always have $(amit(A)\cdot B)(\emptyset)=(amit(A)\cdot C)(\emptyset)=0$,
we can drop the requirement $\d\ne\emptyset$ under the sum, and therefore obtain exactly
$$mu\bigl(B,amit(A)\cdot C)+mu\bigl(amit(A)\cdot B,C),$$
as desired.
\vskip .1cm
\noindent {\bf Exercise.} Show similarly that $anit$ is a derivation.
\vskip .8cm
\noindent {\bf \S A.2. Proofs of (2.4.7) and (2.4.8)}
\vskip .5cm
To prove these two key identities, we need
the following explicit expressions for the flexions occurring
in the definitions of the derivations, and the effect of $swap$:
$$\a\lceil\c=
\pmatrix{u_1&\cdots &u_k\cr v_1&\cdots &v_k}
\pmatrix{u_{k+1}+\cdots+u_{k+l+1}&\cdots& u_r\cr v_{k+l+1}&\cdots &v_r},$$
$$\b\rfloor=
\pmatrix{u_{k+1}&\cdots& u_{k+l}\cr v_{k+1}-v_{k+l+1}&\cdots &v_{k+l}-v_{k+l+1}}$$
$$\a\rceil\c=
\pmatrix{u_1&\cdots &u_{k-1}&u_k+\cdots+u_{k+l}\cr v_1&\cdots &v_{k-1}&v_k}
\pmatrix{u_{k+l+1}&\cdots& u_r\cr v_{k+l+1}&\cdots &v_r}.$$
$$\lfloor\b=
\pmatrix{u_{k+1}&\cdots& u_{k+l}\cr v_{k+1}-v_k&\cdots &v_{k+l}-v_k}.$$
Setting $SC=swap(C)$ for any mould $C$, we have
$$SC(\a\lceil\c)=SC
\pmatrix{u_1&\cdots &u_k &u_{k+1}+\cdots+u_{k+l+1}&u_{k+l+2}&\cdots& u_r\cr 
v_1&\cdots &v_k &v_{k+l+1}&v_{k+l+2}&\cdots &v_r}$$
$$=C\pmatrix{v_r&v_{r-1}-v_r&\cdots& v_{k+l+1}-v_{k+l+2}&v_k-v_{k+l+1}&
v_{k-1}-v_k&\cdots &v_1-v_2\cr u_1+\cdots+u_r&u_1+\cdots+u_{r-1}&
\cdots&u_1+\cdots+u_{k+l+1}&u_1+\cdots+u_k&u_1+\cdots+u_{k-1}&\cdots&u_1}$$
$$SC(\b\rfloor)=SC
\pmatrix{u_{k+1}&\cdots& u_{k+l}\cr v_{k+1}-v_{k+l+1}&\cdots &v_{k+l}-v_{k+l+1}}$$
$$=C\pmatrix{v_{k+l}-v_{k+l+1}&v_{k+l-1}-v_{k+l}&\cdots&v_{k+1}-v_{k+2}\cr
u_{k+1}+\cdots+u_{k+l}&u_{k+1}+\cdots+u_{k+l-1}&\cdots&u_{k+1} }$$ 
$$SC(\a\rceil\c)=SC
\pmatrix{u_1&\cdots& u_{k-1}&u_k+\cdots+u_{k+l} &u_{k+l+1}&\cdots& u_r\cr 
v_1&\cdots &v_{k-1}&v_k& v_{k+l+1}&\cdots &v_r}$$
$$=C\pmatrix{v_r&v_{r-1}-v_r&\cdots&v_{k+l+1}-v_{k+l+2}&v_k-v_{k+l+1}&v_{k-1}-v_k&\cdots&v_1-v_2\cr
u_1+\cdots+u_r&u_1+\cdots+u_{r-1}&\cdots&u_1+\cdots+u_{k+l+1} &u_1+\cdots+u_{k+l}&u_1+\cdots+
u_{k-1}&\cdots&u_1}$$
$$SC(\lfloor\b)=SC
\pmatrix{u_{k+1}&\cdots& u_{k+l}\cr v_{k+1}-v_k&\cdots &v_{k+l}-v_k}$$
$$=C\pmatrix{v_{k+l}-v_k& v_{k+l-1}-v_{k+l}&\cdots&v_{k+1}-v_{k+2}\cr
u_{k+1}+\cdots+u_{k+l}&u_{k+1}+\cdots+u_{k+l-1}&\cdots& u_{k+1}}$$
Applying the swap 
$$\pmatrix{u_1&u_2&\cdots&u_r\cr v_1&v_2&\cdots&v_r}\mapsto
\pmatrix{v_r&v_{r-1}-v_r&\cdots&v_1-v_2\cr u_1+\cdots+u_r&u_1+\cdots+u_{r-1}&
\cdots&u_1},$$ 
i.e.
$$\cases{u_1\mapsto v_r\cr
u_i\mapsto v_{r-i+1}-v_{r-i+2},\ {\rm if}\ i>1\cr
u_1+\cdots+u_i\mapsto v_{r-i+1}\cr
u_i+\cdots+u_j\mapsto -v_{r-i+2}+v_{r-j+1}\ \ {\rm if}\ i<j\cr
v_i\mapsto u_1+\cdots+u_{r-i+1}\cr
v_i-v_{i+1}\mapsto u_{r-i+1}\cr
v_i-v_j\mapsto u_{r-j+2}+\cdots+u_{r-i+1}\ \ {\rm if}\ i<j\cr
v_i-v_j\mapsto -u_{r-i+2}-\cdots-u_{r-j+1}\ \ {\rm if}\ i>j }$$
to these four terms yields
$$C\pmatrix{u_1&u_2&\cdots& u_{r-k-l}&u_{r-k-l+1}+\cdots+u_{r-k+1}&
u_{r-k+2}&\cdots &u_r\cr v_1&v_2&
\cdots&v_{r-k-l}&v_{r-k+1}&v_{r-k+2}&\cdots&v_r}$$
$$C\pmatrix{u_{r-k-l+1}&u_{r-k-l+2}&\cdots&u_{r-k}\cr
v_{r-k-l+1}-v_{r-k+1}&v_{r-k-l+2}-v_{r-k}&\cdots&v_{r-k}-v_{r-k+1}}$$ 
$$C\pmatrix{u_1&u_2&\cdots&u_{r-k-l}&u_{r-k-l+1}\cdots+u_{r-k+1}&
u_{r-k+2}&\cdots&u_r\cr
v_1&v_2&\cdots&v_{r-k-l}&v_{r-k-l+1}&v_{r-k+2}&\cdots&v_r}$$
$$C\pmatrix{-u_{r-k-l+2}-\cdots-u_{r-k+1}& u_{r-k-l+2}&\cdots&u_{r-k}\cr
v_{r-k-l+1}-v_{r-k+1}&v_{r-k-l+2}-v_{r-k+1}&\cdots& v_{r-k}-v_{r-k+1}}$$
Setting $m=r-k-l$, they can be written as
$$C\pmatrix{u_1&u_2&\cdots& u_{m}&u_{m+1}+\cdots+u_{r-k+1}&
u_{r-k+2}&\cdots &u_r\cr v_1&v_2&
\cdots&v_{m}&v_{r-k+1}&v_{r-k+2}&\cdots&v_r}$$
$$C\pmatrix{u_{m+1}&u_{m+2}&\cdots&u_{r-k}\cr
v_{m+1}-v_{r-k+1}&v_{m+2}-v_{r-k}&\cdots&v_{r-k}-v_{r-k+1}}$$ 
$$C\pmatrix{u_1&u_2&\cdots&u_{m}&u_{m+1}\cdots+u_{r-k+1}&
u_{r-k+2}&\cdots&u_r\cr
v_1&v_2&\cdots&v_{m}&v_{m+1}&v_{r-k+2}&\cdots&v_r}$$
$$C\pmatrix{-u_{m+2}-\cdots-u_{r-k+1}& u_{m+2}&\cdots&u_{r-k}\cr
v_{m+1}-v_{r-k+1}&v_{m+2}-v_{r-k+1}&\cdots& v_{r-k}-v_{r-k+1}}$$
Now putting $r-k=m+l$ gives
$$C\pmatrix{u_1&u_2&\cdots& u_{m}&u_{m+1}+\cdots+u_{m+l+1}&
u_{m+l+2}&\cdots &u_r\cr v_1&v_2&
\cdots&v_{m}&v_{m+l+1}&v_{m+l+2}&\cdots&v_r}$$
$$C\pmatrix{u_{m+1}&u_{m+2}&\cdots&u_{m+l}\cr
v_{m+1}-v_{m+l+1}&v_{m+2}-v_{m+l+1}&\cdots&v_{m+l}-v_{m+l+1}}$$ 
$$C\pmatrix{u_1&u_2&\cdots&u_{m}&u_{m+1}\cdots+u_{m+l+1}&
u_{m+l+2}&\cdots&u_r\cr
v_1&v_2&\cdots&v_{m}&v_{m+1}&v_{m+l+2}&\cdots&v_r}$$
$$C\pmatrix{-u_{m+2}-\cdots-u_{m+l+1}& u_{m+2}&\cdots&u_{m+l}\cr
v_{m+1}-v_{m+l+1}&v_{m+2}-v_{m+l+1}&\cdots& v_{m+l}-v_{m+l+1}}$$
\vskip .3cm
Using all these, we can now prove (2.4.7) and (2.4.8).
\vskip .3cm
\noindent {\bf Proof of (2.4.7).} We have
$$swap\Bigl(amit\bigl(swap(B)\bigr)\cdot swap(A)\Bigr)=
swap\Bigl(\sum_{{{\w=\a\b\c}\atop{\b,\c\ne\emptyset}}} SA(\a\lceil\c)SB(\b\rfloor)\Bigr)$$
$$=swap\Bigl[\sum_{l=1}^{r-1}\sum_{m=1}^{r-l}$$
$$A\pmatrix{v_r&v_{r-1}-v_r&\cdots& v_{k+l+1}-v_{k+l+2}&v_k-v_{k+l+1}&
v_{k-1}-v_k&\cdots &v_1-v_2\cr u_1+\cdots+u_r&u_1+\cdots+u_{r-1}&
\cdots&u_1+\cdots+u_{k+l+1}&u_1+\cdots+u_k&u_1+\cdots+u_{k-1}&\cdots&u_1}$$
$$\cdot B\pmatrix{v_{k+l}-v_{k+l+1}&v_{k+l-1}-v_{k+l}&\cdots&v_{k+1}-v_{k+2}\cr
u_{k+1}+\cdots+u_{k+l}&u_{k+1}+\cdots+u_{k+l-1}&\cdots&u_{k+1} }\Bigr]$$ 
$$=\sum_{l=1}^{r-1}\sum_{m=1}^{r-l}
A\pmatrix{u_1&u_2&\cdots& u_{m}&u_{m+1}+\cdots+u_{m+l+1}&
u_{m+l+2}&\cdots &u_r\cr v_1&v_2&
\cdots&v_{m}&v_{m+l+1}&v_{m+l+2}&\cdots&v_r}$$
$$\cdot B\pmatrix{u_{m+1}&u_{m+2}&\cdots&u_{m+l}\cr
v_{m+1}-v_{m+l+1}&v_{m+2}-v_{m+l+1}&\cdots&v_{m+l}-v_{m+l+1}}$$ 
$$=\sum_{l=1}^{r-1}\sum_{k=1}^{r-l}
A\pmatrix{u_1&u_2&\cdots& u_{k}&u_{k+1}+\cdots+u_{k+l+1}&
u_{k+l+2}&\cdots &u_r\cr v_1&v_2&
\cdots&v_{k}&v_{k+l+1}&v_{k+l+2}&\cdots&v_r}$$
$$\cdot B\pmatrix{u_{k+1}&u_{k+2}&\cdots&u_{k+l}\cr
v_{k+1}-v_{k+l+1}&v_{k+2}-v_{k+l+1}&\cdots&v_{k+l}-v_{k+l+1}}$$ 
$$=\sum_{l=1}^{r-1}\sum_{k=0}^{r-l-1}
A\pmatrix{u_1&u_2&\cdots& u_{k}&u_{k+1}+\cdots+u_{k+l+1}&
u_{k+l+2}&\cdots &u_r\cr v_1&v_2&
\cdots&v_{k}&v_{k+l+1}&v_{k+l+2}&\cdots&v_r}$$
$$\cdot B\pmatrix{u_{k+1}&u_{k+2}&\cdots&u_{k+l}\cr
v_{k+1}-v_{k+l+1}&v_{k+2}-v_{k+l+1}&\cdots&v_{k+l}-v_{k+l+1}}$$ 
$$-\sum_{l=1}^{r-1}
A\pmatrix{u_{1}+\cdots+u_{l+1}& u_{l+2}&\cdots &u_r\cr 
v_{l+1}&v_{l+2}&\cdots&v_r}
\cdot B\pmatrix{u_{1}&u_{2}&\cdots&u_{l}\cr
v_{1}-v_{l+1}&v_{2}-v_{l}&\cdots&v_{l}-v_{l+1}}$$ 
$$+\sum_{l=1}^{r-1}
A\pmatrix{u_1&u_2&\cdots& u_{r-l}\cr
v_1&v_2& \cdots&v_{r-l}}
\cdot B\pmatrix{u_{r-l+1}&u_{r-l+2}&\cdots&u_{r}\cr
v_{r-l+1}&v_{r-l+2}&\cdots&v_{r}}$$ 
$$=amit(B)\cdot A-swap\Bigl(mu\bigl(swap(A),swap(B)\bigr)\Bigr)+mu(A,B).$$

\vskip .5cm
\noindent {\bf Proof of (2.4.8).} We have
$$swap\Bigl(anit\bigl(swap(B)\bigr)\cdot swap(A)\Bigr)=
swap\Bigl(\sum_{{{\w=\a\b\c}\atop{\a,\b\ne\emptyset}}} SA(\a\rceil\c)SB(\lfloor\b)\Bigr)$$
$$=swap\Bigl[\sum_{l=1}^{r-1}\sum_{k=1}^{r-l}$$
$$A\pmatrix{v_r&v_{r-1}-v_r&\cdots&v_{k+l+1}-v_{k+l+2}&v_k-v_{k+l+1}&v_{k-1}-v_k&\cdots&v_1-v_2\cr
u_1+\cdots+u_r&u_1+\cdots+u_{r-1}&\cdots&u_1+\cdots+u_{k+l+1} &u_1+\cdots+u_{k+l}&u_1+\cdots+
u_{k-1}&\cdots&u_1}$$
$$\cdot B\pmatrix{v_{k+l}-v_k& v_{k+l-1}-v_{k+l}&\cdots&v_{k+1}-v_{k+2}\cr
u_{k+1}+\cdots+u_{k+l}&u_{k+1}+\cdots+u_{k+l-1}&\cdots& u_{k+1}} \Bigr]$$
$$=\sum_{l=1}^{r-1}\sum_{m=0}^{r-l-1}A\pmatrix{u_1&u_2&\cdots&u_{m}&u_{m+1}\cdots+u_{m+l+1}&
u_{m+l+2}&\cdots&u_r\cr
v_1&v_2&\cdots&v_{m}&v_{m+1}&v_{m+l+2}&\cdots&v_r}$$
$$\cdot B\pmatrix{-u_{m+2}-\cdots-u_{m+l+1}& u_{m+2}&\cdots&u_{m+l}\cr
v_{m+1}-v_{m+l+1}&v_{m+2}-v_{m+l+1}&\cdots& v_{m+l}-v_{m+l+1}}$$
$$=\sum_{l=1}^{r-1}\sum_{m=0}^{r-l-1}A\pmatrix{u_1&u_2&\cdots&u_{m}&u_{m+1}\cdots+u_{m+l+1}&
u_{m+l+2}&\cdots&u_r\cr
v_1&v_2&\cdots&v_{m}&v_{m+1}&v_{m+l+2}&\cdots&v_r}$$
$$\cdot push(B)\pmatrix{u_{m+2}&u_{m+3}&\cdots&u_{m+l+1}\cr
v_{m+2}-v_{m+1}&v_{m+3}-v_{m+1}&\cdots& v_{m+l+1}-v_{m+1}}$$
$$=\sum_{l=1}^{r-1}\sum_{m=1}^{r-l}A\pmatrix{u_1&u_2&\cdots&u_{m-1}&u_{m}\cdots+u_{m+l}&
u_{m+l+1}&\cdots&u_r\cr
v_1&v_2&\cdots&v_{m-1}&v_{m}&v_{m+l+1}&\cdots&v_r}$$
$$\cdot push(B)\pmatrix{u_{m+1}&u_{m+2}&\cdots&u_{m+l}\cr
v_{m+1}-v_{m}&v_{m+2}-v_{m}&\cdots& v_{m+l}-v_{m}}$$
$$=\sum_{l=1}^{r-1}\sum_{k=1}^{r-l}A\pmatrix{u_1&u_2&\cdots&u_{k-1}&u_{k}\cdots+u_{k+l}& u_{k+l+1}&\cdots&u_r\cr
v_1&v_2&\cdots&v_{k-1}&v_{k}&v_{k+l+1}&\cdots&v_r}$$
$$\cdot push(B)\pmatrix{u_{k+1}&u_{k+2}&\cdots&u_{k+l}\cr
v_{k+1}-v_{k}&v_{k+2}-v_{k}&\cdots& v_{k+l}-v_{k}}$$
\vskip .2cm
\noindent {\bf \S A.3. Proof of Lemma 3.2.1.}
\vskip .3cm
We first prove (3.2.8), then (3.2.7). By (3.2.5), we have 
$mi_f(v_1,\ldots,v_r)=\iota_Y(f_Y^r)$.  Since $mi$
is additive, we may assume that $f$ is a monomial,
$f=x^{a_0-1}y\cdots yx^{a_r-1}$.  Then 
$$\pi_Y(f)=\cases{f&if $a_0=1$\cr 0&otherwise.}$$
So
$$\ret_X(\pi_Y(f))=\cases{x^{a_r-1}y\cdots x^{a_1-1}y&if $a_0=1$\cr 0&otherwise.}$$
and
$$f_Y=\cases{y_{a_r}\cdots y_{a_1}&if $a_0=1$\cr 0&otherwise.}$$
Thus 
$$mi_f(v_1,\ldots, v_r)=\iota_Y(f_Y)=
\cases{v_1^{a_r-1}\cdots v_r^{a_1-1}&if $a_0=1$\cr 0&otherwise.}$$
Now by (3.2.4), we have
$$vimo_f(z_0,\ldots,z_r)=z_0^{a_0-1}z_1^{a_1-1}\cdots z_r^{a_r-1},$$
so as desired, we have
$$mi_f(v_1,\ldots, v_r)=
vimo_f(0,v_r,\ldots,v_1)=\cases{v_r^{a_1-1}\cdots v_1^{a_r-1}&if $a_0=1$\cr
0&otherwise.}$$
This settles the proof of (3.2.8) for $mi$.  

The case of $ma$ is a little more complicated.  Again,
by additivity, we can assume that
$f$ is a monomial $C_{a_1}\cdots C_{a_r}$ in the $C_i$.  We will prove it
by induction on $r$ (though there might be a better way).  For the base
case, $r=1$, we have $n=a_1$ and
$$f=C_{a_1}=\sum_{i=0}^{a_1-1} (-1)^{i}C_{a_1-1}^ix^{a_1-1-i}yx^{i},$$
$$vimo_f(z_0,z_1)=\sum_{i=0}^{a_1-1} (-1)^{i}C_{a_1-1}^iz_0^{a_1-1-i}z_1^{i},$$
$$vimo_f(0,u_1)=(-1)^{a_1-1}u_1^{a_1-1}=(-1)^{r+n}u_1^{a_1-1}=ma_f(u_1)$$
using Ecalle's definition, and comparing with (3.2.5), we also have
$$ma_f(u_1)=(-1)^{r+n}\iota_C(C_{a_1})=(-1)^{r+n}u_1^{a_1-1},$$
which settles the base case.

Now make the induction hypothesis that (3.2.7) holds
up to depth $r-1$, and let $f=C_{a_1}\cdots C_{a_{r-1}}C_{a_r}$.  Using
(3.2.5), we have
$$ma_f(u_1,\ldots,u_r)=(-1)^{r+n}\iota_C(f)=
(-1)^{r+n}u_1^{a_1-1}\cdots u_r^{a_r-1}.$$
Let us write $g=C_{a_1}\cdots C_{a_{r-1}}$.  Then again from (3.2.5), we have
$$ma_f(u_1,\ldots,u_r)=ma_g(u_1,\ldots,u_{r-1})ma_{C_{a_r}}(u_r).$$
By the induction hypothesis, we have 
$$\cases{ma_{C_{a_r}}(u_r)=vimo_{C_{a_r}}(0,u_r)=(-1)^{a_r-1}u_r^{a_r-1}\cr
ma_g(u_1,\ldots,u_{r-1})=
vimo_g(0,u_1,\ldots,u_1+\cdots+u_{r-1}).}$$
So to prove (3.2.7), we have to show that
$$\eqalign{vimo_f(0,u_1,\ldots,u_1+\cdots+u_r)&=vimo_g(0,u_1,\ldots,u_1+\cdots+u_{r-1})
\,vimo_{C_{a_r}}(0,u_r)\cr
&=(-1)^{a_r-1}vimo_g(0,u_1,\ldots,u_1+\cdots+u_{r-1})u_r^{a_r-1}.}\eqno(A.3.1)$$
Write
$$g=\sum_{{\bf a}=(a_0,\ldots,a_{r-1})} c_{\bf a}x^{a_0-1}y\cdots yx^{a_{r-1}-1}.$$
Then 
$$vimo_g(z_0,\ldots,z_{r-1})= \sum_{{\bf a}=(a_0,\ldots,a_{r-1})}
c_{\bf a}z_0^{a_0-1}z_1^{a_1-1}\cdots z_{r-1}^{a_{r-1}-1},$$
and
$$vimo_g(0,u_1,\ldots,u_1+\cdots+u_{r-1})=
\sum_{{\bf a}=(1,a_1,\ldots,a_{r-1})} c_{\bf a}u_1^{a_1-1}
(u_1+u_2)^{a_2-1}\cdots 
(u_1+\ldots+u_{r-1})^{a_{r-1}-1}.$$
Thus the second term in (A.3.1) is given by
$$vimo_g(0,u_1,\ldots,u_1+\cdots+u_{r-1})\,vimo_{C_{a_r}}(0,u_r)$$
$$=
(-1)^{a_r-1}\sum_{{\bf a}=(1,a_1,\ldots,a_{r-1})} c_{\bf a}u_1^{a_1-1}
(u_1+u_2)^{a_2-1}\cdots 
(u_1+\ldots+u_{r-1})^{a_{r-1}-1}u_r^{a_r-1}.\eqno(A.3.2)$$

But also 
$$f=gC_{a_r}=
\sum_{{\bf a}=(a_0,\ldots,a_{r-1})}\sum_{j=0}^{a_r-1} (-1)^{j}\Bigl({{a_r-1}\atop{j}}\Bigr)c_{\bf a}x^{a_0-1}y\cdots yx^{a_{r-1}-1}x^{a_r-1-j}yx^j,$$
so
$$vimo_f(z_0,\ldots,z_r)=
\sum_{{\bf a}=(a_0,\ldots,a_{r-1})}\sum_{j=0}^{a_r-1} (-1)^{j}\Bigl({{a_r-1}\atop{j}}\Bigr)c_{\bf a}z_0^{a_0-1}z_1^{a_1-1}\cdots z_{r-1}^{a_{r-1}-2+a_r-j}z_r^{j},$$
so
$$vimo_f(0,z_1,\ldots,z_r)=
\sum_{\a=(1,a_1,\ldots,a_r)}\sum_{j=0}^{a_r-1} (-1)^{j}\Bigl({{a_r-1}\atop{j}}\Bigr)c_{\bf a}z_1^{a_1-1}z_2^{a_2-1}\cdots z_{r-1}^{a_{r-1}-2+a_r-j}z_r^{j},$$
so finally the first term in (A.3.1) is given by
$$vimo_f(0,u_1,\ldots,u_1+\cdots+u_r)=$$
$$\sum_{\a=(1,a_1,\ldots,a_r)}\sum_{j=0}^{a_r-1} (-1)^{j}\Bigl({{a_r-1}\atop{j}}\Bigr)c_{\bf a}u_1^{a_1-1}(u_1+u_2)^{a_2-1}\cdots (u_1+\cdots+u_{r-1})^{a_{r-1}-2+a_r-j}(u_1+\cdots+u_r)^{j}$$
$$=\sum_{\a=(1,a_1,\ldots,a_{r-1})}
c_{\bf a}u_1^{a_1-1}(u_1+u_2)^{a_2-1}\cdots (u_1+\cdots+u_{r-1})^{a_{r-1}-1}
\cdot $$
$$\Biggl(\sum_{j=0}^{a_r-1} (-1)^{j}\Bigl({{a_r-1}\atop{j}}\Bigr)
(u_1+\cdots+u_{r-1})^{a_r-j}(u_1+\cdots+u_r)^{j}\Biggr)$$
$$=(-1)^{a_r-1}\sum_{\a=(1,a_1,\ldots,a_{r-1})}
c_{\bf a}u_1^{a_1-1}(u_1+u_2)^{a_2-1}\cdots (u_1+\cdots+u_{r-1})^{a_{r-1}-1}
\cdot u_r^{a_r-1}$$
since the factor between large parenthesis is just the binomial expansion of
$$\bigl((u_1+\ldots +u_{r-1})-(u_1+\cdots+u_r)\bigr)^{a_r-1}
=(-1)^{a_r-1}u_r^{a_r-1}.$$
But this is equal to the second term as given in (A.3.2), so (A.3.1) holds, thus
proving (3.2.7).
\vskip .8cm
\noindent {\bf \S A.4. Proof of Proposition 3.3.2}
\vskip .5cm
We need to show that
$$arit(A)(BC)=arit(A)(B)C+Barit(A)(C).\eqno(A.4.1)$$
Using the definition of $S_A(B)$ from (4.1),
$$\bigl(S_A(B)\bigr)(\w)=\sum_{\w=\a\b\c} B(\a\c')A(\b)
-\sum_{{{\w=\a\b\c}\atop{\a\ne \emptyset}}} B(\a''\c)A(\b),$$
and $arit(A)(B)=S_A(B)-BA$, we write
$$\bigl(arit(A)(B)\bigr)(\w)=\sum_{\w=\a\b\c} B(\a\c')A(\b)
-\sum_{{{\w=\a\b\c}\atop{\a\ne \emptyset}}} B(\a''\c)A(\b)
-\sum_{\w=\a\b} B(\a)A(\b).$$
Splitting the first sum over $\c=\emptyset$ and $\c\ne\emptyset$, and
recalling that $\c'=\emptyset$ when $\c=\emptyset$, this is equal to
$$\eqalign{\bigl(arit(A)(B)\bigr)(\w)&=\sum_{{{\w=\a\b\c}\atop{\c\ne\emptyset}}} 
B(\a\c')A(\b) +\sum_{\w=\a\b} B(\a)A(\b)
-\sum_{{{\w=\a\b\c}\atop{\a\ne \emptyset}}} B(\a''\c)A(\b)
-\sum_{\w=\a\b} B(\a)A(\b)\cr
&=\sum_{{{\w=\a\b\c}\atop{\c\ne\emptyset}}} 
B(\a\c')A(\b) -\sum_{{{\w=\a\b\c}\atop{\a\ne \emptyset}}} B(\a''\c)A(\b).}
\eqno(A.4.2) $$
Thus we can write the right-hand side of (A.4.1) as
$$\bigl(arit(A)(B)C+Barit(A)(C)\bigr)(\w)=
\sum_{\w=\u\v} \Bigl(\sum_{{{\u=\a\b\c}\atop{\c\ne\emptyset}}} 
B(\a\c')A(\b)C(\v) -\sum_{{{\u=\a\b\c}\atop{\a\ne \emptyset}}} B(\a''\c)A(\b)
C(\v)\Bigr)$$
$$\qquad\qquad\qquad\qquad \qquad\qquad\qquad\ \ 
+ \Bigl(\sum_{{{\v=\a\b\c}\atop{\c\ne\emptyset}}} B(\u)C(\a\c')A(\b) -
\sum_{{{\v=\a\b\c}\atop{\a\ne \emptyset}}} B(\u)C(\a''\c)A(\b)\Bigr),$$
or again as
$$\bigl(arit(A)(B)C+Barit(A)(C)\bigr)(\w)=
\sum_{{{\w=\a\b\c\v}\atop{\c\ne\emptyset}}} 
B(\a\c')A(\b)C(\v) 
-\sum_{{{\w=\a\b\c\v}\atop{\a\ne \emptyset}}} B(\a''\c)A(\b) C(\v)$$
$$\qquad\qquad\qquad\qquad \qquad\qquad\qquad\ \ 
+ \sum_{{{\w=\u\a\b\c}\atop{\c\ne\emptyset}}} 
B(\u)C(\a\c')A(\b) -
\sum_{{{\w=\u\a\b\c}\atop{\a\ne \emptyset}}} B(\u)C(\a''\c)A(\b).
\eqno(A.4.3)$$

By (A.4.2), the left-hand side of (A.4.1) can be written
$$arit(A)(BC)=
\sum_{{{\w=\a\b\c}\atop{\c\ne\emptyset}}} 
BC(\a\c')A(\b) -\sum_{{{\w=\a\b\c}\atop{\a\ne \emptyset}}} BC(\a''\c)A(\b)$$
$$=\sum_{{{\w=\a\b\c}\atop{\c\ne\emptyset}}} 
\sum_{\a\c'=\u\v} B(\u)C(\v)A(\b) -
\sum_{{{\w=\a\b\c}\atop{\a\ne \emptyset}}} \sum_{\a''\c=\u\v}
B(\u)C(\v)A(\b)\eqno(A.4.4)$$
$$=\sum_{{{\w=\a_1\a_2\b\c}\atop{\c\ne\emptyset}}} 
B(\a_1)C(\a_2\c')A(\b) 
+\sum_{{{\w=\a\b\c_1\c_2}\atop{\c_1\ne\emptyset}}} 
B(\a\c'_1)C(\c_2)A(\b) $$
$$\ \ \ \ \ \ \ \ \ \ \ \ \ \ \ -\sum_{{{\w=\a_1\a_2\b\c}\atop{\a_2\ne \emptyset}}} 
B(\a_1)C(\a''_2\c)A(\b)
-\sum_{{{\w=\a\b\c_1\c_2}\atop{\a\ne \emptyset}}} 
B(\a''\c_1)C(\c_2)A(\b).\eqno(A.4.5)$$
The passage from (A.4.4) to (A.4.5) is obtained by separating the first term
into two terms according to whether the decomposition 
$\a\c'=\u\v$ is of the form $\u=\a_1$, $\v=\a_2\c'$ or of the form
$\u=\a\c'_1$, $\v=\c_2$ with $\c_1\ne\emptyset$ (otherwise the case
$\u=\a$, $\v=\c'$ is counted twice).  The second term is separated into
two terms according to whether the decomposition $\a''\c=\u\v$ is of the form
$\u=\a''\c_1$, $v=\c_2$ or of the form $\u=\a_1$, $\v=\a''_2\c$ with
$\a_2\ne\emptyset$ (otherwise the term $\u=\a$, $v=\c$ is counted twice).

Relabeling the indices in the first term of (A.4.5) by $\a_1\mapsto \u$, 
$\a_2\mapsto \a$, we see that this term is equal to the third
term of (A.4.3). 

Relabeling the indices in the second term of (A.4.5) by $\c_1\mapsto \c$, 
$\c_2\mapsto \v$, we see that this term is equal to the first
term of (A.4.3). 

Relabeling the indices in the third term of (A.4.5) by $\a_1\mapsto \u$, 
$\a_2\mapsto \a$, we see that this term is equal to the fourth
term of (A.4.3). 

Relabeling the indices in the fourth term of (A.4.4) by $\c_1\mapsto \c$,
$\c_2\mapsto \v$, we see that this term is equal to the second term of (A.4.3).

So (A.4.3) is equal to (A.4.5), i.e. $arit(A)(B)(C)+Barit(A)(C)=arit(A)(BC)$, proving 
that $arit(A)$ is a derivation.\hfill{$\square$}
$$=anit\bigl(push(B)\bigr)\cdot A.$$
\noindent {\bf \S A.5. Proof of Lemma 3.4.1}
\vskip .5cm
 (i) Let $f\in \F_n$.  We show that $f$ satisfies shuffle if and only if 
$ma(f)$ is alternal.  We know that $f$ satisfies shuffle if and only if
$f\in {\rm Lie}[x,y]$, so $f$ satisfies shuffle if and only if
$$f\in \F_n\cap {\rm Lie}[x,y]={\rm Lie}[C_1,C_2,\ldots]$$
where $C_i=ad(x)^{i-1}(y)$.  Thus the shuffle relations on $f$ written
in $x,y$ are equivalent to the shuffle conditions written in the $C_i$.
I.e., assuming by additivity that $f$ is of homogeneous depth $r$, we can
write 
$$f=\sum_{{\bf a}=(a_1,\ldots,a_r)} c_{\bf a}C_{a_1}\cdots C_{a_r},\eqno(A.5.1)$$
and the shuffle relations are
$$\sum_{w\in sh\bigl((C_{a_1},\cdots,C_{a_i}),(C_{a_{i+1}},\ldots,C_{a_r})\bigr)} (f|w)=0.\eqno(A.5.2)$$
It is convenient to write the shuffle using the set $Sh(i,r)\subset S_r$
of permutations $\sigma$ of $\{1,\ldots,r\}$ satisfying
$$\sigma(1)<\cdots<\sigma(i)\ \ {\rm and}\ \ \sigma(i+1)<\cdots<\sigma(r).$$
Then (A.5.2) can be rewritten 
$$\sum_{\sigma\in Sh(i,r)} (f|C_{a_{\sigma^{-1}(1)}}\cdots C_{a_{\sigma^{-1}(r)}})=
\sum_{\sigma\in Sh(i,r)} c_{a_{\sigma^{-1}(1)},\ldots,a_{\sigma^{-1}(r)}}=0.\eqno(A.5.3)$$
Let us compare this property with the alternality condition on 
$$ma(f)=\sum_{{\bf a}} c_{\bf a}u_1^{a_1-1}\cdots u_r^{a_r-1}.$$
The alternality conditions are given by
$$\eqalign{0&=\sum_{w\in sh\bigl((u_1,\ldots,u_i),(u_{i+1},\ldots,u_r)\bigr)}
ma(f)(w)\cr
&=\sum_{\sigma\in Sh(i,r)} \sum_{\bf a} c_{a_1,\ldots,a_r} u_1^{a_{\sigma(1)}-1}\cdots u_r^{a_{\sigma(r)}-1}\cr
&=\sum_{\sigma\in Sh(i,r)}
\sum_{\bf a} c_{a_{\sigma^{-1}(1)},\ldots,a_{\sigma^{-1}(r)}} 
u_{\sigma^{-1}(1)}^{a_1-1}\cdots u_{\sigma^{-1}(r)}^{a_r-1},}$$
which monomial by monomial implies that
$$\sum_{\sigma\in Sh(i,r)} c_{a_{\sigma^{-1}(1)},\ldots,a_{\sigma^{-1}(r)}}=0,$$
which is identical to (A.5.3).

\vskip .4cm
(ii) The proof is identical to (i), with $u_i$ replaced by $v_i$ and
$C_{a_i}$ replaced by $y_{a_i}$.
\vskip .4cm
(iii) As in \S 2.3, we write $st(r,s)$ for the set of words in the stuffle sum
$$st\bigl((a_1,\ldots,a_r),(a_{r+1},\ldots,a_{r+s})\bigr).$$ 
We saw in \S 2.3 that each stuffle sum $st(r,s)$ corresponds to an alternility 
sum associated to a mould $A$, containing one term for each word in the 
stuffle set. Let $A_{r,s}$ denote the alternality sum associated to $A$
corresponding to $st(r,s)$ as in \S 2.3; recall for example that
$st(1,2)=(a,b,c)+(b,a,c)+(b,c,a)+(a+b,c)+(b,a+c)$ and
$$A_{1,2}(v_1,v_2,v_3)=A(v_1,v_2,v_3)+A(v_2,v_1,v_3)+A(v_2,v_3,v_1)+$$
$${{1}\over{(v_1-v_2)}}\bigl(A(v_1,v_3)-A(v_2,v_3)\bigr)
+{{1}\over{(v_1-v_3)}}\bigl(A(v_2,v_1)-A(v_2,v_3)\bigr).$$
Assume that $A$ is a polynomial-valued mould, i.e. $A=mi(f)=swap(ma(f))$ 
for a power series $f\in \Q$ with constant term 1.  We will show that 
$A$ is symmetril if and only if $f$ satisfies the stuffle relations in the
sense of (1.3.3).  To do this, we write
$$A_r(v_1,\ldots,v_r)=\sum_{\bf a=(a_1,\ldots,a_r)} c_{\bf a}
v_1^{a_1-1}\cdots v_r^{a_r-1},$$
and compute the coefficient of a given monomial $w=v_1^{b_1-1}\cdots 
v_{r+s}^{b_{r+s}-1}$ in each term of the alternility sum $A_{r,s}$.
For the shuffle-type terms in the alternility sum
$$A(v_{\sigma^{-1}(1)},\ldots, v_{\sigma^{-1}(r+s)})=
\sum_{\bf a=(a_1,\ldots,a_{r+s})} c_{\bf a}
v_{\sigma^{-1}(1)}^{a_1-1}\cdots v_{\sigma^{-1}(r+s)}^{b_{r+s}-1},$$
the coefficient of $w$ is the single coefficient 
$c_{b_{\sigma^{-1}(1)},\ldots,b_{\sigma^{-1}(r+s)}}$ of $A$. But also in the
case of the terms with denominators in the alternility sum, the coefficient
of the monomial $w$ is a single coefficient of $A$. Indeed, since $A$
is polynomial-valued, these terms simplify into polynomials whose
monomials each have one coefficient from $A$ as coefficient. We give the
example of the depth 4 term corresponding to $(a+c,b+d)$:
$${{1}\over{(v_1-v_3)(v_2-v_4)}}\Bigl(A(v_1,v_2)-A(v_3,v_2)-A(v_1,v_4)+A(v_3,v_4)\Bigr)$$
$$={{1}\over{(v_1-v_3)(v_2-v_4)}}\sum_{a,b} c_{a,b}(v_1^a-v_3^a)(v_2^b-v_4^b)$$
$$=\sum c_{a,b} (v_1^{a-1}+v_1^{a-1}v_3+\cdots +v_3^{a-1})(v_2^{b-1}+
v_2^{b-2}v_4+\cdots+v_4^{b-2});$$
thus, the coefficient of a given monomial $w=v_1^{b_1-1}\cdots v_4^{b_4-1}$ 
is equal to $c_{b_1+b_3,b_2+b_4}$. Thus, the coefficient of a single
monomial in the alternility sum $A_{r,s}$ is exactly equal to the
stuffle sum on the coefficients of the power series $f$ such that $A=mi(f)$.
\vskip .3cm
(iv) This assertion follows directly from the fact that if a polynomial
$f\in {\rm Lie}_n[x,y]$ is such that $f_Y$ satisfies the stuffle relations
in depths $1\le r<n$, then there exists a unique term in $y^n$, namely
$a_y={{-1}\over{n}}(f|x^{n-1}y)y^n$, such that $f_Y+a_y$ satisfies 
the stuffle relations in all depths $1\le r\le n$.  (Cf. [SC, Theorem 2]). 

\vskip 1cm
\noindent {\bf \S A.6. Proof of Proposition 4.2.6.}
\vskip .3cm
The proof consists in putting together a bunch of
niggly lemmas, following Ecalle's indications in [Eupolars]. Let $I$ be the
mould concentrated in depth 1 defined by $I(u_1)=1$, and $Pa$ the mould concentrated
in depth 1 defined by $Pa(u_1)=1/u_1$.
\vskip .3cm
\noindent {\bf Lemma A.6.1.} {\it We have $dupal(u_1)=I$, and for $r\ge 1$, 
$$dupal(u_1,\ldots,u_r)={{B_r}\over{r!}}lu\bigl(lu(\cdots lu(I,Pa),\cdots,Pa),Pa\bigr).
\eqno(A.6.1)$$}
\noindent {\bf Proof.} Let us use the notation $lu^r(I,Pa,\ldots,Pa)$ for
the bracket $lu(lu(\cdots lu(I,Pa),\cdots,Pa),Pa)$ where $lu$ is iterated $r$ times.
By the definition (4.2.4) of $dupal$, we certainly have $dupal(u_1)=1$. Let us
use induction on $r$.  Assume that
$${{(r-1)!}\over{B_{r-1}}}dupal(u_1,\ldots,u_{r-1})=lu^{r-2}(I,Pa,\ldots,Pa).\eqno(A.6.2)$$
We then have 
$$\eqalign{lu^{r-1}&(I,Pa,\ldots,Pa)(u_1,\ldots,u_r)\cr
&=
{{(r-1)!}\over{B_{r-1}}}\Bigl(dupal(u_1,\ldots,u_{r-1})Pa(u_r)-Pa(u_1)dupal(u_2,\ldots,u_r)\Bigr)\cr
&={{1}\over{u_1\cdots u_r}}\Bigl(\sum_{i=0}^{r-1}
(-1)^i\Bigl({{r-1}\atop{i}}\Bigr) (u_{i+1}-u_{i+2})\Bigr)\cr
&={{1}\over{u_1\cdots u_r}}\Bigl(\sum_{i=0}^{r}
(-1)^i\Bigl({{r}\atop{i}}\Bigr) u_{i+1}\Bigr)\cr
&={{r!}\over{B_r}}dupal(u_1,\ldots,u_r).}$$
This concludes the proof.\hfill{$\square$}
\vskip .3cm
Since $dapal=swap(dipil)$, it is given by
$$dapal(u_1,\cdots,u_r)=-{{1}\over{(r+1)!}}swap(re_r),\eqno(A.6.3)$$
where we see explicitly from the definition of $re_r$ in (4.1.3) that
$$swap(re_r)(u_1,\ldots,u_r)={{ru_1+(r-1)u_2+\cdots +2u_{r-1}+u_r}\over{u_1\cdots u_r
(u_1+\cdots +u_r)}}.\eqno(A.6.4)$$  
Let $mu_q(Pa)=mu(\underbrace{Pa,\ldots,Pa}_q)$.
The following lemma concerns the mould $swap(re_r)$. 
\vskip .2cm
\noindent {\bf Lemma A.6.2.} {\it For $r\ge 1$, the mould $swap(re_r)$ satisfies
\vskip .2cm
\noindent (i) $swap(re_r)+anti\cdot swap(re_r)=(r+1)\ mu_r(Pa)$\hfill{\rm (A.6.5)}
\vskip .2cm
and
\vskip .2cm
\noindent (ii) $-push\cdot swap(re_r)=anti\cdot swap(re_r).$\hfill{\rm (A.6.6)}}
\vskip .3cm
\noindent {\bf Proof.} (i) By (A.6.4), we have
$$swap(re_r)+anti\cdot swap(re_r)=(r+1){{1}\over{u_1\cdots u_r}},$$
and this is nothing other than $r+1$ times $mu_r(Pa)$.
\vskip .1cm
(ii) This is trivial; indeed the right-hand side is just
$${{u_1+2u_2+\cdots+ru_r}\over{u_1\cdots u_r(u_1+\cdots+u_r)}},\eqno(A.6.7)$$
whereas $push\cdot swap(re_r)$ is given by
$$-{{r(-u_1-\cdots -u_r)+(r-1)u_1+\cdots+2u_{r-2}+u_{r-1}}\over{u_1\cdots u_r
(u_1+\cdots+u_r)}},$$
which is nothing but the negative of (A.6.7).  \hfill{$\square$}
\vskip .3cm
We need one more lemma that will help us compute the key term $irat(dapal)\cdot
dupal$ of (4.2.9). 
\vskip .3cm
\noindent {\bf Lemma A.6.3.} {\it We have
$$irat(swap(re_r))\cdot mu_q(Pa)=-(r-q+1)mu_{r+q}(Pa)+\qquad\qquad\qquad$$
$$\qquad\qquad\qquad mu\bigl(swap(re_r),mu_q(Pa)\bigr)+mu\bigl(mu_q(Pa),anti\cdot swap(re_r)\bigr).\eqno(A.6.8)$$}\par
\noindent {\bf Proof.} Thanks to (4.2.3), we can replace $irat$ by $iwat$ in
(A.6.8), since by definition $irat(B)=iwat(B)$ whenever $B$ is a mould such that
$anti(B)=-push(B)$.  Using $iwat$ makes it easier to prove (A.6.8).  We will
do it by induction on $q$.
\vskip .3cm
\noindent {\bf Base case $q=1$.}  We first compute the mould 
$iwat(swap(re_r))\cdot Pa$, which is concentrated in depth $r+1$. By definition,
we have $iwat(swap(re_r))=amit(swap(re_r))+anit(anti(swap(re_r))).$
We check directly using (2.2.1) that
$$\eqalign{amit(swap(re_r))\cdot Pa(u_1,\ldots,u_{r+1})&=
swap(re_r)(u_1,\ldots,u_r){{1}\over{u_1+\cdots+u_{r+1}}}\cr
&= {{ru_1+\cdots+2u_{r-1}+u_r}\over{u_1\cdots u_r(u_1+\cdots+u_r)
(u_1+\cdots+u_{r+1})}}.}\eqno(A.6.9)$$
Similarly, we check directly from (2.2.2) that
$$anit(anti(swap(re_r)))\cdot Pa(u_1,\ldots,u_{r+1})=
{{ru_{r+1}+\cdots+2u_3+u_2}\over{u_2\cdots u_{r+1}(u_2+\cdots+u_{r+1})
(u_1+\cdots+u_{r+1})}}.\eqno(A.6.10)$$
Putting (A.6.10) and (A.6.11) together immediately yields 
$$iwat(swap(re_r))\cdot Pa(u_1,\ldots,u_{r+1})=$$
$${{u_1u_2+2u_1u_3+\cdots+(r-1)u_1u_r+ru_1u_{r+1}+(r-1)u_2u_{r+1}+\cdots
+2u_{r-1}u_{r+1}+u_ru_{r+1}}\over{u_1\cdots u_r(u_1+\cdots+u_{r-1})(u_2+\cdots+u_r)}}.
\eqno(A.6.11)$$
Now, the right-hand of (A.6.8) for $q=1$ is given by
$${{-r}\over{u_1\cdots u_{r+1}}}+{{ru_1+\cdots +u_r}\over{u_1\cdots u_{r+1}(u_1+\cdots
+u_r)}}+{{u_2+\cdots+ru_{r+1}}\over{u_1\cdots u_{r+1}(u_2+\cdots +u_{r+1})}},$$
and putting this over a common denominator yields exactly (A.6.11).  This settles the
base case.
\vskip .2cm
\noindent {\bf Induction step.} Assume that (A.6.8) holds up to $q$.  We compute
$$\eqalign{&irat(swap(re_r))\cdot mu_{q+1}(Pa)(u_1,\ldots,u_{r+q+1})\cr
&=mu\bigl(irat(swap(re_r))\cdot mu_q(Pa)\,,\,Pa\bigr)+
mu\bigl(mu_q(Pa)\,,\,irat(swap(re_r))\cdot Pa\bigr)\cr
&=-(r-q+1){{1}\over{u_1\cdots u_{r+q+1}}}+
mu\bigl(swap(re_r),mu_q(Pa)){{1}\over{u_{r+q+1}}}\cr
&\quad +mu(mu_q(Pa),anti\cdot swap(re_r)){{1}\over{u_{r+q+1}}}\cr
&\quad +{{1}\over{u_1\cdots u_q}}\cdot \Bigl(irat(swap(re_r))\cdot Pa\Bigr)(u_{q+1},\ldots,u_{q+r+1})\cr
&=-(r-q+1){{1}\over{u_1\cdots u_{r+q+1}}}+
{{ru_1+\cdots+u_r}\over{u_1\cdots u_{r+q+1}(u_1+\cdots +u_r)}} \cr
&\quad +{{u_{q+1}+\cdots +ru_{q+r}}\over{u_1\cdots u_{r+q+1}(u_{q+1}+\cdots +u_{q+r})}}
+{{-r}\over{u_1\cdots u_{r+q+1}}}\cr
&\quad +{{ru_{q+1}+\cdots +u_{r+q}}\over{u_1\cdots u_{r+q+1}(u_{q+1}+\cdots
+u_{r+q})}}+{{u_{q+2}+\cdots+ru_{r+q+1}}\over{u_1\cdots u_{r+q+1}(u_{q+2}+\cdots +u_{r+q+1})}}\cr
&=-(r-q){{1}\over{u_1\cdots u_{r+q+1}}}+
{{ru_1+\cdots+u_r}\over{u_1\cdots u_{r+q+1}(u_1+\cdots +u_r)}} \cr
&\quad +{{u_{q+2}+\cdots+ru_{r+q+1}}\over{u_1\cdots u_{r+q+1}(u_{q+2}+\cdots +u_{r+q+1})}}\cr
&=-(r-q)mu_{r+q+1}(Pa)+mu\bigl(swap(re_r),mu_{q+1}(Pa)\bigr)+mu\bigl(mu_{q+1}(Pa),anti\cdot swap(re_r)\Bigr),
}$$
proving the induction step. This concludes the proof of Lemma A.6.3.\hfill{$\square$}
\vskip .5cm
We will now compute the term $irat(dapal)\cdot dupal$ of (4.2.9).
We have
$$\eqalign{irat(dapal)\cdot dupal&=irat\Bigl(\sum_{r\ge 1} {{-1}\over{(r+1)!}}swap(re_r)\Bigr)\cdot
\Bigl(\sum_{s\ge 1} {{B_s}\over{s!}}lu^{s-1}(I,Pa,\ldots,Pa)\Bigr)\cr
&=\sum_{r,s\ge 1} {{-1}\over{(r+1)!}}{{B_s}\over{s!}} irat\Bigl(swap(re_r)\Bigr)\cdot lu^{s-1}(I,Pa,\ldots,Pa).}$$
Writing $lu^{s-1}(I,Pa,\ldots,Pa)=\sum_{i=0}^s (-1)^i\Bigl({{s-1}\atop{i}}\Bigr)mu\bigl(mu_i(Pa),I,mu_{s-1-i}(Pa)\bigr)$, this gives
$$\sum_{r,s\ge 1} \sum_{i=0}^{s-1} 
{{-1}\over{(r+1)!}}{{B_s}\over{s!}} (-1)^i \Bigl({{s-1}\atop{i}}\Bigr) 
irat\bigl(swap(re_r)\bigr)\cdot mu\bigl(mu_i(Pa),I,mu_{s-1-i}(Pa)\bigr).$$
Since $irat(swap(re_r))$ is a derivation for $mu$, this is equal to
\vfill\eject
$$\sum_{r,s\ge 1} \sum_{i=0}^{s-1} E_{r,s,i}
\biggl(mu\Bigl(irat\bigl(swap(re_r)\bigr)\cdot mu_i(Pa),I,mu_{s-1-i}(Pa)\Bigr)$$
$$\qquad+mu\Bigl(mu_i(Pa),I,irat\bigl(swap(re_r)\bigr)\cdot mu_{s-1-i}(Pa)\Bigr)$$
$$\qquad+mu\Bigl(mu_i(Pa),irat\bigl(swap(re_r)\bigr)\cdot I,mu_{s-1-i}(Pa)\Bigr)\biggr),$$
where $E_{r,s,i}={{-1}\over{(r+1)!}}{{B_s}\over{s!}} (-1)^i \Bigl({{s-1}\atop{i}}\Bigr)$. 
Using (A.6.8), this becomes
$$\eqalign{\sum_{r,s\ge 1} &\sum_{i=0}^{s-1} E_{r,s,i}\biggl(-(r-i+1)\,mu\bigl(mu_{r+i}(Pa),I,mu_{s-1-i}(Pa)\bigr)
\cr
&+mu\bigl(swap(re_r),mu_i(Pa),I,mu_{s-1-i}(Pa)\bigr)\cr
&+mu\bigl(mu_i(Pa),anti\cdot swap(re_r),I,mu_{s-1-i}(Pa)\bigr)\cr
&-(r-s+i+2)\,mu\bigl(mu_i(Pa),I,mu_{r+s-1-i}(Pa)\bigr)\cr
&+mu\bigl(mu_i(Pa),I,swap(re_r),mu_{s-1-i}(Pa)\bigr)\cr
&+mu\bigl(mu_i(Pa),I,mu_{s-1-i}(Pa),anti\cdot swap(re_r)\bigr)\cr
&+mu\bigl(mu_i(Pa),irat\bigl(swap(re_r)\bigr)\cdot I,mu_{s-1-i}(Pa)\bigr)\biggr).
}\eqno(A.6.12)$$
Let us use the following substitution in the two terms containing $anti\cdot swap(re_r)$:
$$anti\cdot swap(re_r)=(r+1)\,mu_r(Pa)-swap(re_r).$$
Then (A.6.12) becomes
$$\eqalign{\sum_{r,s\ge 1} &\sum_{i=0}^{s-1} E_{r,s,i}\biggl(-(r-i+1)\,mu\bigl(mu_{r+i}(Pa),I,mu_{s-1-i}(Pa)\bigr)
\cr
&+mu\bigl(swap(re_r),mu_i(Pa),I,mu_{s-1-i}(Pa)\bigr)\cr
&+(r+1)mu\bigl(mu_{r+i}(Pa),I,mu_{s-1-i}(Pa)\bigr)\cr
&-mu\bigl(mu_i(Pa),swap(re_r),I,mu_{s-1-i}(Pa)\bigr)\cr
&-(r-s+i+2)\,mu\bigl(mu_i(Pa),I,mu_{r+s-1-i}(Pa)\bigr)\cr
&+mu\bigl(mu_i(Pa),I,swap(re_r),mu_{s-1-i}(Pa)\bigr)\cr
&+(r+1)mu\bigl(mu_i(Pa),I,mu_{r+s-1-i}(Pa)\bigr)\cr
&-mu\bigl(mu_i(Pa),I,mu_{s-1-i}(Pa),swap(re_r)\bigr)\cr
&+mu\bigl(mu_i(Pa),irat\bigl(swap(re_r)\bigr)\cdot I,mu_{s-1-i}(Pa)\bigr)\biggr).
}\eqno(A.6.13)$$
Putting like terms together, this becomes
$$\eqalign{\sum_{r,s\ge 1} &\sum_{i=0}^{s-1} E_{r,s,i}\biggl(i\,mu\bigl(mu_{r+i}(Pa),I,mu_{s-1-i}(Pa)\bigr) \cr
&+(s-i-1)\,mu\bigl(mu_i(Pa),I,mu_{r+s-1-i}(Pa)\bigr)\cr
&+mu\bigl(swap(re_r),mu_i(Pa),I,mu_{s-1-i}(Pa)\bigr)\cr
&-mu\bigl(mu_i(Pa),swap(re_r),I,mu_{s-1-i}(Pa)\bigr)\cr
&+mu\bigl(mu_i(Pa),I,swap(re_r),mu_{s-1-i}(Pa)\bigr)\cr
&-mu\bigl(mu_i(Pa),I,mu_{s-1-i}(Pa),swap(re_r)\bigr)\cr
&+mu\bigl(mu_i(Pa),irat\bigl(swap(re_r)\bigr)\cdot I,mu_{s-1-i}(Pa)\bigr)\biggr).
}\eqno(A.6.14)$$
We will compare (A.6.14)=$irat(dapal)\cdot dupal$ with the other
crucial term $lu(dapal,dupal)$ from (4.2.9). We have
$$\eqalign{lu&(dapal,dupal)=mu(dapal,dupal)-mu(dupal,dapal)\cr
&=\sum_{r,s\ge 1} \Biggl({{-1}\over{(r+1)!}}{{B_{s}}\over{s!}}mu\Bigl(swap(re_r),
lu^{s-1}(I,Pa,\ldots,Pa)\Bigr)\cr
&\qquad\qquad \ -{{-1}\over{(r+1)!}}{{B_{s}}\over{s!}}mu\Bigl(lu^{s-1}(I,Pa,\ldots,Pa),
swap(re_r)\Bigr)\Biggr)\cr
&=\sum_{r,s\ge 1} \sum_{i=0}^{s-1}E_{r,s,i}\Biggl(mu\Bigl(swap(re_r),
mu_i(Pa),I,mu_{s-1-i}(Pa)\Bigr)\cr
&\qquad\qquad\qquad\qquad\ -mu\Bigl(mu_i(Pa),I,mu_{s-1-i}(Pa),swap(re_r)\Bigr) \Biggr).
}\eqno(A.6.15)$$
Let us rewrite (4.2.9) as
$$irat(dapal)\cdot dupal - lu(dapal,dupal)=der\cdot dupal-dur\cdot dapal.\eqno(A.6.16)$$ 
We note that this equality holds in depth $d=1$ since the depth 1 part of the
left-hand side is zero, and the depth one part of $der\cdot dupal$ is
equal to that of $dur\cdot dapal$, namely $-1/2$. Thus from now on we 
work in depth $d>1$.

The left-hand side is (A.6.14) - (A.6.15), which we compute as
$$\eqalign{\sum_{r,s\ge 1} &\sum_{i=0}^{s-1} E_{r,s,i}\biggl(i\,mu\bigl(mu_{r+i}(Pa),I,mu_{s-1-i}(Pa)\bigr) \cr
&+(s-i-1)\,mu\bigl(mu_i(Pa),I,mu_{r+s-1-i}(Pa)\bigr)\cr
&-mu\bigl(mu_i(Pa),swap(re_r),I,mu_{s-1-i}(Pa)\bigr)\cr
&+mu\bigl(mu_i(Pa),I,swap(re_r),mu_{s-1-i}(Pa)\bigr)\cr
&+mu\bigl(mu_i(Pa),irat\bigl(swap(re_r)\bigr)\cdot I,mu_{s-1-i}(Pa)\bigr)\biggr).
}\eqno(A.6.17)$$
\vfill\eject
Setting $d=r+s$ and $ru_r=swap(re_r)$, we rewrite the sum as
$$\eqalign{\sum_{d\ge 1}\sum_{s=1}^{d-1} &\sum_{i=0}^{s-1} E_{d-s,s,i}\biggl(i\,mu\bigl(mu_{d-s+i}(Pa),I,mu_{s-1-i}(Pa)\bigr) \cr
&+(s-i-1)\,mu\bigl(mu_i(Pa),I,mu_{d-1-i}(Pa)\bigr)\cr
&-mu\bigl(mu_i(Pa),ru_{d-s},I,mu_{s-1-i}(Pa)\bigr)\cr
&+mu\bigl(mu_i(Pa),I,ru_{d-s},mu_{s-1-i}(Pa)\bigr)\cr
&+mu\bigl(mu_i(Pa),irat\bigl(ru_{d-s}\bigr)\cdot I,mu_{s-1-i}(Pa)\bigr)\biggr),
}\eqno(A.6.18)$$
which is useful because $d$ gives the depth of the mould.
Let us consider the first two lines of (A.6.18), whose simple expressions
are easy to compute directly. For given indices $d,s,i$, we have
$$i\,mu\bigl(mu_{d-s+i}(Pa),I,mu_{s-1-i}(Pa)\bigr) 
+(s-i-1)\,mu\bigl(mu_i(Pa),I,mu_{d-1-i}(Pa)\bigr)$$
$$={{(s-i-1)u_{i+1}+iu_{d-s+i+1}}\over{u_1\cdots u_d}}.
\eqno(A.6.19)$$
The next three lines taken together are even simpler, since for given 
$d,s,i$ we have
$$-mu\bigl(mu_i(Pa),ru_{d-s},I,mu_{s-1-i}(Pa)\bigr)
+mu\bigl(mu_i(Pa),I,ru_{d-s},mu_{s-1-i}(Pa)\bigr)$$
$$+mu\bigl(mu_i(Pa),irat\bigl(ru_{d-s}\bigr)\cdot I,mu_{s-1-i}(Pa)\bigr)
={{(d-s+1)u_{i+1}}\over{u_1\cdots u_d}}.  \eqno(A.6.20)$$
Using (A.6.19) and (A.6.20) we see that
in given depth $d$, (A.6.18) is equal to
$$\sum_{s=1}^{d-1}\sum_{i=0}^{s-1} E_{d-s,s,i}{{(d-i)u_{i+1}+iu_{d-s+i+1}}\over{u_1\cdots u_d}}$$
$$={{1}\over{u_1\cdots u_d}}\sum_{s=1}^{d-1}\sum_{i=0}^{s-1} (-1)^{i+1}{{1}\over{(d-s+1)!}}{{B_s}\over{s!}}\Bigl({{s-1}\atop{i}}\Bigr)\bigl((d-i)u_{i+1}+iu_{d-s+i+1}\bigr)$$
$$={{1}\over{u_1\cdots u_d}}\sum_{s=1}^{d-1}\sum_{i=0}^{s-1} (-1)^{i+1}{{1}\over{(d-s+1)!}}{{B_s}\over{s!}}{{(s-1)!}\over{i!(s-1-i)!}}\bigl((d-i)u_{i+1}+iu_{d-s+i+1}\bigr)$$
$$={{1}\over{u_1\cdots u_d}}\sum_{s=1}^{d-1}\sum_{i=0}^{s-1} (-1)^{i+1}{{1}\over{(d-s+1)!}}{{B_s}\over{s}}{{1}\over{i!(s-1-i)!}}\bigl((d-i)u_{i+1}+iu_{d-s+i+1}\bigr)$$
\vfill\eject
The coefficient of a given $u_j$ for $j\in \{1,\ldots,d\}$ in the linear
factor is thus given by
$$\sum_{s=j}^{d-1}(-1)^j{{1}\over{(d-s+1)!}}{{B_s}\over{s}}
{{d-j+1}\over{(j-1)!(s-j)!}}+
\sum_{s=d-j+2}^{d-1} (-1)^{j+s-d}{{1}\over{(d-s+1)!}}{{B_s}\over{s}}{{1}\over{(j-d+s-2)!(d-j)!}}\eqno(A.6.21)$$
Let us compare this with the depth $d$ part of $der(dupal)-dur(dapal)$, which 
is explicitly given by
$${{1}\over{u_1\cdots u_d}} \Biggl({{B_d}\over{(d-1)!}}
\bigl(\sum_{i=0}^{d-1} (-1)^i\Bigl({{d-1}\atop{i}}\Bigr)u_{i+1}\bigr)
+{{1}\over{(d+1)!}}
\bigl(du_1+(d-1)u_2+\cdots +2u_{d-1}+u_d\bigr)\Biggr).$$
In particular the coefficient of $u_j$ in the linear factor
for $j\in \{1,\ldots,d\}$ is given by
$${{B_d}\over{(d-1)!}}(-1)^{j-1}\Bigl({{d-1}\atop{j-1}}\Bigr)+
{{d-j+1}\over{(d+1)!}}.\eqno(A.6.22)$$
Let us show that (A.6.21)=(A.6.22).  Recall from the remark after (A.6.16) that
we may assume that $d>1$. We first assume that $d$ is odd, so (A.6.22)
reduces to $(d-j+1)/(d+1)!$. The equality with (A.6.22) can be simplified 
(thanks to H. Gangl) to the equality
$$\sum_{n=0}^d  {d+1 \choose n} B_n  \Biggl( {n-1 \choose n-j}  +  (-1)^{n-1} {n-1 \choose n-d+j} \Biggr)  = (-1)^j\eqno(A.6.23)$$
for odd $d>1$ and $1\le j\le d$. 
The remarkable, elegant proof of (A.6.23) was provided to us by D. Zagier.
\vskip .5cm

$(n,k)=(-1)^k(-n+k-1,k)$ if $k$ nonnegative

$(n,k)=(-1)^(n-k)(-k-1,n-k)$ if $k <=n$

$(n,k)=0$ otherwise

Interlude: $d=3, j=1$

$n=0$: $({-1} \choose 3)-({{-1}\choose {-3}})=-(3,3)-(2,2)=-2$

$n=1$: $(0 \choose 3)+({0\choose {-3}})=0$

$n=2$: $(({1\choose 1})-({1\choose 0}))=0$

$n=3$: 0

Total 0
\vskip .5cm
Interlude: $d=5, j=2$

$n=0$: $({-1 \choose -2})-({-1\choose -3})=2$

$n=1$: $({0 \choose 3})+({0\choose -3})=0$

$n=2$: $(15)(1/6)\cdot (({1\choose 0})-({1\choose -1}))=5$

$n=3$: 0

$n=4$: $(15)(-1/30)(({3\choose 2})-({3\choose 1}))=0$

$n=5$: 0

Total 4+4=8
\vskip .5cm

\vskip .2cm
Let $f_{d,j}=\sum_{n=0}^{d} \Bigl({{d+1}\atop{n}}\Bigr)B_n\Bigl({{n-1}\atop{j}}\Bigr)$ for $0\le j\le d$. Then the desired expression follows from the
following more general equality, valid for all $d\ge 0$:
$$f_{d,j}+(-1)^{n-1}f_{d,d-j}=(-1)^j+\delta_{d,0}.\eqno(A.6.24)$$
$$F(x,y)=\sum_{d\ge j\ge 0} {{1}\over{(d+1)!}}f_{d,j}x^{d-j}y^j.$$
\noindent {\bf Claim.} We have
$$F(x,y)+F(-y,-x)=1+\sum_{i,j\ge 0} (-1)^j{{1}\over{(i+j+1)!}}x^iy^j,
\eqno(A.6.25)$$
and thus in particular 
$$f_{d,j}+(-1)^df_{d,d-j}=\delta_{d,0}+(-1)^j.$$
\noindent Proof. We have
$$\eqalign{F(x,y)&=\sum_{d=0}^\infty {{x^d}\over{(d+1)!}}\Biggl(\sum_{n=0}^d\Bigl({{d+1}\over{n}}\Bigr)B_n\bigl(1+{{y}\over{x}}\bigr)^{n-1}-
\Bigl(-{{y}\over{x}}\Bigr)\Bigl(1+{{y}\over{x}}\Bigr)^{-1}\Biggr)\cr
&=\sum_{n\ge 0,r\ge 1} {{1}\over{n!r!}}B_n(x+y)^{n-1}x^r+{{1}\over{x+y}}(1-e^{-y})\cr
&={{e^x-1}\over{e^{x+y}-1}}+{{1-e^{-y}}\over{x+y}}.}$$
This expression makes the calculation of $F(x,y)+F(-y,-x)$ trivial and
proves (A.6.24).  \hfill{$\square$}
\vfill\eject
\noindent {\bf References}
\vskip .5cm
\noindent The motivation for this work comes from the (published and
unpublished) works of Jean \'Ecalle, which can be consulted
on his web page. Many of the results and suggestions announced in
his papers were completely proved and published elsewhere.
We list here the two main published articles by \'Ecalle that served
as sources for this work, along with a number of articles by other 
authors, some of which, like this book, contain the only written proofs 
of some of \'Ecalle's statements.
\vskip .3cm
\noindent [BS] S.Baumard, L. Schneps, On the derivation representation of the fundamental Lie algebra of mixed elliptic motives, {\it Ann. Math. Qu\'ebec} {\bf 41 (1)} (2014), 43-62.
\vskip .2cm
\noindent [CS] S.~Carr, L..Schneps,  in Galois-Teichm\"uller theory and Arithmetic Geometry, H.~Nakamura, F.~Pop, L.~Schneps, A. Tamagawa, eds., Adv. Stud. Pure Math. 63, Mathematical Society of Japan, 2012, 59-89.
\vskip .2cm
\noindent [E1] J. Ecalle, The flexion structure of dimorphy: flexion units, singulators,
generators, and the enumeration of multizeta irreducibles, in {\it Asymptotics
in Dynamics, Geometry and PDEs; Generalized Borel Summation II}, O. Costin,
F. Fauvet, F. Menous, D. Sauzin, eds., Edizioni della Normale, Pisa, 2011.
\vskip .1cm
\noindent [E2] J. Ecalle, Eupolars and their bialternality grid,
{\it Acta Math. Vietnam.} {\bf 40} no. 4 (2015), 545-636.
\vskip .2cm
\noindent [F] H.~Furusho, The multiple zeta algebra and the stable derivation
algebra, Publ. RIMS Kyoto Univ. {\bf 39} (2003), 695-720.
\vskip .2cm
\noindent [FK] H. Furusho and N. Komiyama, Kashiwara-Vergne and dihedral
bigraded Lie algebras in mould theory, Ann. Fac. Sci. Toulouse Math. (6)
32 (2023), no. 4, 655-725.
\vskip .2cm
\noindent [FKRS] H.~Furusho, N.~Komiyama, E.~Raphael, L.~Schneps,
On linearised and elliptic versions of the Kashiwara-Vergne Lie algebra,
preprint 2025.
\vskip .2cm
\noindent [K] N.~Komiyama, On properties of $adari(pal)$ and $ganit_v(pic)$,
arXiv:2110.04834v2, 2021.
\vskip .2cm
\noindent [R] G.~Racinet, S\'eries g\'en\'eratrices non-commutatives de
polyz\^etas et associateurs de Drinfel'd, Ph.D. dissertation, Paris, France, 2000.
\vskip .2cm
\noindent [S1] L. Schneps, Double shuffle and Kashiwara-Vergne Lie algebras, {\it J. Algebra} {\bf 367} (2012), 54-74. 
\vskip .2cm
\noindent [S2] L.~Schneps, Elliptic double shuffle, Grothendieck-Teichm\"uller and mould theory, {\it Annales Math. Qu\'ebec} {\bf 44 (2)} (2020), 261-289.  

\vskip .2cm
\noindent [SS] A.~Salerno, L.~Schneps,  Mould theory and the double shuffle
Lie algebra structure, in {\it Periods in Quantum Field Theory and Arithmetic},
J.~Burgos Gil, K.~Ebrahimi-Fard, H.~Gangl, eds., Springer Proc. Math. Stat. 2020.

\bye